\documentclass[11pt]{article}

\usepackage{latexsym}
\usepackage{amsfonts}
\usepackage{amsmath}
\usepackage{amsthm}
\usepackage{amssymb}
\usepackage{mathrsfs}
\usepackage{dsfont}    
\usepackage{bbold}     
\usepackage[english]{babel}
\usepackage{caption}
\usepackage{epsfig}
\usepackage{float}
\usepackage{subfigure}
\usepackage{psfrag}
\usepackage{graphicx}
\usepackage{epsfig}

\DeclareMathOperator{\Val}{\matV}

\newtheorem{theorem}{Theorem}
\newtheorem*{prop*}{Theorem}
\newtheorem{theo}[theorem]{Theorem}

\newtheorem{defi}[theorem]{Definition}
\newtheorem{lemma}[theorem]{Lemma}
\newtheorem{prop}[theorem]{Proposition}
\newtheorem{rmk}[theorem]{Remark}

\newtheorem{hyp}{Hypothesis}

\newcommand{\zerarcounters}{\setcounter{equation}{0}\setcounter{theorem}{0}\setcounter{figure}{0}}

\newcommand{\ZZZ}{\mathds{Z}}

\newcommand{\NNN}{\mathds{N}}

\newcommand{\RRR}{\mathds{R}}
\newcommand{\TTT}{\mathds{T}}
\newcommand{\uno}{\mathds{1}}

\newcommand{\calA}{{\mathcal A}}
\newcommand{\BB}{{\mathcal B}}
\newcommand{\CCCC}{{\mathcal C}}

\newcommand{\calF}{{\mathcal F}}
\newcommand{\calG}{{\mathcal G}}

\newcommand{\calL}{{\mathcal L}}
\newcommand{\MM}{{\mathcal M}}

\newcommand{\calP}{{\mathcal P}}

\newcommand{\RR}{{\mathcal R}}

\newcommand{\gotd}{{\mathfrak d}}

\newcommand{\gotn}{{\mathfrak n}}

\newcommand{\gotA}{{\mathfrak A}}
\newcommand{\gotB}{{\mathfrak B}}
\newcommand{\gotC}{{\mathfrak C}}

\newcommand{\gotF}{{\mathfrak F}}

\newcommand{\gotL}{{\mathfrak L}}
\newcommand{\gotM}{{\mathfrak M}}
\newcommand{\gotN}{{\mathfrak N}}

\newcommand{\gotR}{{\mathfrak R}}
\newcommand{\gotS}{{\mathfrak S}}

\newcommand{\matV}{{\mathscr V}}

\newcommand{\ol}{\overline}
\newcommand{\oln}{\overline{n}}

\newcommand{\Fullbox}{{\rule{2.0mm}{2.0mm}}}

\newcommand{\EP}{\hfill\Fullbox\vspace{0.2cm}}
\newcommand{\prova}{\noindent{\it Proof. }}
\newcommand{\io}{\infty}
\newcommand{\e}{\varepsilon}
\newcommand{\al}{\alpha}
\newcommand{\de}{\delta}
\newcommand{\be}{\beta}
\newcommand{\ze}{\zeta}

\newcommand{\x}{\xi}
\newcommand{\p}{\pi}

\newcommand{\g}{\gamma}
\newcommand{\om}{\omega}

\newcommand{\la}{\lambda}
\newcommand{\La}{\Lambda}

\newcommand{\s}{\sigma}

\newcommand{\oo}{\boldsymbol{\omega}}

\newcommand{\aaa}{\boldsymbol{\alpha}}

\newcommand{\nn}{\boldsymbol{\nu}}
\newcommand{\pps}{\boldsymbol{\psi}}

\newcommand{\vzero}{\boldsymbol{0}}
\newcommand{\xx}{\boldsymbol{x}}
\newcommand{\yy}{\boldsymbol{y}}

\newcommand{\AAA}{\boldsymbol{A}}

\newcommand{{\resonance}}{relevant self-energy cluster }
\newcommand{\resonances}{relevant self-energy clusters }
\newcommand{\ii}{{\rm i}}
\newcommand{\To}{{\mathring{T}}}

\newcommand{\Ga}{\Gamma}
\newcommand{\De}{\Delta}

\def\ins#1#2#3{\vbox to0pt{\kern-#2 \hbox{\kern#1 #3}\vss}\nointerlineskip}

\begin{document}


\title{\bf Resonant motions in the presence of degeneracies
\\for quasi-periodically perturbed systems}

\author
{\bf Livia Corsi and Guido Gentile
\vspace{2mm}
\\ \small 
Dipartimento di Matematica, Universit\`a di Roma Tre, Roma,
I-00146, Italy
\\ \small 
E-mail: lcorsi@mat.uniroma3.it, gentile@mat.uniroma3.it}


\date{}

\maketitle

\begin{abstract}
We consider one-dimensional systems in the presence of a quasi-periodic
perturbation, in the analytical setting, and study the problem of existence
of quasi-periodic solutions which are resonant with the frequency vector of
the perturbation. We assume that the unperturbed system is locally integrable
and anisochronous, and that the frequency vector of the perturbation satisfies
the Bryuno condition. Existence of resonant solutions is related to the zeroes of a
suitable function, called the Melnikov function -- by analogy with the periodic case.
We show that, if the Melnikov function has a zero of odd order and under
some further condition on the sign of the perturbation parameter, then
there exists at least one resonant solution which continues an unperturbed
solution. If the Melnikov function is identically zero then one can push
perturbation theory up to the order where a counterpart of Melnikov function
appears and does not vanish identically: if such a function has a zero of odd order
and a suitable positiveness condition is met, again the same persistence result is
obtained. If the system is Hamiltonian, then the procedure can be
indefinitely iterated and no positiveness condition must be required:
as a byproduct, the result follows that at least one resonant quasi-periodic solution
always exists with no assumption on the perturbation.
Such a solution can be interpreted as a (parabolic) lower-dimensional torus.
\end{abstract}


\zerarcounters
\section{Introduction}
\label{sec:1}

Melnikov theory studies the fate of homoclinic and periodic orbits of
two-dimensional dynamical systems when they are periodically perturbed;
see for instance \cite{GH} for an introduction to the subject.
The problem can be stated as follows. Consider in $\RRR^{2}$ a dynamical system of the form
\begin{equation} \label{eq:1.1}
\begin{cases}
\dot{x} = f_{1}(x,y) +\e g_{1}(x,y,t),\\
\dot{y} = f_{2}(x,y) + \e g_{2}(x,y,t),
\end{cases}
\end{equation}
with $f_{1},f_{2},g_{1},g_{2}$ `sufficiently smooth', $g_{1}$ and
$g_{2}$ $T$-periodic in $t$ for some $T>0$ and $\e$ a small parameter,
called the \emph{perturbation parameter}. 
If $f_{1}=\partial_{y}h$ and $f_{2}=-\partial_{x}h$, for a suitable
function $h$, the unperturbed system is Hamiltonian.
Assume that for $\e=0$ the system
(\ref{eq:1.1}) admits a homoclinic orbit $u_{1}(t)$ to a hyperbolic
saddle point $p$ and that the bounded region of the phase space delimited by
$\{u_{1}(t):t\in\RRR\}\cup\{p\}$ is filled with a continuous family of
periodic orbits $u_{\de}(t)$, $\de\in(0,1)$,
whose periods tend monotonically to $\io$ as $\de\to1$.
Because of the assumptions, it is easy to see (as an application of the implicit function theorem)
that for $\e\ne0$ small enough, the system (\ref{eq:1.1}) admits
a hyperbolic periodic orbit $\widetilde{u}(t,\e)=p+O(\e)$. Then one can ask
whether the stable and unstable manifolds of $\widetilde{u}(t,\e)$ intersect
transversely (in turn if this happens it can be used to prove that chaotic motions occur).
Another natural question is what happens to the periodic orbits
$u_{\de}(t)$ when $\e\ne0$. In particular one can investigate under
which conditions `periodic orbits persist', that is there are periodic orbits
which are close to the unperturbed ones and reduce to them
when the perturbation parameter is set equal to zero.
If such orbits exist, they are called \emph{subharmonic} or \emph{resonant} orbits.

Both the existence of transverse intersection of the stable and unstable manifolds
of $\widetilde{u}(t,\e)$ and the persistence of periodic orbits
are related to the zeroes of suitable functions.
More precisely if one define the \emph{Melnikov function} as
\begin{equation} \label{eq:1.2}
M(t_{0}):=\int_{-\io}^{\io}{\rm d}t\,\bigl(
f_{2}(u_{1}(t-t_{0}))g_{2}(u_{1}(t-t_{0}),t) - f_{1}(u_{1}(t-t_{0}))g_{1}(u_{1}(t-t_{0}),t) \bigr),
\end{equation}
then if $M(t_{0})$ has simple zeroes the stable and unstable manifolds of $\widetilde{u}(t,\e)$ 
intersect transversely, while if $M(t_{0})\ne0$ for all $t_{0}\in\RRR$ no intersection occurs;
essentially $M(t_{0})$ measures the distance between the two
manifolds along the normal to the homoclinic orbit at $u_{1}(t_{0})$.
Concerning the periodic orbits, if the period $T_{\de}$ of $u_{\de}(t)$
is not commensurable with the period $T$ of the functions $g_{1},g_{2}$, in general
such an orbit will not persist under perturbations. Otherwise, set $T_{\de}=mT/n$ and
define the \emph{(subharmonic) Melnikov function} as
\begin{equation} \label{eq:1.3}
M_{m/n}(t_{0}):=\int_{0}^{mT}{\rm d}t\,\bigl(
f_{2}(u_{\de}(t-t_{0}))g_{2}(u_{\de}(t-t_{0}),t)-
f_{1}(u_{\de}(t-t_{0}))g_{1}(u_{\de}(t-t_{0}),t)
\bigr).
\end{equation}
If $M_{m/n}(t_{0})$ has a simple zero then \eqref{eq:1.1} admits
a subharmonic orbit $\ol{u}(t,\e)$ with period $mT$; in particular, if
the functions $f_{1},f_{2},g_{1},g_{2}$ are analytic, then $\ol{u}(t,\e)$ is analytic in both
$\e$ and $t$. If there are no zeroes at all, no periodic solution persists.
The proof of the claims above is rather standard and it is essentially
based on the application of the implicit function theorem.
A possible approach for the case of subharmonic orbits
consists in splitting the equations of motion
into two separate sets of equations, the so-called \emph{range equations}
and \emph{bifurcation equations}: one can solve the range equations
in terms of the free parameter $t_{0}$ and then fix the latter by solving
the bifurcation equations, which represent an implicit function problem.

The assumption that the zeroes of the Melnikov function are
simple corresponds to a (generic) non-degeneracy condition on the perturbation.
When the zeroes are not simple, the situation is slightly more complicated.
In the case of subharmonic orbits, the same result of persistence extends
to the more general case of zeroes of odd order \cite{ACE}, and interesting new
analytical features of the solutions appear \cite{ALGM,GBD,CG1};
indeed the subharmonic solutions  turn out to be analytic in a suitable
fractional power of $\e$ rather than $\e$ itself. On the other hand
if the zeroes are of even order one cannot predict a priori the persistence
of periodic orbits. Finally, if the Melnikov function is identically zero,
one has to consider higher order generalisations of it and study
the existence and multiplicity of their zeroes to deal with the problem \cite{CG1}.

If one considers a quasi-periodic perturbation instead of a periodic one,
that is $g_{k}(x,y,t)=G_{k}(x,y,\oo t)$, with $G_{k}:\RRR^{2}\times\TTT^{d}
\to\RRR^{2}$ and $\oo\in\RRR^{d}$, $d\ge2$, one can still ask whether
there exist hyperbolic sets run by quasi-periodic solutions with
stable and unstable manifolds which intersect transversely
and one can still study the existence of quasi-periodic solutions
which are ``resonant'' with the \emph{frequency vector}
$\oo$ of the perturbation; see below -- after (\ref{eq:1.4}) -- for
a formal definition of resonant solution for quasi-periodic forcing.

Also in the quasi-periodic case, non-degeneracy assumptions are essential
to prove transversality of homoclinic intersections.
Existence of a quasi-periodic hyperbolic orbit close to the
unperturbed saddle point and of its stable and unstable manifolds follows from general
arguments, such as the invariant manifold theorem \cite{Fe,HPS},
without even assuming any condition on the frequency vector $\oo$.
Palmer generalises Melnikov's method to the case of bounded
perturbations \cite{Pa} using the theory of exponential dichotomies \cite{C}.
A suitable generalisation of the Melnikov function for quasi-periodic forcing
is also introduced by Wiggins \cite{W1}. He shows that
if such a function has a simple zero then the stable and unstable manifolds intersect transversely.
Then, generalising the Smale-Birkhoff homoclinic theorem to the case of orbits homoclinic
to normally hyperbolic tori, he finds that there is an invariant set on
which the dynamics of a suitable Poincar\'e map is conjugate to a
subshift of finite type; in turn this yields the existence of chaos.
Similar results hold also for more general almost periodic perturbations (which include
the quasi-periodic ones as a special case): Meyer and Sell show that
also in that case the dynamics near transverse homoclinic orbits
behaves as a subshift of finite type \cite{MS} and
Scheurle, relying on Palmer's results, finds particular solutions
which have a random structure \cite{Sch};
again, to obtain transversality the Melnikov function is assumed to have simple zeroes.
Such assumption can be weakened to an assumption of ``topological non-degeneracy''
(i.e. the existence of an isolated minimum or maximum of the primitive of the
Melnikov function) as in the case of subharmonic orbits, and one can deal with the problem
by use of a variational approach; see for instance \cite{CES,S,AB,BB2}.

A natural application for the study of homoclinic intersections,
widely studied in the literature, is the quasi-periodically
forced Duffing equation \cite{W2,MS,Y}. Often, especially in applications, the frequency vector
is taken to be two-dimensional, with the two components which are nearly resonant with the proper
frequency of the unperturbed system (see for instance \cite{BH,Y}  and references therein).
Then a different approach with respect to \cite{W1} is proposed by Yagasaki \cite{Y}:
first, through a suitable change of coordinates, one arrives at a system with two frequencies,
one fast and one slow, and then one uses averaging to reduce the analysis of the original system
to that of a perturbation of a periodically forced system for which the standard Melnikov's method applies:
the persistence of hyperbolic periodic orbits and their stable and unstable manifolds
for the original system is then obtained as a consequence of the invariant manifold theorem. 
Transversality of homoclic intersections plays also a crucial role in the phenomenon of
Arnold diffusion \cite{A,LL}: non-degeneracy assumptions on the perturbation are heavily used
in the proofs existing in the literature (see e.g. \cite{DGJS,GGM2,DLS}) in order to find
lower bounds on the transversality, which in turn are fundamental to compute
the diffusion times along the heteroclinic chains (see e.g. \cite{Be,Ga2,BB,BBB,T}).
A physically relevant case, studied within the context of Arnold diffusion,
is that with frequency vectors with two fast components \cite{DGJS,GGM1,St} or
with one component much faster and one component much slower than the proper
frequency (`three scale system') \cite{GGM2,GGM3,Pr}. In such cases the
homoclinic splitting is exponentially small in the perturbation
parameter and this makes the analysis rather delicate, as one has to check
that the first order contribution to the splitting (the Melnikov function)
really dominates; in particular non-degeneracy conditions on the perturbation are needed once more.

The problem of existence of quasi-periodic orbits close to the center of the unperturbed system
is harder and does not follow from the invariant manifold theorem.
 Second-order approximations for the quasi-periodic solutions close to
the centers of a forced oscillator are studied in \cite{BH},
using the multiple scale technique for asymptotic expansions \cite{KC,NM}.
But if one wants to really prove the existence of the solution, one must
require additional assumption on $\oo$ to deal with the presence of \emph{small divisors}.
In \cite{M1}, Moser considers Duffing's equation with a quasi-periodic
driving term and assumes that (i) the system is \emph{reversible},
i.e. it can be written in the form $\dot x=f(x)$, with $f\!\!:\RRR^{n}\to\RRR^{n}$,
and there exists an involution $I\!\!:\RRR^{n}\to\RRR^{n}$ such that $f(Ix)=-If(x)$
(so that with $x(t)$ also $Ix(-t)$ is a solution), and (ii) the frequency
vector of the driving satisfies some Diophantine condition involving
also the proper frequency of the unperturbed system linearised around its center.
Then he shows that there exists a quasi-periodic solution,
with the same frequency vector as the driving, to a slightly modified
equation, in which  the coefficient of the linear term is suitably corrected.
If one tried to remove the correction then one should deal with an implicit function problem
(see \cite{BG} for a similar situation), which, without assuming any non-degeneracy condition
on the perturbation, would have the same kind of problems as in the present paper.
Quasi-periodically forced Hamiltonian oscillators are also considered in \cite{BHJVW},
where the persistence of quasi-periodic solutions close to the centers
of the unperturbed system is studied, including the case of resonance
between the frequency vector of the forcing and the proper frequency. 
However, again, non-degeneracy conditions are assumed. 

On the contrary the problem of persistence of quasi-periodic solutions
far from the stationary points, corresponding to the subharmonic solutions of the periodic case,
does not seem to have been studied a great deal.
(We can mention a paper by Xu and Jing \cite{XJ}, who consider Duffing's equation
with a two-frequency quasi-periodic perturbation and follow the approach in \cite{Y} to
reduce the analysis to a one-dimensional backbone system;
however the argument used to show the persistence
of the two-dimensional tori is incomplete and requires further hypotheses.)
Again the existence of resonant solutions is related to the zeroes of a suitable
function, still called Melnikov function by analogy with the periodic case.
If the zeroes are simple, assuming some Diophantine condition on $\oo$,
the analysis can be carried out so as to reach conclusions similar to the periodic case,
that is the persistence of resonant solutions.
In this paper we study the same problem in the case of zeroes of odd
order and additionally investigate what can still be said when the Melnikov function is identically zero. 
As remarked before, considering non-simple zeroes means removing
non-degeneracy -- and hence genericity -- conditions on the perturbation.
This introduces nontrivial technical complications,
because one is no longer allowed to separate the small divisor problem
plaguing the range equations from the implicit function problem
represented by the bifurcation equations. The method we use is based on
the analysis and resummation of the perturbation series through
renormalisation group techniques \cite{GGG,G1,G2,GG1};
for other renormalisation group approaches to small divisors problems
in dynamical systems see for instance \cite{BGK,K,L,KLM,KK}. As in \cite{CG2},
the frequency vector of the perturbation will be assumed to satisfy
the \emph{Bryuno condition}; such a condition, originally introduced
by Bryuno \cite{Br}, has been studied recently in several small divisor
problems arising in dynamical systems \cite{G1,G2,G3,KK,L,Po}.
With respect to \cite{CG2}, we consider here also non-Hamiltonian systems:
what is required on the unperturbed system is a non-degeneracy condition on the
frequency map of the periodic solutions (\emph{anisochrony condition}).
In the Hamiltonian case, such a condition becomes a convexity condition 
on the unperturbed Hamiltonian function, analogously to Cheng's paper
\cite{Ch}, where the fate of resonant tori is studied. In the
Hamiltonian case, the main difference with respect to \cite{Ch} -- and
what prevents us from simply relying on that result --  is that
we consider isochronous perturbations (while in \cite{Ch} the unperturbed
Hamiltonian is convex in all action variables) and assume a weaker
Diophantine condition on the frequency vector of the perturbation
(the Bryuno condition instead of the standard one). Furthermore, 
as we said, our method covers also the non-Hamiltonian case,
where Cheng's approach, based on a sequence of
canonical transformations \textit{\`a la} KAM, does not apply.
In the Hamiltonian case, we do not require any further assumption
on the perturbation (besides analyticity), as in \cite{Ch}.
In the non-Hamiltonian case we shall make some further assumptions.
More precisely we shall require that some zeroes of odd order appear at
some level of perturbation theory and a suitable positiveness condition holds;
see Section \ref{sec:1} -- in particular Hypotheses \ref{hyp3} and \ref{hyp4}
-- for a more formal statement.

Of course, one could also investigate what happens if the non-degeneracy
condition on the unperturbed system is completely removed too. However, this
would be a somewhat different problem and very likely a non-degeneracy
condition could become necessary for the perturbation. Not even in the KAM
theory for maximal tori, the fully degenerate case (no assumption on the
unperturbed integrable system and no assumption on the perturbation, besides
analyticity) has ever been treated in the literature -- as far as we know.

The paper is organised as follows. We consider systems of the form
(\ref{eq:1.1}) and assume that for $\e=0$ there is a family of periodic solutions
satisfying the same hypotheses as in the case of periodic forcing. In fact
the analysis we will be interested in will be essentially local, so we can allow
a more general setting and assume that, in suitable coordinates $(\be,B)\in\TTT\times\gotB$,
with $\gotB$ an open subset of $\RRR$, the unperturbed system reads
\begin{equation} \label{eq:1.4}
\begin{cases} \dot\be=\om_{0}(B) , & \\ \dot B = 0 , \end{cases}
\end{equation}
with $\om_{0}$ analytic and $\partial_{B}\om_{0}(B) \neq0$ (anisochrony condition).
As a particular case we can consider that $(B,\be)$ are canonical coordinates
(action-angle coordinates), but the formulation we are giving here is more general
and applies also to non-Hamiltonian unperturbed systems; see also \cite{BDG,GBD}.
Then we add to the vector field a small analytic quasi-periodic forcing term
with frequency vector $\oo=(\om_{1},\ldots,\om_{d})$ which satisfy
some weak Diophantine condition (Bryuno condition)
and concentrate on a periodic solution of the unperturbed system
which is resonant with $\oo$, that is a solution with $B=\ol{B}_{0}$
such that $\om(\ol{B}_{0})\nu_{0}+\om_{1}\nu_{1}+\ldots+\om_{d}\nu_{d}=0$ 
for suitable integers $\nu_{0},\nu_{1},\ldots,\nu_{d}$.
In Section \ref{sec:2} we state formally
our two main results on the persistence of such a solution:
Theorem \ref{thm:2.2} takes into account the case in which the 
system is not assumed to be  Hamiltonian and a zero of odd
order appears at some order of perturbation theory, while
Theorem \ref{thm:2.3} deals with the case in which the system is Hamiltonian
and no further assumption is made on the perturbation.
In Sections \ref{sec:3}--\ref{sec:5} we shall prove Theorem \ref{thm:2.2}.
As we shall see, the quasi-periodic solution will be only continuous in
the perturbation parameter. In fact, in contrast to the case of periodic
perturbations, in general the quasi-periodic solution is not expected to be analytic
in $\e$ nor in some fractional power of $\e$; already in the non-degenerate
anisochronous Hamiltonian case the solution has been proved only to be $C^{\io}$ smooth
in $\e$ \cite{GG1} and analyticity is very unlikely. 
In Section \ref{sec:6} we shall prove Theorem \ref{thm:2.3}:
we shall see that either (a) one is able to reduce the analysis to Theorem \ref{thm:2.2}
or (b) suitable ``cancellations'' occur to all orders in the perturbation series formally
defining the solution. In particular we shall see how the Hamiltonian structure of the
equations of motion is fundamental in order to prove such cancellations.
In turn this will imply, in case (b), the convergence of the perturbation series
and hence the existence of a solution which is analytic in the perturbation parameter:
we stress since now that this is a highly non-generic -- and hence very unlikely -- possibility.
The cancellation mechanism turns out to be quite similar to the one performed in \cite{CGP},
where Moser's modifying terms theorem \cite{M2} is proved by using Cartesian
coordinates instead of action-angle coordinates.
It would be interesting to understand the deep reason of such a similarity.

\zerarcounters
\section{Results}
\label{sec:2}

Let us consider the system
\begin{equation}\label{eq:2.1}
\left\{\begin{aligned}
&\dot{\be}=\om_{0}(B)+\e F(\oo t,\be,B),\\
&\dot{B}=\e G(\oo t,\be,B),
\end{aligned}\right.
\end{equation}
where $(\be,B)\in \TTT\times \gotB$, with $\gotB$ an open subset of $\RRR$,
$F,G\!:\TTT^{d+1}\times\gotB\to\RRR$ and $\om_{0}\!:\gotB\to\RRR$ are
real analytic functions, $\oo \in \RRR^{d}$ and $\e$ is a real parameter
called the \emph{perturbation parameter}; hence the \emph{perturbation}
$(F,G)$ is quasi-periodic in $t$ with \emph{frequency vector} $\oo$.
Without loss of generality we can assume that $\oo$ has rationally independent components.

Denote by $\cdot$ the standard scalar product in $\RRR^{d}$, i.e.
$\xx\cdot\yy=x_1y_1+\ldots+x_dy_d$ for $\xx,\yy\in\RRR^d$,
and set $|\xx|:=\|\xx\|_{1}=|x_{1}|+\ldots+|x_{d}|$. If $f\!:\RRR^{n}\to\RRR$,
$n\ge 1$, is a differentiable function, we shall denote (when no ambiguity
arises) by $\partial_jf$ the derivative of $f$ with respect to the $j$-th argument,
i.e. $\partial_{j}f(x_1,\ldots,x_n)=\partial_{x_{j}}f(x_1,\ldots,x_n)$;
if $n=1$ we shall write also $f'(x_1)=\partial_1f(x_1)=\partial f(x_1)$. Finally, for any
finite set $S$ we denote by $|S|$ its cardinality.

Take the solution for the  unperturbed system given by
$(\be(t),B(t))=(\be_{0}+\om_{0}(\ol{B}_{0})t,\ol{B}_{0})$, with $\ol{B}_{0}$ such that
$\om_{0}(\ol{B}_{0})$ is \emph{resonant} with $\oo$, i.e. such that there exists
$(\ol{\nu}_{0},\ol{\nn})\in\ZZZ^{d+1}$ for which $\om_{0}(\ol{B}_{0})\,\ol{\nu}_{0}+\oo\cdot\ol{\nn}=0$.
We want to study whether for some value of $\be_{0}$, that is for
a suitable choice of the \emph{initial phase}, such a solution can be continued under perturbation.

The resonance condition between $\om_{0}(\ol{B}_{0})$ and $\oo$ yields a ``simple resonance'' 
(or resonance of order 1) for the vector $(\om_{0}(\ol{B}_{0}),\oo)$. 
The main assumptions on (\ref{eq:2.1}) are a Diophantine condition
on the frequency vector of the perturbation and a non-degeneracy condition
on the unperturbed system. More precisely we shall require
that the vector $(\om_{0}(\ol{B}_{0}),\oo)$ satisfies the condition
\begin{equation}
\sum_{n\ge0}\frac{1}{2^{n}}\log \Big( \!\!\!\!\!\!\!\!\!\! \inf_{
\substack{(\nu_{0},\nn)\in\ZZZ^{d+1} \\ (\nu_{0},\nn) \nparallel
(\ol{\nu}_{0},\ol{\nn}) , 0<|(\nu_{0},\nn)| \le 2^{n}}}
\!\!\!\!\!\!\!\!\!\! \!\!\!\!\!\! |\om_{0}(\ol{B}_0 )\nu_{0} +
\oo\cdot\nn|\Big)^{-1}   < \infty
\label{eq:2.0} \end{equation}
and that $\om'_{0}(\ol{B}_{0}) \neq 0$. Note that the condition (\ref{eq:2.0}) is
weaker than requiring that the vector $(\om_{0}(\ol{B}_{0}),\oo)$
satisfies the standard Diophantine condition
$|\om_{0}(\ol{B}_0 )\nu_{0} +\oo\cdot\nn|> \g \, (|\nu_{0}|+|\nn|)^{-\tau}$
for suitable positive constants $\g,\tau$ and all $(\nu_{0},\nn)$
non-parallel to $(\ol{\nu}_{0},\ol{\nn})$.

Up to a linear change of coordinates, we can (and shall) assume
$\om_{0}(\ol{B}_{0})=0$, so that the vector $\ol{\nn}$,
such that $\om_{0}(\ol{B}_{0})\,\ol{\nu}_{0}+\oo\cdot\ol{\nn}=0$, must be the null vector.
Therefore it is not restrictive to formulate the assumptions on $\ol{B}_{0}$ and $\oo$ as follows.

\begin{hyp}\label{hyp1}
$\om(\ol{B}_{0})=0$ and $\oo$ satisfies the Bryuno condition $\BB(\oo)<\io$, where
\begin{equation}\label{eq:2.2}
\BB(\oo)=\sum_{n\ge0}\frac{1}{2^{n}}\log\frac{1}{\al_{n}(\oo)},
\qquad
\al_{n}(\oo)=\inf_{\substack{\nn\in\ZZZ^{d} \\ 0<|\nn| \le 2^{n}}}
|\oo\cdot\nn| .
\end{equation}
\end{hyp}

\begin{hyp}\label{hyp2}
$\om_{0}'(\ol{B}_{0}) \ne 0$.
\end{hyp}

Let us write
\begin{equation}\label{eq:2.3}
F(\pps,\be,B)=\sum_{\nn\in\ZZZ^{d}} {\rm e}^{\ii\nn\cdot\pps}F_{\nn}(\be,B),
\qquad
G(\pps,\be,B)=\sum_{\nn\in\ZZZ^{d}} {\rm e}^{\ii\nn\cdot\pps}G_{\nn}(\be,B),
\end{equation}
and note that, as $F$ and $G$ are real-valued functions, one has
\begin{equation}\label{eq:2.4}
F_{-\nn}(\be,B)=F_{\nn}(\be,B)^{*} ,\qquad
G_{-\nn}(\be,B)=G_{\nn}(\be,B)^{*} .
\end{equation}
Here and henceforth $*$ denotes complex conjugation. By analogy with the periodic case,
$G_{\vzero}(\be,\ol{B}_0)$ will be called the (first order) \textit{Melnikov function}.

\begin{hyp}\label{hyp3}
$\ol{\be}_{0}$ is a zero of order $\gotn$ for $G_{\vzero}(\be_{0},\ol{B}_{0})$,
with $\gotn$ odd, and $\e \, \om_{0}'(\ol{B}_{0}) \, \partial_{\be_{0}}^{\gotn}G_{\vzero}(\ol{\be}_{0},
\ol{B}_{0})>0$.
\end{hyp}

We look for a quasi-periodic solution to (\ref{eq:2.1}) with frequency
vector $\oo$, that is a solution of the form
$(\be(t),B(t))=(\be_{0}+b(t),B_{0}+\widetilde{B}(t))$, with
\begin{equation}\label{eq:2.5}
b(t)=\sum_{\nn\in\ZZZ^{d}_{*}}{\rm e}^{\ii\nn\cdot\oo t}b_{\nn},
\qquad
\widetilde{B}(t)=\sum_{\nn\in\ZZZ^{d}_{*}}{\rm e}^{\ii\nn\cdot\oo t}B_{\nn},
\end{equation}
where $\ZZZ^{d}_{*}=\ZZZ^{d}\setminus\{\vzero\}$.
Note that the existence of a quasi-periodic solution with frequency
$\oo$ in the variables in which $\om_{0}(\ol{B}_{0})=0$ implies
the existence of a quasi-periodic solution with frequency resonant
with $\oo$ in terms of the original variables (that is, before
performing the change of variables leading to $\om_{0}(\ol{B}_{0})=0$).

If we set $\Phi(t):=\om_{0}(B(t))+\e F(\oo t,\be(t),B(t))$ and $\Ga(t)=
\e G(\oo t,\be(t),B(t))$ and write
\begin{equation}\label{eq:2.6}
\Phi(t)=\sum_{\nn\in\ZZZ^{d}} {\rm e}^{\ii\nn\cdot\oo t}\Phi_{\nn},
\qquad
\Ga(t)=\sum_{\nn\in\ZZZ^{d}} {\rm e}^{\ii\nn\cdot\oo t}\Ga_{\nn},
\end{equation}
in Fourier space (\ref{eq:2.1}) becomes
\begin{subequations}
\begin{align}
&(i\oo\cdot\nn)b_{\nn}=\Phi_{\nn},\quad \nn\ne\vzero,
\label{eq:2.7a}\\
&(i\oo\cdot\nn)B_{\nn}=\Ga_{\nn},\quad \nn\ne\vzero,
\label{eq:2.7b}\\
&\Phi_{\vzero}=0,
\label{eq:2.7c}\\
&\Ga_{\vzero}=0
\label{eq:2.7d}.
\end{align}
\label{eq:2.7}
\end{subequations}
\vskip-.3truecm
According to the usual terminology, we shall call (\ref{eq:2.7a}) and (\ref{eq:2.7b})
the \emph{range equations}, while (\ref{eq:2.7c}) and (\ref{eq:2.7d})
will be referred to as the \emph{bifurcation equations}.

Our first result will be the following.

\begin{theo}\label{thm:2.1}
Consider the system (\ref{eq:2.1}) and assume Hypotheses \ref{hyp1},
\ref{hyp2} and \ref{hyp3} to be satisfied. Then for $\e$ small enough
there exists at least one quasi-periodic solution $(\be(t),B(t))$ with frequency vector
$\oo$ such that $(\be(t),B(t))\to
(\ol{\be}_{0},\ol{B}_{0})$ for $\e\to0$.
\end{theo}

Actually we shall prove a more general result which can be stated as follows.
We look for a formal solution $(\be(t),B(t))$, with 
\begin{equation}\label{eq:2.8}
\begin{aligned}
\be(t)=\be(t;\e,\be_{0})&=\be_{0}+\sum_{k\ge 1}\e^{k}b^{(k)}(t;\be_{0})= 
\be_{0}+\sum_{k\ge 1}\e^{k}\sum_{\nn\in\ZZZ^{d}_{*}} 
{\rm e}^{\ii\nn\cdot\oo t}b^{(k)}_{\nn}(\be_{0}), \\
B(t)=B(t;\e,\be_{0})&=\ol{B}_{0}+\sum_{k\ge 1}\e^{k}B^{(k)}(t;\be_{0})= 
\ol{B}_{0}+ \sum_{k\ge 1}\e^{k}\sum_{\nn\in\ZZZ^{d}} 
{\rm e}^{\ii\nn\cdot\oo t}B^{(k)}_{\nn}(\be_{0}) ,
\end{aligned}
\end{equation} 
and set $U(t):=\om_{0}(B(t))-\om_{0}'(\ol{B}_{0})(B(t)-\ol{B}_{0})$ and
$\phi(t)=U(t)+\e F(\oo t,\be(t),B(t))$.
Then define recursively for $k\ge1$
\begin{equation}\label{eq:2.9} 
\begin{aligned}
b_{\nn}^{(k)}(\be_{0})&=\frac{1}{(\ii\oo\cdot\nn)} 
\phi^{(k)}_{\nn} (\be_{0}) + \frac{\om_{0}'(\ol{B}_{0})}{(\ii\oo\cdot\nn)^{2}}
\Ga^{(k)}_{\nn} (\be_{0}) ,
\qquad \nn\ne\vzero \\
B_{\nn}^{(k)}(\be_{0})&=\frac{1}{(\ii\oo\cdot\nn)} 
\Ga^{(k)}_{\nn} (\be_{0}) ,\qquad \nn\ne\vzero \\
B_{\vzero}^{(k)}(\be_{0})&=-\frac{1}{\om_{0}'(\ol{B}_{0})}\phi^{(k)}_{\vzero}(\be_{0}),
\end{aligned}
\end{equation} 
where
$\Ga^{(k)}_{\nn}(\be_{0})=[G(\oo t,\be(t),B(t))]^{(k-1)}_{\nn}$ and
$\phi^{(k)}_{\nn}(\be_{0})=[U(t)]^{(k)}_{\nn}+[F(\oo t,\be(t),B(t))]^{(k-1)}_{\nn}$,
with $U^{(1)}_{\nn}(\be_{0})=0$, so that
$\Ga^{(1)}_{\nn}(\be_{0})=G_{\nn}(\be_{0},\ol{B}_{0})$
and $\phi^{(1)}_{\nn}(\be_0)=F_{\nn}(\be_{0},\ol{B}_{0})$,
while, for $k\ge2$,
\begin{equation}\label{eq:2.10}
[U(t)]^{(k)}_{\nn}=\sum_{s\ge2}\frac{1}{s!}\partial_{B}^{s}\om_{0}(\ol{B}_{0})
\sum_{\substack{\nn_{1}+\ldots+\nn_{s}=\nn \\ 
\nn_{i} \in \ZZZ^{d},\, i=1,\ldots,s}}
\sum_{\substack{k_{1}+\ldots+k_{s}=k, \\ k_{i}\ge 1}} 
\prod_{i=1}^{s} B_{\nn_{i}}^{(k_{i})}(\be_{0}) ,
\end{equation}
and
\begin{equation}\label{eq:2.11}
\begin{aligned}
&[P(\oo t,\be(t),B(t))]^{(k-1)}_{\nn}=\sum_{s\ge 1} \sum_{p+q=s}
\sum_{\substack{\nn_{0}+\ldots+\nn_{s}=\nn \\ \nn_{0},\nn_{j}\in
\ZZZ^{d}\,j=p+1,\ldots,s \\ 
\nn_{i} \in \ZZZ^{d}_{*},\, i=1,\ldots,p}} 
\frac{1}{p!q!}\partial^{p}_{\be}\partial^{q}_{B} 
P_{\nn_{0}}(\be_{0},\ol{B}_{0}) \; \times \\
&\qquad\qquad\times \!\!\!\!\!\!\!
\sum_{\substack{k_{1}+\ldots+k_{s}=k-1, \\ k_{i}\ge 1}} 
\prod_{i=1}^{p} b_{\nn_{i}}^{(k_{i})}(\be_{0}) 
\prod_{i=p+1}^{s} B_{\nn_{i}}^{(k_{i})}(\be_{0}) , 
\qquad \qquad P = F,G .
\end{aligned}
\end{equation} 

The series (\ref{eq:2.8}), with the coefficients defined as above and arbitrary $\be_0$,
turn out to be a formal solution of (\ref{eq:2.7a})-(\ref{eq:2.7c}):
the coefficients $b^{(k)}_{\nn}(\be_{0})$, $B_{\vzero}^{(k)}(\be_{0})$ and 
$B^{(k)}_{\nn}(\be_{0})$ are well defined for all 
$k\ge1$ and all $\nn\in\ZZZ^{d}_{*}$, by Hypothesis \ref{hyp1},
and solve (\ref{eq:2.7a})-(\ref{eq:2.7c}) order by order -- as it is
straightforward to check (for instance by using the formalism
introduced below in Section \ref{sec:3}) --; moreover the functions
$b^{(k)}(t;\be_{0})$ and $B^{(k)}(t;\be_{0})$ are analytic and
quasi-periodic in $t$ with frequency vector $\oo$.
 
Assume that there exists $k_{0}\in\NNN$ such that all functions 
${\Ga}^{(k)}_{\vzero}(\be_{0})$ are identically zero for
$0\le k \le k_{0}-1$; then we can solve the
equation of motion up to order $k_{0}-1$ without fixing the parameter $\be_{0}$
and moreover
$\Ga^{(k_{0})}_{\vzero}$ is a well-defined function of $\be_{0}$.

\begin{hyp}\label{hyp4} 
There exist $k_0\in\NNN$ and $\ol{\be}_{0}$ such that
${\Ga}^{(k)}_{\vzero}(\be_{0})$ vanish identically for $k<k_{0}$ and
$\ol{\be}_{0}$ is a zero of order $\bar{\gotn}$ for 
${\Ga}^{(k_{0})}_{\vzero}(\be_{0})$, with $\bar{\gotn}$ odd. Moreover one has
$\e^{k_{0}}\om_{0}'(\ol{B}_{0})
\partial^{\bar{\gotn}}_{\be_{0}}{\Ga}^{(k_{0})}_{\vzero} (\ol{\be}_{0})>0$. 
\end{hyp} 

We shall prove the following result.

\begin{theo}\label{thm:2.2}
Consider the system (\ref{eq:2.1}) and assume Hypotheses \ref{hyp1},
\ref{hyp2} and \ref{hyp4} to be satisfied. Then for $\e$ small enough
there exists at least one quasi-periodic solution $(\be(t),B(t))$ with frequency vector
$\oo$ such that $(\be(t),B(t))\to(\ol{\be}_{0},
\ol{B}_{0})$ for $\e\to0$.
\end{theo}

Note that Hypothesis \ref{hyp4} reduces to Hypothesis \ref{hyp3} if $k_{0}=1$. 
Therefore it will be enough to prove Theorem \ref{thm:2.2}. The  proof will be  organised as follows.
Besides the system (\ref{eq:2.7}) we shall consider first the system described
by the range equations
\begin{subequations}
\begin{align}
&(i\oo\cdot\nn)b_{\nn}=\Phi_{\nn},\quad \nn\ne\vzero,
\label{eq:2.12a}\\
&(i\oo\cdot\nn)B_{\nn}=\Ga_{\nn},\quad \nn\ne\vzero,
\label{eq:2.12b}
\end{align}
\label{eq:2.12}
\end{subequations}
\vskip-.3truecm
\noindent
i.e. with no condition for $\nn=\vzero$. In Sections \ref{sec:3} and \ref{sec:4}
we shall prove that, if some further conditions (to be specified later on)
are found to be satisfied, it is possible to find, for $\e$ small enough and
arbitrary $\be_{0},B_{0}$, a solution
\begin{equation} \label{eq:2.13}
(\be_{0}+b(t),B_{0}+\widetilde{B}(t)),
\end{equation}
to the system (\ref{eq:2.12}), with $b(t)$ and $\widetilde B(t)$
as in (\ref{eq:2.5}) depending on the free parameters $\e,\be_{0},B_{0}$;
such a solution is obtained via a `resummation procedure', starting from
the formal solution of the range equations \eqref{eq:2.12}.
The conditions mentioned above can be illustrated as follows.
The resummation procedure turns out to be well-defined if the small divisors
of the resummed series can be bounded proportionally to the square of
the small divisors of the formal series. However, it is not obvious at all
that this is possible, since the latter are of the form
$(\ii\oo\cdot\nn)^{-1}$ with $\nn\in\ZZZ^{d}_{*}$, while the small divisors
of the resummed series are of the form $(\det((\ii\oo\cdot\nn)\uno-
\MM^{[n]}(\oo\cdot\nn;\e,\be_{0},B_{0})))^{-1}$,
for suitable $2\times2$ matrices $\MM^{[n]}$ (see Section \ref{sec:3}).
The bound on the small divisors of the resummed series is difficult to
check without assuming any
non-degeneracy condition on the perturbation.
Therefore we replace
$\MM^{[n]}(x;\e,\be_{0},B_{0})$ with $\MM^{[n]}(x;\e,\be_{0},B_{0})
\xi_{n}(\det(\MM^{[n]}(0;\e,\be_{0},B_{0})))$, for suitable
`cut-off functions' $\xi_{n}$, in such a way that the bound automatically
holds. The introduction
of the cut-offs changes the series in such a way that
if on the one hand the modified series are well-defined, on the other hand
in principle they no longer solve the range equations: this turns out to be
the case only if one can prove that the cut-offs can be removed. So, the
last part of the proof consists in showing that, by suitably choosing
the parameters $\be_{0},B_{0}$ as continuous functions of $\e$, this occurs
and moreover, for the same choice of $\be_{0},B_{0}$, the bifurcation
equations (\ref{eq:2.7c}) and (\ref{eq:2.7d}) hold;
hence for such $\be_{0},B_{0}$, the function
(\ref{eq:2.13}) is a solution of the whole system (\ref{eq:2.1}).
Once Theorem \ref{thm:2.2} is proved, Theorem \ref{thm:2.1} will
immediately follow taking $k_{0}=1$.

Next we shall see that if the system is Hamiltonian one can prove
the same result as in Theorem \ref{thm:2.2} with the only assumptions
in Hypotheses \ref{hyp1} and \ref{hyp2}. More precisely,
consider the Hamiltonian function
\begin{equation}\label{eq:2.14}
H(\aaa,\be,\AAA,B):=\oo\cdot\AAA + h_{0}(B)+\e f(\aaa,\be,B),
\end{equation}
where $(\aaa,\be)\in\TTT^{d+1}$ and $(\AAA,B)\in\RRR^{d}\times\gotB$,
with $\gotB$ an open subset of $\RRR$, are
canonically conjugate (action-angle) variables and $f\!:\TTT^{d+1}\times\gotB\to\RRR$
and $h_{0}\!:\gotB\to\RRR$ are real analytic functions.
Set $\om_{0}(B)=\partial_{B}h_{0}(B)$. Then the corresponding
Hamilton equations for the variables $(\be,B)$ are given by
\begin{equation}\label{eq:2.15}
\left\{\begin{aligned}
&\dot{\be}=\om_{0}(B)+\e \partial_{B}f(\oo t,\be,B),\\
&\dot{B}=-\e \partial_{\be}f(\oo t,\be,B).
\end{aligned}\right.
\end{equation}
We shall prove in Section \ref{sec:6} the following result.

\begin{theo}\label{thm:2.3} 
Consider the system (\ref{eq:2.15}) and
assume Hypotheses \ref{hyp1} and \ref{hyp2} to be satisfied. Then for $\e$ small enough
there exists at least one quasi-periodic solution $(\be(t),B(t))$ with frequency vector $\oo$.
Such a solution depends continuously on $\e$.
\end{theo}

Note that Hypothesis \ref{hyp2} is tantamount  to requiring $h_0$ to be convex.
Quasi-periodic solutions to (\ref{eq:2.15}) with frequency vector $\oo$ describe
lower-dimensional tori ($d$-dimensional tori for a system with $d+1$ degrees of freedom).
Such tori are parabolic in the sense that the ``normal frequency'' vanishes for $\e=0$.
Theorem \ref{thm:2.3} can be seen as the counterpart of Cheng's result \cite{Ch}
in the case in which all ``proper frequencies'' are fixed (isochronous
case) and the perturbation does not depend on the actions
conjugated to the ``fast angles'' (otherwise one should add a correction
like in \cite{M2}); moreover, with respect to \cite{Ch},
a weaker Diophantine condition is assumed on the proper frequencies.

\zerarcounters
\section{Diagrammatic rules}
\label{sec:3}

Let us consider the range equations \eqref{eq:2.12} and start by looking
for a quasi-periodic solution which can be formally written as
\begin{equation}\label{eq:3.1bis}
\begin{aligned}
&\be(t;\e,\be_{0},B_{0})=\be_{0}+\sum_{k\ge1}\e^{k}
b^{\{k\}}(t;\be_{0},B_{0})=
\be_{0}+\sum_{k\ge1}\e^{k}
\sum_{\nn\in\ZZZ^{d}_{*}}{\rm e}^{\ii\oo\cdot\nn t}b^{\{k\}}(\be_{0},B_{0}),\\
&B(t;\e,\be_{0},B_{0})=B_{0}+\sum_{k\ge1}\e^{k}
B^{\{k\}}(t;\be_{0},B_{0})=
B_{0}+\sum_{k\ge1}\e^{k}
\sum_{\nn\in\ZZZ^{d}_{*}}{\rm e}^{\ii\oo\cdot\nn t}B^{\{k\}}(\be_{0},B_{0}),
\end{aligned}
\end{equation}
where a different notation for the Taylor coefficients has been
used with respect to \eqref{eq:2.8} to stress that now we are
considering $B_{0}=B_{\vzero}$ as a parameter. If we define recurively for $k\ge1$
and $\nn\in\ZZZ^{d}_{*}$
\begin{equation} \nonumber
b^{\{k\}}_{\nn}(\be_{0},B_{0}):=\Phi^{\{k\}}_{\nn}(\be_{0},B_{0}), \qquad
B^{\{k\}}_{\nn}(\be_{0},B_{0}):=\Ga^{\{k\}}_{\nn}(\be_{0},B_{0}),
\end{equation}
where we have set
\begin{equation} \nonumber
\begin{aligned}
&\Phi^{\{k\}}_{\nn}(\be_{0},B_{0}):= \!\!
\sum_{s\ge1}\frac{1}{s!}\partial_{B}\om_{0}
(B_{0})\sum_{\substack{\nn_{1}+\ldots+\nn_{s}=\nn\\ \nn_{i}\in\ZZZ^{d}_{*}}}
\sum_{\substack{k_{1}+\ldots+k_{s}=k\\ k_{i}\ge1}}\prod_{i=1}^{s}B^{\{k_{i}\}}_{\nn_{i}}
(\be_{0},B_{0}) \\
&\qquad +\sum_{\substack{s\ge1 \\p+q=s}} \sum_{\substack{\nn_{0}+\ldots+\nn_{s}=\nn\\
\nn_{0}\in\ZZZ^{d},\nn_{i}\in\ZZZ^{d}_{*}}}\frac{1}{p!q!}\partial_{\be}^{p}
\partial_{B}^{q}F_{\nn_{0}}(\be_{0},B_{0}) 
\!\!\!\!\!\!\!\!\!\!\!
\sum_{\substack{k_{1}+\ldots+k_{s}=k-1\\ k_{i}\ge1}}
\prod_{i=1}^{p}b^{\{k_{i}\}}_{\nn_{i}}(\be_{0},B_{0}) \!\!\!\!
\prod_{i=p+1}^{s}B^{\{k_{i}\}}_{\nn_{i}}(\be_{0},B_{0}),\\
& \Ga^{\{k\}}_{\nn}(\be_{0},B_{0}):= \!\!\!\!
\sum_{\substack{s\ge1 \\ p+q=s}}\sum_{\substack{\nn_{0}+\ldots+\nn_{s}=\nn\\
\nn_{0}\in\ZZZ^{d},\nn_{i}\in\ZZZ^{d}_{*}}}\frac{1}{p!q!}\partial_{\be}^{p}
\partial_{B}^{q}G_{\nn_{0}}(\be_{0},B_{0}) 
\!\!\!\!\!\!\!\!\!\!\!
\sum_{\substack{k_{1}+\ldots+k_{s}=k-1\\ k_{i}\ge1}}
\prod_{i=1}^{p}b^{\{k_{i}\}}_{\nn_{i}}(\be_{0},B_{0}) \!\!\!\!
\prod_{i=p+1}^{s}B^{\{k_{i}\}}_{\nn_{i}}(\be_{0},B_{0}),
\end{aligned}
\end{equation}
for all $k\ge1$ and all $\nn\in\ZZZ^{d}$, then \eqref{eq:3.1bis}
turns out to be a formal solution to the range equations \eqref{eq:2.12}.
Note that we can see the formal expansion \eqref{eq:2.8} as obtained
from \eqref{eq:3.1bis} by solving the bifurcation equation
\eqref{eq:2.7c} and further expanding $B_{0}=B_{\vzero}(\e)$.

One could easily prove that \eqref{eq:3.1bis}
is formally well defined,  that is that the coefficients
$b^{\{k\}}(t;\be_{0},B_{0})$ and $B^{\{k\}}(t;\be_{0},B_{0})$ are
well defined to all orders $k\ge 1$; for instance one could
adapt the forthcoming diagrammatic formalism
(which would rather simplify with respect to the discussion below).
Unfortunately the power series may not be  convergent -- as far as we
know --, so we have to look for a different approach:
we shall see how to  construct a series, convergent if
$\be_{0},B_{0}$ are suitably chosen,
whose formal expansion coincides with (\ref{eq:2.8}). To this aim we shall
introduce a convenient graphical representation for the coefficients of
such a series. We start by introducing some notations.

\begin{figure}[ht] 
\centering 
\includegraphics[width=4in]{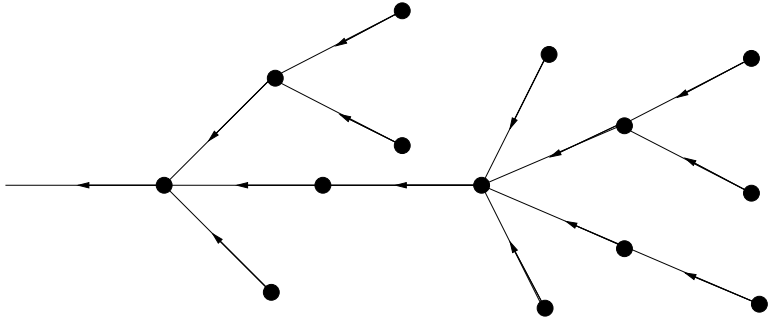} 
\vskip.2truecm 
\caption{A tree with 14 nodes.}
\label{fig:31} 
\end{figure}

A graph is a set of points and lines connecting them. 
A \emph{tree} $\theta$ is a connected graph with no cycle, 
such that all the lines are oriented toward a unique point 
(\emph{root}) which has only one incident line $\ell_{\theta}$ 
(\emph{root line}). 
All the points in a tree except the root are called \emph{nodes}. 
The orientation of the lines in a tree induces a partial ordering 
relation ($\preceq$) between the nodes and the lines: we can 
imagine that each line carries an arrow pointing toward the root. 
Given two nodes $v$ and $w$, 
we shall write $w \prec v$ every time $v$ is along the path 
(of lines) which connects $w$ to the root. 
When drawing a tree, we shall put the root to the extreme left
so that all the lines (and the corresponding arrows) will be directed
from right to left; see Figure \ref{fig:31}.

We denote by $N(\theta)$ and $L(\theta)$ the sets of nodes and 
lines in $\theta$ respectively. 
Since a line $\ell\in L(\theta)$ is uniquely identified 
by the node $v$ which it leaves, we may write $\ell = \ell_{v}$. 
We write $\ell_{w} \prec \ell_{v}$ if $w\prec v$ and $w\prec\ell=\ell_{v}$
if $w\preceq v$; if $\ell$ and $\ell'$ are two distinct comparable lines, i.e. 
$\ell' \prec \ell$, we denote by $\calP(\ell,\ell')$ the 
(unique) path of lines connecting $\ell'$ to $\ell$, with $\ell$ and 
$\ell'$ not included (in particular $\calP(\ell,\ell')=\emptyset$ 
if $\ell'$ enters the node $\ell$ exits). 

Given a tree $\theta$ we associate labels with the nodes and the
lines of $\theta$, as follows.

With each node $v\in N(\theta)$ we associate a \emph{mode label} 
$\nn_{v}\in \ZZZ^{d}$, a \emph{component label} $h_{v}\in\{\be,B\}$
and an \emph{order label} $k_{v}\in\{0,1\}$ with
the constraint that $k_{v}=1$ if $h_{v}=B$ or $\nn_{v}\ne\vzero$.
With each line $\ell=\ell_{v}$ we associate a pair of \emph{component labels}
$(e_{\ell},u_{\ell})\in\{\be,B\}^{2}$, with the constraint that $u_{\ell}
=h_{v}$, and a \emph{momentum} $\nn_{\ell}\in \ZZZ^{d}_{*}$, 
except for the root line which can have either zero momentum or not, 
i.e. $\nn_{\ell_{\theta}}\in\ZZZ^{d}$.
For any line $\ell$, we call $e_{\ell}$ and $u_{\ell}$ the
$e$-\emph{component} and the $u$-\emph{component} of $\ell$, respectively

We denote by $p_{v}$ and $q_{v}$ the numbers of lines with $e$-component
$\be$ and $B$, respectively, entering the node $v$ and set $s_{v}=p_{v}+q_{v}$. 
If $k_{v}=0$ for some $v \in N(\theta)$ we force also $p_{v}=0$ and $q_{v}\ge1$.

We impose the \emph{conservation law}
\begin{equation}\label{eq:3.1}
\nn_{\ell}=\sum_{v\prec \ell}\nn_{v}
\end{equation}
and we call \emph{order} of $\theta$ the number
\begin{equation}\label{eq:3.2}
k(\theta)=\sum_{v\in N(\theta)}k_{v}.
\end{equation}

Finally, we associate with each line $\ell$ also a \emph{scale label} $n_{\ell}$
such that $n_{\ell}=-1$ if $\nn_{\ell}=\vzero$, while 
$n_{\ell}\in\ZZZ_{+}$ if $\nn_{\ell}\neq\vzero$. 
Note that one can have $n_{\ell}=-1$ only if $\ell$ is the root line of
$\theta$. 

In the following we shall call simply trees the trees with labels and we shall use
the term \emph{unlabelled tree}  for the trees without labels. 
 
We shall say that two trees are \emph{equivalent} if they can be 
transformed into each other by continuously deforming the lines in 
such a way that these do not cross each other and also labels match. 
This provides an equivalence relation on the set of the trees.
 From now on we shall call trees tout court such equivalence classes. 
 
A subset $T \subset \theta$ will be called a \emph{subgraph} of
$\theta$ if it is formed by a set of nodes $N(T)\subseteq N(\theta)$
and a set of of lines $L(T)\subseteq L(\theta)$ connecting them (possibly
including the root line of $\theta$) in such a way that $N(T) \cup L(T)$ is
connected. We call \emph{order} of $T$ the number
\begin{equation}\label{eq:3.3}
k(T)=\sum_{v\in N(T)}k_{v}.
\end{equation}
We say that a line enters $T$ if it 
connects a node $v\notin N(T)$ to a node $w\in N(T)$
and we say  that a line exits $T$ if it connects a node $v\in N(T)$ 
to a node $w\notin N(T)$. Of course, if a line $\ell$ enters or exits $T$,
then $\ell\notin L(T)$. 
If $\theta$ is a labelled tree and $T$ a
subgraph of $\theta$, then $T$ inherits the labels of $\theta$.

\begin{figure}[ht] 
\centering 
\ins{020pt}{-050pt}{(a)}
\ins{174pt}{-015pt}{$0$}
\ins{159pt}{-029pt}{$0$}
\ins{136pt}{-034pt}{$0$}
\ins{195pt}{-061pt}{$0$}
\ins{118pt}{-086pt}{$2$}
\ins{156pt}{-092pt}{$2$}
\ins{079pt}{-066pt}{$4$}
\ins{081pt}{-086pt}{$4$}
\ins{098pt}{-066pt}{$4$}
\ins{142pt}{-098pt}{$4$}
\ins{045pt}{-086pt}{$6$}
\ins{110pt}{-048pt}{$9$}
\ins{164pt}{-076pt}{$9$}
\ins{240pt}{-050pt}{(b)}
\ins{404pt}{-018pt}{$0$}
\ins{428pt}{-049pt}{$0$}
\ins{362pt}{-084.5pt}{$2$}
\ins{340pt}{-106pt}{$4$}
\includegraphics[width=6.0in]{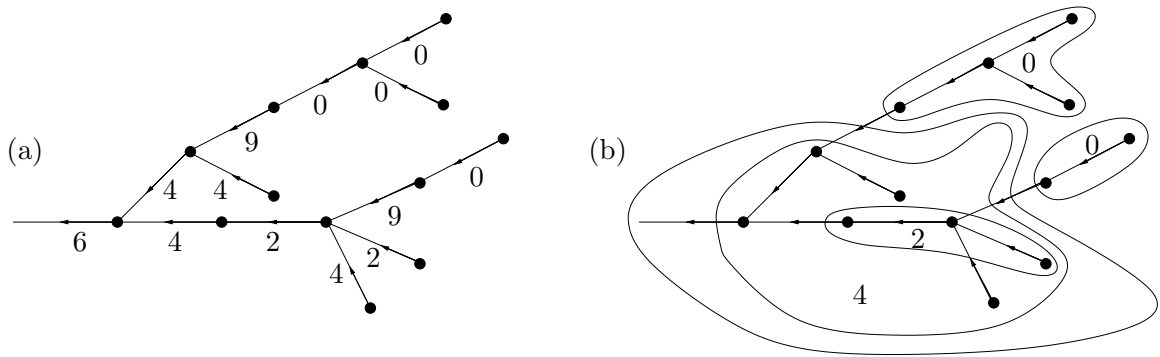}
\caption{Example of clusters: (a) a tree $\theta$ is represented with
(only) the scale labels associated with its lines; (b) the clusters in $\theta$,
corresponding to the same assignment of scale labels, are drawn.}
\label{fig:32}
\end{figure}

A \emph{cluster} $T$ on scale $n$ is a maximal subgraph 
of a tree $\theta$ such that all the lines have scales 
$n'\le n$ and there is at least one line with scale $n$. 
The lines entering the cluster $T$ and the line coming 
out from it (unique if existing at all) are called the \emph{external}
lines of $T$. An example of clusters is in Figure \ref{fig:32}.

\begin{figure}[ht] 
\centering 
\ins{135pt}{-14pt}{$n_{1}$} 
\ins{237pt}{-22pt}{$n$} 
\ins{330pt}{-14pt}{$n_{2}$} 
\ins{135pt}{-35pt}{$\ell_{T}$}
\ins{290pt}{-50pt}{$T$}
\ins{330pt}{-35pt}{$\ell_{T}'$}
\includegraphics[width=4in]{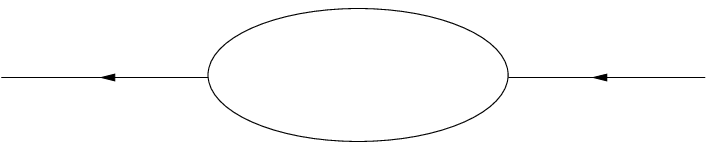} 
\vskip.2truecm 
\caption{Graphical representation of the self-energy clusters $T$ on scale $n$;
by construction $n_{1},n_{2} \ge n+1$ and $\nn_{\ell_{T}}=\nn_{\ell_{T}'}$.
Note that neither $\ell_{T}'$ (the entering line) nor $\ell_{T}$
(the exiting line) belong to $L(T)$.}
\label{fig:33} 
\end{figure}

\begin{figure}[ht] 
\centering 
\ins{010pt}{-02pt}{$(a)$} 
\ins{160pt}{-04pt}{$(b)$} 
\ins{360pt}{-06pt}{$(c)$}
\ins{032pt}{-44pt}{$\ell_{T}$}
\ins{083pt}{-57pt}{$T$}
\ins{097pt}{-44pt}{$\ell_{T}'$}
\ins{178pt}{-44pt}{$\ell_{T}$}
\ins{262pt}{-57pt}{$T$}
\ins{282pt}{-44pt}{$\ell_{T}'$}
\ins{360pt}{-44pt}{$\ell_{T}$}
\ins{398pt}{-57pt}{$T$}
\ins{430pt}{-46pt}{$\ell_{T}'$}
\includegraphics[width=6in]{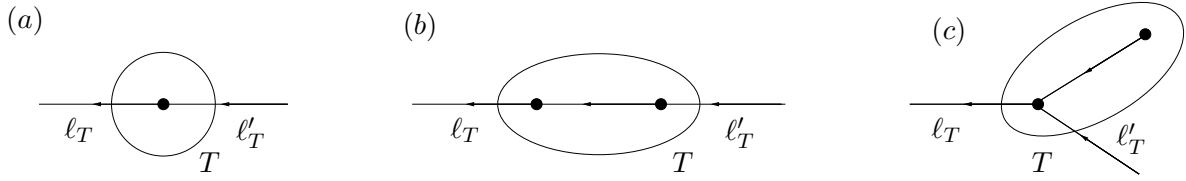} 
\vskip.2truecm 
\caption{Self-energy clusters $T$ with their external lines $\ell_{T}$
(exiting line) and $\ell_{T}'$ (entering line),
for $|N(T)|=1$ (case (a)) and $|N(T)|=2$ (cases (b) and (c)),
Note that in case (a) one has either $k(T)=0$ or $k(T)=1$, while in cases (b)
and (c) one has either $k(T)=1$ or $k(T)=2$; in all cases, when
$k(T)<|N(T)|$, one node $v\in N(T)$ has $k_{v}=0$. The scale of $T$ is
$n=-1$ in case (a) and can be any $n\in\ZZZ_{+}$ in cases (b) and (c).}
\label{fig:34} 
\end{figure}

A \emph{self-energy cluster} is a cluster $T$ such that 
(i) $T$ has only one entering line $\ell'_{T}$ and 
one exiting line $\ell_{T}$, (ii) one has $\nn_{\ell_{T}}= 
\nn_{\ell'_{T}}$ and hence 
\begin{equation}\label{eq:3.4} 
\sum_{v\in N(T)}\nn_{v}=\vzero. 
\end{equation} 
Self-energy clusters will be represented graphically as in Figure \ref{fig:33}. 
Examples of low order self-energy clusters are given in Figure \ref{fig:34}. 

For any self-energy cluster $T$, set $\calP_{T}= 
\calP(\ell_{T},\ell'_{T})$. More generally, if $T$ is a subgraph 
of $\theta$ with only one entering line $\ell'$ and one exiting line 
$\ell$, we set $\calP_{T}=\calP(\ell,\ell')$. 
We shall say that a self-energy cluster $T$ is on 
scale $-1$, if $N(T)=\{v\}$ (with of course $\nn_{v}=\vzero$),
so that $\calP_{T}=\emptyset$. If a self-energy cluster is on a scale $n\ge0$
then $|N(T)|\ge 2$ and $k(T)\ge 1$, as is easy to check.

\begin{rmk}\label{rmk:3.1} 
\emph{ 
Given a self-energy cluster $T$, the momenta of the lines in $\calP_{T}$ 
depend on $\nn_{\ell'_{T}}$ because of the conservation law (\ref{eq:3.1}). 
More precisely, for all $\ell\in\calP_{T}$ one has 
$\nn_{\ell}=\nn_{\ell}^{0}+\nn_{\ell'_{T}}$ with 
\begin{equation}\nonumber
\nn_{\ell}^{0}=\sum_{\substack{w\in N(T) \\ w\prec \ell}} 
\nn_{w},
\end{equation} 
while all the other momenta in $T$ do not depend on $\nn_{\ell'_{T}}$.
} 
\end{rmk} 
 
We shall say that two self-energy clusters $T_{1},T_{2}$ have the same 
\emph{structure} if setting $\nn_{\ell'_{T_{1}}}=\nn_{\ell'_{T_{2}}}=\vzero$ one has 
$T_{1}=T_{2}$. This provides an equivalence relation on the 
set of all self-energy clusters.
From now on we shall call self-energy
clusters tout court such equivalence  classes. 
 
A \emph{renormalised tree} is a tree in which no self-energy 
clusters appear; analogously a \emph{renormalised subgraph} is a 
subgraph of a tree $\theta$ which does not contains any self-energy cluster. 
Note that if $T$ is a renormalised self-energy cluster and $N(T)\ge2$ then $k(T)\ge2$.

Given a tree $\theta$ we call \emph{total momentum} of 
$\theta$ the momentum associated with $\ell_{\theta}$ and
\emph{total component} of $\theta$ the $e$-component of $\ell_{\theta}$.
We denote by $\Theta_{k,\nn,h}^{\RR}$ the set of renormalised trees
with order $k$, total momentum $\nn$ and total component $h$;
the set of renormalised self-energy clusters $T$ on scale $n$ such that
$u_{\ell_{T}}=u$ and $e_{\ell'_{T}}=e$ will be denoted by $\gotR_{n,u,e}$.

\begin{lemma}\label{lem:3.2}
Let $T$ be a subgraph of any tree $\theta$. Then one has $|N(T)|\le 3k(T)-1$.
\end{lemma}

\prova
We shall prove the result by induction on $k=k(T)$.
For $k=1$ the bound is trivially satisfied as a direct check shows.
Assume then the bound to hold for all $k'<k$. Call $v$ the node
which $\ell_{T}$ (possibly $\ell_{\theta}$) exits, $\ell_{1},\ldots,
\ell_{s_{v}}$ the lines
entering $v$ and $T_{1},\ldots,T_{s_{v}}$ the subgraphs of $T$ with
exiting lines $\ell_{1},\ldots,\ell_{s_{v}}$.
If $k_{v}=1$ then by the inductive hypothesis one has
$$
|N(T)|=1+\sum_{i=1}^{s_{v}}|N(T_{i})|\le 1+3(k-1)-s_{v}.
$$
If $k_{v}=0$ one has $s_{v}=q_{v}\ge2$ and hence
$$
|N(T)|=1+\sum_{i=1}^{q_{v}}|N(T_{i})|\le 1+3k-q_{v}.
$$
so that in both cases the bound follows.\EP

\begin{rmk}\label{rmk:3.3}
\emph{For any subgraph $T$, one has
$$
\sum_{v\in N(T)}s_{v}\le |L(T)|+1\le |N(T)|+1\le 3k(T).
$$
In particular $\sum_{v\in N(\theta)}s_{v}\le 3k(\theta)$.
}
\end{rmk}

For any $\theta\in\Theta^{\RR}_{k,\nn,h}$ we associate 
with each node $v\in N(\theta)$ a \emph{node factor} 
\begin{equation} \label{eq:3.5}
\calF_{v}=\calF_{v}(\be_{0},B_{0}):=\left\{\begin{aligned}
&\frac{1}{p_{v}!q_{v}!}\partial_{\be}^{p_{v}}\partial_{B}^{q_{v}}
F_{\nn_{v}}(\be_{0},B_{0}), &h_{v}=\be, k_{v}=1,\\
&\frac{1}{q_{v}!}\partial^{q_{v}}\om_{0}(B_{0}),&h_{v}=\be, k_{v}=0,\\
&\frac{1}{p_{v}!q_{v}!}\partial_{\be}^{p_{v}}\partial_{B}^{q_{v}}
G_{\nn_{v}}(\be_{0},B_{0}),&h_{v}=B, k_{v}=1.\\
\end{aligned}\right.
\end{equation}
With each line $\ell=\ell_{v}$ we associate a \emph{propagator}
$\calG^{[n_{\ell}]}_{e_{\ell},u_{\ell}}(\oo\cdot\nn_{\ell};\e,\be_{0},B_{0})$
defined recursively as follows. 

\begin{figure}[ht] 
\centering 
\ins{370pt}{-132pt}{$x$} 
\ins{320pt}{-130pt}{$\displaystyle{\frac{\al_{0}}{8}}$} 
\ins{220pt}{-130pt}{$\displaystyle{\frac{\al_{0}}{16}}$} 
\ins{195pt}{-130pt}{$\displaystyle{\frac{\al_{m_{1}}}{8}}$} 
\ins{150pt}{-130pt}{$\displaystyle{\frac{\al_{m_{1}}}{16}}$} 
\ins{110pt}{-130pt}{$\displaystyle{\frac{\al_{m_{2}}}{16}}$} 
\ins{135pt}{-10pt}{$\Psi_{2}(x)$} 
\ins{210pt}{-10pt}{$\Psi_{1}(x)$} 
\ins{330pt}{-10pt}{$\Psi_{0}(x)$} 
\includegraphics[width=4in]{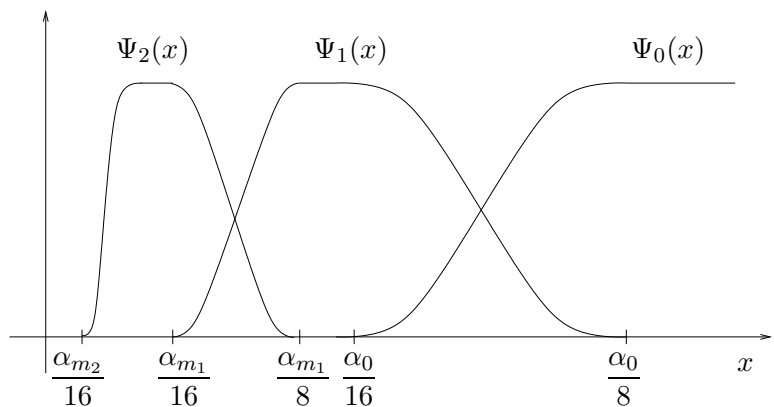} 
\vskip.2truecm 
\caption{Graphs of some of the $C^{\io}$ functions $\Psi_{n}(x)$
partitioning the unity in $\RRR\setminus\{0\}$; here $\al_{m}=\al_{m}(\oo)$. 
The function $\chi_{0}(x)=\chi(8x/\al_{0})$ is given by the sum of 
all functions $\Psi_{n}(x)$ for $n\ge 1$.} 
\label{fig:36} 
\end{figure}

Let us introduce the sequences $\{m_{n},p_{n}\}_{n \ge 0}$, with $m_{0}=0$ 
and, for all $n\ge 0$, $m_{n+1}=m_{n}+p_{n}+1$,  where
$p_{n}:=\max\{q\in\ZZZ_{+}\,:\,\al_{m_{n}}(\oo)<2\al_{m_{n}+q}(\oo)\}$. Then the 
subsequence $\{\al_{m_{n}}(\oo)\}_{n\ge 0}$ of $\{\al_{m}(\oo)\}_{m\ge0}$ is decreasing. 
Let $\chi:\RRR\to\RRR$ be a $C^{\io}$ function, non-increasing for $x\ge0$ and
non-decreasing for $x<0$, such that 
\begin{equation}\label{eq:3.6} 
\chi(x)=\left\{ 
\begin{aligned} 
&1,\qquad |x| \le 1/2, \\ 
&0,\qquad |x| \ge 1. 
\end{aligned}\right. 
\end{equation} 
Set $\chi_{-1}(x)=1$ and $\chi_{n}(x)=\chi(8x/\al_{m_{n}}(\oo))$ for $n\ge0$. 
Set also $\psi(x)=1-\chi(x)$, $\psi_{n}(x)=\psi(8x/\al_{m_{n}}(\oo))$, 
and $\Psi_{n}(x)=\chi_{n-1}(x)\psi_{n}(x)$, for $n\ge 0$; see Figure
\ref{fig:36}. 
 
\begin{lemma}\label{lem:3.4} 
For all $x\neq 0$ and for all $p\ge0$ one has 
\begin{equation} \nonumber 
\psi_{p}(x)+\sum_{n\ge p+1}\Psi_{n}(x)=1. 
\end{equation} 
\end{lemma} 

\prova For fixed $x\neq 0$ let $N=N(x):=\min\{n\,:\,\chi_{n}(x)= 0\}$ and 
note that $\max\{n\,:\,\psi_{n}(x)=0\} \le N-1$. Then if $p\le N-1$ 
\begin{equation} \nonumber 
\psi_{p}(x)+\sum_{n\ge p+1}\Psi_{n}(x)= \psi_{N-1}(x)+\chi_{N-1}(x)=1, 
\end{equation} 
while if $p\ge N$ one has 
\begin{equation} \nonumber 
\psi_{p}(x)+\sum_{n\ge p+1}\Psi_{n}(x)= \psi_{p}(x)=1. 
\end{equation} 
\vskip-1.1truecm\EP 

Then for $n\ge0$ we set formally
\begin{equation}\label{eq:3.7}
\begin{aligned}
\calG^{[n]}(x;\e,\be_{0},B_{0})&=
\begin{pmatrix}
\calG^{[n]}_{\be,\be}(x;\e,\be_{0},B_{0}) &
\calG^{[n]}_{\be,B}(x;\e,\be_{0},B_{0}) \cr
& \cr
\calG^{[n]}_{B,\be}(x;\e,\be_{0},B_{0}) &
\calG^{[n]}_{B,B}(x;\e,\be_{0},B_{0})
\end{pmatrix} \\ &
:=\Psi_{n}(x)\left((\ii x)\uno-\MM^{[n-1]}
(x;\e,\be_{0},B_{0})\right)^{-1},
\end{aligned}
\end{equation}
where (here and henceforth) $\uno$ is the $2\times2$ identity matrix and
\begin{equation}\label{eq:3.8}
\MM^{[n-1]}(x;\e,\be_{0},B_{0}):=\sum_{q=-1}^{n-1}\chi_{q}(x)M^{[q]}(x;\e,\be_{0},B_{0}),
\end{equation}
where, for $n\ge -1$, $M^{[n]}(x;\e,\be_{0},B_{0})$ is the $2\times2$ matrix
\begin{equation}\label{eq:3.9}
M^{[n]}(x;\e,\be_{0},B_{0}):=\begin{pmatrix}
M^{[n]}_{\be,\be}(x;\e,\be_{0},B_{0}) &
M^{[n]}_{\be,B}(x;\e,\be_{0},B_{0}) \cr
& \cr
M^{[n]}_{B,\be}(x;\e,\be_{0},B_{0}) &
M^{[n]}_{BB}(x;\e,\be_{0},B_{0})\end{pmatrix},
\end{equation}
with formally
\begin{equation}\label{eq:3.10}
M^{[n]}_{u,e}(x;\e,\be_{0},B_{0}):=\sum_{T\in\gotR_{n,u,e}}\e^{k(T)}\Val_{T}
(x;\e,\be_{0},B_{0}),
\end{equation}
and $\Val_{T}(x;\e,\be_{0},B_{0})$ is the \emph{renormalised value} of $T$, defined as
\begin{equation}\label{eq:3.11}
\Val_{T}(x;\e,\be_{0},B_{0}):=\Biggl(\prod_{v\in N(T)}\calF_{v}\Biggr)
\Biggl(\prod_{\ell\in L(T)}\calG^{[n_{\ell}]}_{e_{\ell},u_{\ell}}(\oo\cdot\nn_{\ell};\e,\be_{0},B_{0})\Biggr).
\end{equation}

Here and henceforth, the sums and the products over empty sets have to be
considered as zero and $1$, respectively.
Note that $\Val_{T}$ depends on $\e$ -- because the propagators do --; 
moreover it depends on $x=\oo\cdot\nn_{\ell'_{T}}$ only through the propagators 
associated with the lines $\ell\in\calP_{T}$ (see Remark \ref{rmk:3.1}).

Set $\MM:=\{\MM^{[n]}(x;\e,\be_{0},B_{0})\}_{n\ge-1}$. We call
\emph{self-energies} the matrices $\MM^{[n]}(x;\e,\be_{0},B_{0})$.

\begin{rmk}\label{rmk:3.5} 
\emph{ 
One has 
\begin{equation}\nonumber
\partial_{c} \calG^{[n]}_{e,u}(x;\e,\be_{0},B_{0})=
\biggl(\calG^{[n]}(x;\e,\be_{0},B_{0})
\partial_{c}\MM^{[n-1]}(x;\e,\be_{0},B_{0})
\left( (\ii x)\uno-\MM^{[n-1]}(x;\e,\be_{0},B_{0})\right)^{-1}\biggr)_{e,u}
\end{equation}
for both $c=\be_{0},B_{0}$.
}
\end{rmk} 

Setting also $\calG^{[-1]}(0;\e,\be_{0},B_{0})=\uno$, for any subgraph $S$
of any $\theta\in \Theta^{\RR}_{k,\nn,h}$
define the \emph{renormalised value} of $S$ as 
\begin{equation}\label{eq:3.12} 
\Val(S;\e,\be_{0},B_{0}) := \Biggl(\prod_{v\in N(S)}\calF_{v} \Biggr)
\Biggl( \prod_{\ell\in L(S)} \calG^{[n_{\ell}]}_{e_{\ell},u_{\ell}}(\oo\cdot\nn_{\ell}; \e,\be_{0},B_{0})\Biggr). 
\end{equation} 

We define
\begin{equation}\label{eq:3.13}
b^{[k]}_{\nn}(\e,\be_{0},B_{0}):=
\sum_{\theta\in\Theta_{k,\nn,\be}^{\RR}}\Val(\theta;\e,\be_{0},B_{0}),
\qquad
B^{[k]}_{\nn}(\e,\be_{0},B_{0}):=
\sum_{\theta\in\Theta_{k,\nn,B}^{\RR}}\Val(\theta;\e,\be_{0},B_{0}),
\end{equation}
for any $\nn\ne\vzero$, and
\begin{equation}\label{eq:3.14}
\Phi^{[k]}_{\vzero}(\e,\be_{0},B_{0}):=
\sum_{\theta\in\Theta_{k,\vzero,\be}^{\RR}}\Val(\theta;\e,\be_{0},B_{0}),
\qquad
\Ga^{[k]}_{\vzero}(\e,\be_{0},B_{0}):=
\sum_{\theta\in\Theta_{k,\vzero,B}^{\RR}}\Val(\theta;\e,\be_{0},B_{0}) .
\end{equation}
Set (again formally)
\begin{equation}\label{eq:3.15}
\begin{aligned}
&b^{\RR}(t;\e,\be_{0},B_{0}):=\sum_{k\ge1}\e^{k}\sum_{\nn\in\ZZZ^{d}_{*}}
{\rm e}^{\ii\nn\cdot\oo t}b^{[k]}_{\nn}(\e,\be_{0},B_{0}),\\
&\widetilde{B}^{\RR}(t;\e,\be_{0},B_{0}):=\sum_{k\ge1}\e^{k}\sum_{\nn\in\ZZZ^{d}_{*}}
{\rm e}^{\ii\nn\cdot\oo t}B^{[k]}_{\nn}(\e,\be_{0},B_{0}),
\end{aligned}
\end{equation}
and 
\begin{equation}\label{eq:3.16}
\Phi^{\RR}_{\vzero}(\e;\be_{0},B_{0}):=\sum_{k\ge0}\e^{k}
\Phi^{[k]}_{\vzero}(\e,\be_{0},B_{0}),  \qquad
\Ga^{\RR}_{\vzero}(\e;\be_{0},B_{0}):=\sum_{k\ge0}\e^{k}
\Ga^{[k]}_{\vzero}(\e,\be_{0},B_{0}) ,
\end{equation}
and define $\be^{\RR}(t;\e,\be_{0},B_{0})=\be_{0}+b^{\RR}(t;\e,\be_{0},B_{0})$ and
$B^{\RR}(t;\e,\be_{0},B_{0})=B_{0}+\widetilde{B}^{\RR}(t;\e,\be_{0},B_{0})$.
Set also
$\Theta_{k,\nn,h}^{\RR,n}=\{\theta\in\Theta_{k,\nn,h}^{\RR}:n_{\ell}
\le n\mbox{ for all }\ell\in L(\theta) \}$ and define
\begin{equation}\label{eq:3.17}
\begin{aligned}
&\Phi^{\RR,n}_{\vzero}(\e;\be_{0},B_{0}):=\sum_{k\ge0}\e^{k}
\sum_{\theta\in\Theta_{k,\vzero,\be}^{\RR,n}}\Val(\theta;\e,\be_{0},B_{0}),\\
&\Ga^{\RR,n}_{\vzero}(\e;\be_{0},B_{0}):=\sum_{k\ge0}\e^{k}
\sum_{\theta\in\Theta_{k,\vzero,B}^{\RR,n}}\Val(\theta;\e,\be_{0},B_{0}).
\end{aligned}
\end{equation}

\begin{rmk}\label{rmk:3.6}
\emph{
One has
\begin{equation}\nonumber
\MM^{[-1]}(x;\e,\be_{0},B_{0})=M^{[-1]}(x;\e,\be_{0},B_{0})=
\begin{pmatrix}\e\partial_{\be_{0}}F_{\vzero}(\be_{0},B_{0}) &
\om_{0}'(B_{0})+\e\partial_{B_{0}}F_{\vzero}(\be_{0},B_{0}) \cr
& \cr
\e\partial_{\be_{0}}G_{\vzero}(\be_{0},B_{0}) &
\e\partial_{B_{0}}G_{\vzero}(\be_{0},B_{0})
\end{pmatrix},
\end{equation}
where $\om_{0}'({B}_{0}) \neq 0$ for $B_{0}$ close enough to $\ol{B}_{0}$
by Hypothesis \ref{hyp2}.
In particular $\MM^{[-1]}(x;\e,\be_{0},B_{0})$ does not depend on $x$ and
is a real-valued matrix.
}
\end{rmk}

\begin{rmk}\label{rmk:3.7} 
\emph{ 
If $T$ is a renormalised self-energy cluster, then $\Val(T;\e,\be_{0},B_{0})=
\Val_{T}(\oo\cdot\nn_{\ell'_{T}};\e,\be_{0},B_{0})$.} 
\end{rmk} 

\begin{rmk} \label{rmk:3.8} 
\emph{ 
Given a renormalised tree $\theta$ such that $\Val(\theta;\e,\be_{0},B_{0}) 
\neq 0$, for any line $\ell\in L(\theta)$ (except possibly the root 
line) one has $\Psi_{n_{\ell}}(\oo\cdot\nn_{\ell})\neq0$ and hence 
\begin{equation} \label{eq:3.18} 
\frac{\al_{m_{n_{\ell}}}(\oo)}{16} < |\oo\cdot\nn_{\ell}| < 
\frac{\al_{m_{n_{\ell}-1}}(\oo)}{8} , 
\end{equation} 
where $\al_{m_{-1}}(\oo)$ has to be interpreted as $+\io$. 
Note also that $\Psi_{n_{\ell}}(\oo\cdot\nn_{\ell})\neq0$ implies
$$
|\oo\cdot\nn_{\ell}| < \frac{1}{8} \al_{m_{n_{\ell}-1}}(\oo)  < 
\frac{1}{4} \al_{m_{n_{\ell}-1}+p_{n_{\ell}-1}}(\oo)  = 
\frac{1}{4} \al_{m_{n_{\ell}}-1}(\oo) < \al_{m_{n_{\ell}}-1}(\oo) $$
and hence, by definition of $\al_{m}(\oo)$, one has $|\nn_{\ell}|> 2^{m_{n_{\ell}}-1}$.
Moreover, by the definition of $\{\al_{m_{n}}(\oo)\}_{n\ge 0}$, the 
number of scales which can be associated with a line $\ell$ in such a way 
that the propagator does not vanishes is at most 2. 
The same considerations apply to any subgraph of $\theta$ and to 
any renormalised self-energy cluster. 
} 
\end{rmk} 

For any renormalised subgraph $S$ of any tree $\theta$ we denote by
$\gotN_{n}(S)$ the number of lines on scale $\ge n$ in $S$ and set
\begin{equation}\nonumber
K(S):=\sum_{v\in N(S)}|\nn_{v}|.
\end{equation}
%

\begin{lemma}\label{lem:3.9}
For any $h\in\{\be,B\}$, $\nn\in\ZZZ^{d}$, $k\ge 1$ and for any
$\theta\in \Theta_{k,\nn,h}^{\RR}$ such
that $\Val(\theta;\e,\be_{0},B_{0})$ $\ne 0$, one has $\gotN_{n}(\theta)\le 
2^{-(m_{n}-2)}K(\theta)$ for all $n\ge0$.
\end{lemma}

\prova 
First of all we note that if $\gotN_{n}(\theta)\ge 1$, then there is
at least one line $\ell$ with $n_{\ell}=n$ and hence
$K(\theta)\ge |\nn_{\ell}|\ge 2^{m_{n}-1}$ (see Remark \ref{rmk:3.8}). 
Now we prove the bound 
$\gotN_{n}(\theta)\le \max\{2^{-(m_{n}-2)}K(\theta)-1,0\}$ by induction on 
the order. 
 
If the root line of $\theta$ has scale $n_{\ell_{\theta}}<n$ then the bound 
follows by the inductive hypothesis. If $n_{\ell_{\theta}}\ge n$, 
call $\ell_{1},\ldots,\ell_{r}$ the lines 
with scale $\ge n$ closest to $\ell_{\theta}$ (that is such that 
$n_{\ell'} < n$ for all lines $\ell'\in\calP(\ell_{\theta},\ell_{i})$, 
$i=1,\ldots,r$). If $r=0$ then $\gotN_{n}(\theta)=1$ 
and $|\nn| \ge 2^{m_{n}-1}$, so that the bound follows. 
If $r\ge 2$ the 
bound follows once more by the inductive hypothesis. If $r=1$, then 
$\ell_{1}$ is the only entering line of a cluster $T$ which is not a 
renormalised self-energy cluster as $\theta\in\Theta^{\RR}_{k,\nn,h}$ and hence 
$\nn_{\ell_{1}}\neq\nn$. But then 
$$ 
|\oo\cdot(\nn-\nn_{\ell_{1}})|\le|\oo\cdot\nn|+ 
|\oo\cdot\nn_{\ell_{1}}|\le \frac{1}{4}\al_{m_{n-1}}(\oo)< 
\al_{m_{n-1}+p_{n-1}}(\oo)=\al_{m_{n}-1}(\oo) , 
$$ 
as both $\ell_{\theta}$ and $\ell_{1}$ are on scale $\ge n$, so that one 
has $K(T)\ge|\nn-\nn_{\ell_{1}}|\ge 2^{m_{n}-1}$. Now, call 
$\theta_{1}$ the subtree of $\theta$ with root line $\ell_{1}$. Then one has 
$\gotN_{n}(\theta)=1+\gotN_{n}(\theta_{1}) \le 1+ 
\max\{2^{-(m_{n}-2)}K(\theta_{1})-1,0\}$,
so that  $\gotN_{n}(\theta) 
\le 2^{-(m_{n}-2)}(K(\theta)-K(T))\le 2^{-(m_{n}-2)}K(\theta)-1$,
again by induction.\EP 

\begin{lemma}\label{lem:3.10}
For any $e,u\in\{\be,B\}$, $n \ge0$ and for any $T\in\gotR_{n,u,e}$ such that
$\Val_{T}(x;\e,\be_{0},B_{0})\ne 0$, one has $K(T) > 2^{m_{n}-1}$ and
$\gotN_{p}(T)\le 2^{-(m_{p}-2)} K(T)$ for $0\le p\le n$.
\end{lemma}

\prova 
We first prove that for all $n\ge 0$ and all $T\in\gotR_{n,u,e}$, one 
has $K(T)\ge 2^{m_{n}-1}$. In fact if $T\in\gotR_{n,u,e}$ then $T$ contains 
at least a line on scale $n$. If there is $\ell\in L(T)\setminus\calP_{T}$ 
with $n_{\ell}=n$, then $K(T)\ge|\nn_{\ell}|>2^{m_{n}-1}$ (see Remark \ref{rmk:3.8}).  Otherwise, let 
$\ell\in\calP_{T}$ be the line on scale $n$ which is closest to $\ell'_{T}$. 
Call $\widetilde{T}$ the subgraph (actually the cluster) consisting of all 
lines and nodes of $T$ preceding $\ell$. Then 
$\nn_{\ell}\neq \nn_{\ell'_{T}}$, otherwise $\widetilde{T}$ would be a 
renormalised self-energy cluster. Therefore $K(T)>|\nn_{\ell}-\nn_{\ell'_{T}}|>2^{m_{n}-1}$ 
as both $\ell,\ell'_{T}$ are on scale $\ge n$. 

Given a tree $\theta$, call $\CCCC(n,p)$ the set of  
renormalised subgraphs $T$ of 
$\theta$ with only one entering line $\ell'_{T}$ and one exiting line
$\ell_{T}$ both on scale $\ge p$, such that $L(T)\neq\emptyset$ and
$n_{\ell}\le n$ for any $\ell\in L(T)$. 
Note that $\gotR_{n,u,e}\subset \CCCC(n,p)$ for all $n,p\ge 0$ and $u,e\in\{\be_{0},B_{0}\}$.
We prove that $\gotN_{p}(T)\le \max\{K(T)2^{-(m_{p}-2)}-1,0\}$ 
for $0\le p\le n$ and all $T\in\CCCC(n,p)$. The proof is by induction on
the order. 
Call $N(\calP_{T})$ the set of nodes in $T$ connected by lines in 
$\calP_{T}$. If all lines in $\calP_{T}$ are on scale $< p$, then 
$\gotN_{p}(T)=\gotN_{p}(\theta_{1})+\ldots+\gotN_{p}(\theta_{r})$ 
if $\theta_{1},\ldots,\theta_{r}$ are the subtrees with root line 
entering a node in $N(\calP_{T})$ and hence the bound follows 
from (the proof of) Lemma \ref{lem:3.9}. 
If there exists a line $\ell\in\calP_{T}$ on scale $\ge p$, call 
$T_{1}$ and $T_{2}$ the subgraphs of $T$ such that 
$L(T)=\{\ell\}\cup L(T_{1})\cup L(T_{2})$. Note that if $L(T_{1}),L(T_{2})
\neq\emptyset$, then $T_{1},T_{2}\in\CCCC(n,p)$. 
Hence, by the inductive hypothesis one has 
\begin{equation*} 
\gotN_{p}(T) =1+\gotN_{p}(T_{1})+\gotN_{p}(T_{2}) 
\le 1+\max\{2^{-(m_{p}-2)}K(T_{1})-1,0\}+\max\{2^{-(m_{p}-2)}K(T_{2})-1,0\}. 
\end{equation*} 
If both $\gotN_{p}(T_{1}),\gotN_{p}(T_{2})$ are zero the bound
follows as $K(T) \ge 2^{m_{p}-1}$, while if both are non-zero one has
$\gotN_{p}(T)\le2^{-(m_{p}-2)}(K(T_{1})+K(T_{2}))-1 = 2^{-(m_{p}-2)}K(T)-1$.
Finally if only one is zero, say $\gotN_{p}(T_{1})\ne0$ and $\gotN_{p}(T_{2})=0$, then 
$\gotN_{p}(T) \le 2^{-(m_{p}-2)}K(T_{1}) = 2^{-(m_{p}-2)}K(T)-2^{-(m_{p}-2)}K(T_{2})$.
On the other hand, in such a case $T_{2}$ is a cluster and hence
$\nn_{\ell}\neq \nn_{\ell_{T}'}$, which implies $K(T_{2}) \ge 2^{m_{p}-1}$.
The same argument can be used in the case $\gotN_{p}(T_{1})=0$ and
$\gotN_{p}(T_{2}) \neq 0$.\EP 

\begin{rmk}\label{rmk:3.11}
\emph{
Inequality (\ref{eq:3.18}) has been repeatedly used in the proof of
Lemmas \ref{lem:3.9} and \ref{lem:3.10}. In fact the proof works -- as
one can easily check -- under the weaker condition that
\begin{equation} \label{eq:3.19}
\frac{\al_{m_{n_{\ell}}}(\oo)}{32} < |\oo\cdot\nn_{\ell}| < 
\frac{\al_{m_{n_{\ell}-1}}(\oo)}{4}  
\end{equation} 
as long as $\Psi_{n_{\ell}}(\oo\cdot\nn_{\ell}) \neq 0$. This
observation will be used later on (see Lemma \ref{lem:4.11} below).
}
\end{rmk}

\zerarcounters 
\section{Convergence of the resummed series: dimensional bounds} 
\label{sec:4} 

Now we shall prove that, under the assumption that the propagators
$\calG^{[n]}_{e,u}(\oo\cdot\nn;\e,\be_{0},B_{0})$ are
bounded proportionally to $1/(\oo\cdot\nn)^{2}$, the series (\ref{eq:3.15})
converge and solve the range equations (\ref{eq:2.12}). Then,
in the next section, we shall see that the assumption is justified
at least along a curve $(\be_{0}(\e),B_{0}(\e))$ satisfying
also the bifurcation equations (\ref{eq:2.7c}) and (\ref{eq:2.7d}).
We shall not write the dependence on $\e,\be_{0},B_{0}$ unless needed.

\begin{defi}\label{def:4.1} 
We shall say that $\MM$ satisfies property 1 if one has
$$
\Psi_{n+1}(x)\left|{\rm det}\left((\ii x)\uno-\MM^{[n]}(x)
\right)\right|\ge\Psi_{n+1}(x)x^{2}/2,
$$
for all $n\ge -1$.
\end{defi}

\begin{defi}\label{def:4.2}
We shall say that $\MM$ satisfies property 1-$p$ if one has
$$
\Psi_{n+1}(x)\left|{\rm det}\left((\ii x)\uno-\MM^{[n]}(x)
\right)\right|\ge\Psi_{n+1}(x)x^{2}/2.
$$
for $-1\le n<p$.
\end{defi}

\begin{lemma}\label{lem:4.3} 
Assume $\MM$ to satisfy property 1-$p$. Then, for
$0\le n\le p$ and $\e$ small enough, the self-energies are well defined
and one has 
\begin{subequations}
\begin{align}
&|M^{[n]}_{u,e}(x)|\le |\e| \, K_{1} {\rm e}^{-K_{2}2^{m_{n}}}, 
\label{eq:4.1a} \\
&|\partial_{x}^{j}M^{[n]}_{u,e}(x)|\le  |\e| \, C_{j} {\rm e}^{-\ol{C}_{j}2^{m_{n}}},\qquad j=1,2,
\label{eq:4.1b}
\end{align}
\label{eq:4.1}
\end{subequations}
\vskip-.3truecm
\noindent
for some constants $K_{1},K_{2},C_{1},C_{2},\ol{C}_{1}$ and $\ol{C}_{2}$.
\end{lemma}

\prova
We shall prove first (\ref{eq:4.1a}) by induction on $n$.
Let $n\le p$ and $T\in\gotR_{n,u,e}$.
The analyticity of $F,G$ and $\om_{0}$ implies that there exist 
positive constants $F_{1},F_{2},\x$ such that for all $v\in N(T)$ one has 
\begin{equation} \nonumber
|\calF_{v}|
\le F_{1}F_{2}^{s_{v}} {\rm e}^{-\x|\nn_{v}|}. 
\end{equation}
Note that
$$
\prod_{v \in N(T)} {\rm e}^{-\frac14 \x|\nn_{v}|} =
\exp \Big( -\frac14 \x \, K(T) \Big) < \exp \Big( -\frac18 \xi \, 2^{m_{n}} \Big) ,
$$
by Lemma \ref{lem:3.10}.
Moreover by property 1-$p$ and the inductive
hypothesis, one has (for instance)
\begin{equation}\nonumber
\begin{aligned}
|\calG^{[n']}_{\be,\be}(x)|&\le\frac{2}{x^{2}} \Big( |\ii x| +
\big| \MM_{B,B}^{[n'-1]}(x) \big| \Big) \Psi_{n'}(x) \\
&\le \frac{2}{x^{2}}\Big(|x|+P_{1}+|\e|^{2} K_{1}\sum_{q=0}^{n'-1} {\rm e}^{-K_{2}2^{m_{q}}}\Big)
\Psi_{n'}(x) \le \g_{0} \, \al_{m_{n'}}(\oo)^{-2}
\end{aligned}
\end{equation}
for all $0\le n'\le n$ and for a suitable constant $\g_{0}$, where we used
that any renormalised self-energy cluster $T$
on scale $\ge 0$ has at least two nodes and hence $k(T)\ge 2$ and that
there exists $P_{1}\ge 0$ such that $|\MM^{[-1]}_{u,e}|\le P_{1}$ (see
Remark \ref{rmk:3.6}). Of course one can reason analogously
for $\calG^{[n']}_{\be,B}(x)$, $\calG^{[n']}_{B,\be}(x)$ and
$\calG^{[n']}_{\be,B}(x)$, possibly redefining $\g_{0}$.
Hence by Lemmas \ref{lem:3.10} and \ref{lem:3.2} one can bound
\begin{equation*} 
\begin{aligned} 
\prod_{\ell\in L(T)}|\calG^{[n_{\ell}]}_{e_{\ell},u_{\ell}}(\oo\cdot\nn_{\ell})|
&\le  \prod_{q\ge 0}\left(\frac{\g_{0}}{\al_{m_{q}}(\oo) ^{2}}\right)^{\gotN_{q}(T)} 
\le \left(\frac{\g_{0}}{\al_{m_{n_{0}}}(\oo) ^{2}}\right)^{3k(T)-1} \!\!\!\!\!\!\!\!\!\! 
\prod_{q\ge n_{0}+1} \left(\frac{\g_{0}}{\al_{m_{q}}(\oo) ^{2}}\right)^{\gotN_{q}(T)}\\ 
 &\le \left(\frac{\g_{0}}{\al_{m_{n_{0}}}(\oo) ^{2}}\right)^{3k(T)-1} \!\!\!\!\!\!\!\!\!\! 
\prod_{q\ge  n_{0}+1}\left(\frac{\g_{0}^{1/2}}{\al_{m_{q}}(\oo)}\right)^{2^{-(m_{q}-3)} K(T)}\\ 
 &\le D(n_{0})^{3k(T)-1}{\rm exp}(\x(n_{0})K(T)), 
\end{aligned} 
\end{equation*} 
with
$$ 
D(n_{0})=\frac{\g_{0}}{\al_{m_{n_{0}}}(\oo) ^{2}},\qquad 
\x(n_{0})=8\sum_{q\ge n_{0}+1}\frac{1}{2^{m_{q}}}\log\frac{\g_{0}^{1/2}}{\al_{m_{q}} 
(\oo)}. 
$$ 
Then, by Hypothesis \ref{hyp1}, one can choose $n_{0}$ such that 
$\x(n_{0})\le \x/2$. Furthermore, Lemma \ref{lem:3.2} ensures also that
the sum over the other labels
is bounded by a constant to the power $k(T)$ and hence one can bound,
for some positive constants $C$ and $K_{0}$,
\begin{equation} \label{eq:4.2} 
|M^{[n]}_{u,e}(x)|\le
\sum_{T\in\gotR_{n,u,e}}|\e|^{k(T)}|\Val_{T}(x)| \le
\sum_{T\in\gotR_{n,u,e}}|\e|^{k(T)}C^{k(T)} {\rm e}^{-K_{0}K(T)} 
\le\sum_{k\ge 2}|\e|^{k}C^{k} {\rm e}^{-K_{2}2^{m_{n}}}, 
\end{equation}
with $K_{2}=K_{0}/2$, then (\ref{eq:4.1a}) is proved for $\e$ small enough. 
Now we prove (\ref{eq:4.1b}), again by induction on $n$. 
For $n=0$ the bound is obvious.
Assume then (\ref{eq:4.1b}) to hold for 
all $n'<n$. For any $T\in\gotR_{n,u,e}$ such that $\Val_{T}(x)\ne0$ 
one has 
\begin{equation} \label{eq:4.3}
\partial_{x}\Val_{T}(x)= 
\sum_{\ell\in\calP_{T}}\!\!\left(\prod_{v\in N(T)} 
\calF_{v}\!\!\right)\!\! 
\left(\partial_{x}\calG^{[n_{\ell}]}_{e_{\ell},u_{\ell}}(x_{\ell})\!\!\!\!\!\! 
\prod_{\ell'\in L(T)\setminus\{\ell\}} 
\!\!\!\!\!\! 
\calG^{[n_{\ell'}]}_{e_{\ell'},u_{\ell'}}(\oo\cdot\nn_{\ell'})\!\!\right), 
\end{equation} 
where $x_{\ell}=\oo\cdot\nn_{\ell}=x+\oo\cdot\nn_{\ell}^{0}$ and 
\begin{equation} \nonumber 
\begin{aligned} 
\partial_{x}\calG^{[n_{\ell}]} (x_{\ell})= 
\frac{d}{dx}\calG^{[n_{\ell}]} & (\oo\cdot\nn_{\ell}^{0}+x)
=\partial_{x}\Psi_{n_{\ell}}(x_{\ell}) 
\Big((\ii x_{\ell})\uno-\MM^{[n_{\ell}-1]}(x_{\ell})\Big)^{-1} 
\\&\quad
-\Psi_{n_{\ell}}(x_{\ell}) 
\left((\ii x)\uno-\MM^{[n_{\ell}-1]}(x_{\ell})\right)^{-2}
\left(\ii\uno-\partial_{x}\MM^{[n_{\ell}-1]}(x_{\ell}) 
\right). 
\end{aligned} 
\end{equation} 
One has
\begin{equation} \nonumber 
|\partial_{x}\Psi_{n_{\ell}}(x_{\ell})|\le |\partial_{x}\chi_{n_{\ell}-1} 
    (x_{\ell})|+|\partial_{x}\psi_{n_{\ell}}(x_{\ell})|\le 
    \frac{B_{1}}{\al_{m_{n_{\ell}}}(\oo)}, 
\end{equation} 
for some constant $B_{1}$
and, by (\ref{eq:4.1a}), the inductive hypothesis and Hypothesis \ref{hyp1}, 
\begin{equation} \nonumber 
\begin{aligned} 
|\partial_{x}\MM^{[n_{\ell}-1]}_{u,e}(x_{\ell})| 
    &\le \sum_{q=0}^{n_{\ell}-1}|(\partial_{x}\chi_{q}(x_{\ell})) 
    M^{[q]}_{u,e}(x_{\ell})| 
    +\sum_{q=0}^{n_{\ell}-1}|\partial_{x}M^{[q]}_{u,e}(x_{\ell})|\\ 
 &\le |\e|\, B_{1}K_{1}\sum_{q\ge 0}\frac{1}{\al_{m_{q}}(\oo)} 
    {\rm e}^{-K_{2}2^{m_{q}}}+|\e|\, C_{1}\sum_{q\ge 0} {\rm
      e}^{-\ol{C}_{1}2^{m_{q}}} \le |\e|\, B_{2}, 
\end{aligned} 
\end{equation} 
for some constant $B_{2}$. Hence the differentiated propagator
$\partial_{x}\calG^{[n_{\ell}]}_{e_{\ell},u_{\ell}}(x_{\ell})$ can be bounded by
$\g_{1}\al_{m_{n_{\ell}}}(\oo)^{-4}$ for some constant $\g_{1}$.
Possibly redefining the constant $\g_{1}$, also the propagators
of the lines $\ell' \neq\ell$ in (\ref{eq:4.3}) can be bounded
by $\g_{1}\,\al_{m_{n_{\ell'}}}(\oo)^{-4}$, and hence,
at the cost of replacing the previous bound $\g_{0} \al_{m_{n}}(\oo)^{-2}$
for the propagators $\calG^{[n]}(x)$ with $\g_{1} \al_{m_{n}}(\oo)^{-4}$, 
one can reason as in the proof of (\ref{eq:4.1a}) to obtain (\ref{eq:4.1b}) for $j=1$. 
For $j=2$ one can reason analogously.\EP 

\begin{rmk}\label{rmk:4.4}
\emph{
 From the proof of Lemma \ref{lem:4.3} it follows that if $\MM$ satisfies property 1-$p$
the matrices $\MM^{[n]}(x)$ and $\calG^{[n]}(x)$ are well defined
for all $-1\le n\le p$. In particular there exists $\g_{0}>0$ such that
$|\calG^{[n]}_{e,u}(x)|\le \g_{0}\,\al_{m_{n}}(\oo)^{-2}$ for all $0\le n\le p$.
Moreover if $\MM$ satisfies property $1$, the same considerations apply for all $n\ge 0$.
}
\end{rmk}

\begin{lemma}\label{lem:4.10}
Assume $\MM$ to satisfy property 1-$p$. Then one has
\begin{equation} \label{eq:4.4}
\MM^{[n]}(-x)=\MM^{[n]}(x)^{*}
\end{equation}
for all $-1\le n\le p$.
\end{lemma}

\prova
We shall prove the result by induction on $n$. For $n=-1$ the result
is obvious; see Remark \ref{rmk:3.6}.
Assume (\ref{eq:4.4}) to hold up to scale $n-1$. Then, by definition, one has also
$\calG^{[q]}(-x)=\calG^{[q]}(x)^{*}$
for all $0\le q\le n$. For any renormalised self-energy cluster $T$
contributing to
$M^{[n]}(x)$, consider the renormalised self-energy cluster $T'$ obtained
from $T$ by replacing the mode labels $\nn_{v}$ with $-\nn_{v}$ and changing
the sign of the momentum of the entering line. Then
the node factors are changed into their complex conjugated, and this holds
also for the propagators because of the conservation law (\ref{eq:3.1}). Then
$\Val_{T'}(-x)=\Val_{T}(x)^{*}$.
This is enough to prove the assertion.\EP

\begin{lemma}\label{lem:4.11}
Assume $\MM$ to satisfy property 1-$p$. Then, for
$0\le n\le p$ and $\e$ small enough, one has 
\begin{equation} \label{eq:4.5}
\left| M^{[n]}_{u,e}(x) - M^{[n]}_{u,e}(0) - x \, \partial_{x} M^{[n]}_{u,e}(0) \right|
\le  |\e| \,  {K}_{3} {\rm e}^{-\ol{K}_{4} 2^{m_{n}}} x^{2}
\end{equation}
for some constants $K_{3}$ and $K_{4}$.
\end{lemma}

\prova
For $x^{2}>|\e|$ the bound follows trivially from Lemma \ref{lem:4.3}:
thus, we may assume in the following $x^{2}\le |\e|$.
Consider a self-energy cluster $T$ whose value $\Val_{T}(x)$
contributes to $M^{[n]}_{u,e}(x)$ through (\ref{eq:3.10}) and set
$\calA_{T}(x) = \Val_{T}(x) - \Val_{T}(0) - x \, \partial_{x} \Val_{T}(0) $.
Define also
\begin{equation*}
\oln = \min \{ n \in \ZZZ_{+} : K(T) \le 2^{m_{n}} \} .
\end{equation*}
Let us distinguish between the two cases:
(a) $2^{m_{\oln}-1} < K(T) \le 2^{m_{\oln}}$ and (b) $2^{m_{\oln-1}} < K(T) \le 2^{m_{\oln}-1}$.

In case (a), if $\al_{m_{\oln}}(\oo) \le 4|x|$ then one can bound
$|\calA_{T}(x) | \le |\Val_{T}(x) |+|\Val_{T}(0)| +|x \, \partial_{x}\Val_{T}(0)|$.
As soon as $\Psi_{n_{\ell}}(\oo\cdot\nn_{\ell}) \neq 0$ for all
$\ell\in L(T)$, by (the proof of) Lemma \ref{lem:4.3} -- see in
particular (\ref{eq:4.2}) -- each contribution can be bounded as
\begin{equation*} 
\begin{aligned} |\e|^{k(T)} C^{k} {\rm e}^{-K_{0} K(T)} & \le
|\e|^{k(T)} C^{k} {\rm e}^{-(K_{0}/2) K(T)} {\rm e}^{-(K_{0}/2) 2^{m_{\oln}-1}} \\
& \le |\e|^{k(T)} C^{k} {\rm e}^{-(K_{0}/2) K(T)} \al_{m_{\oln}}(\oo)^2\le
16 \, x^{2}|\e|^{k(T)} C^{k} {\rm e}^{-(K_{0}/4) 2^{m_{n}}} .
\end{aligned}
\end{equation*}
If on the contrary $\al_{m_{\ol{n}}}(\oo) > 4|x|$, one can reason as
follows. For any line $\ell\in L(T)$ one has $ |\nn_{\ell}^{0}| \le
K(T) \le 2^{m_{\ol{n}}}$ and hence $|\oo\cdot\nn_{\ell}^{0}| \ge
\al_{m_{\oln}}(\oo)$. Then for all $\tau\in[0,1]$
\begin{equation*}
\frac{5}{4} \left| \oo\cdot\nn_{\ell}^{0} \right| \ge
\left| \oo\cdot\nn_{\ell}^{0} \right| + \left| x \right| \ge
\left| \oo\cdot\nn_{\ell}^{0} + \tau x \right| \ge
\left| \oo\cdot\nn_{\ell}^{0} \right| - \left| x \right| \ge \frac{3}{4}
\left| \oo\cdot\nn_{\ell}^{0} \right| .
\end{equation*}
In particular $(5/4) | \oo\cdot\nn_{\ell}^{0}| \ge |
\oo\cdot\nn_{\ell}| \ge (3/4) | \oo\cdot\nn_{\ell}^{0}|$ and therefore
\begin{equation} \label{eq:4.6}
2 \left| \oo\cdot\nn_{\ell} \right| \ge \left| \oo\cdot\nn_{\ell}^{0}
  + \tau x \right|
\ge \frac{1}{2} \left| \oo\cdot\nn_{\ell} \right| .
\end{equation}
This implies that the sizes of the propagators in $\Val_{T}(tx)$ `do not change too much'
with respect to $\Val_{T}(x)$: in particular (\ref{eq:3.18}) yields the bound (\ref{eq:3.19})
and hence, by Remark \ref{rmk:3.11},
Lemmas \ref{lem:3.9} and \ref{lem:3.10} still hold, so as to obtain
$| \partial_{1}^{2} \Val_{T}(tx) | \le C'(C'')^{k(T)} {\rm e}^{-K_{2} 2^{m_{n}}} $,
where $\partial_{1}$ denotes the derivative with respect to the (only) argument,
for some constants $C'$ and $C''$. Then
\begin{equation} \label{eq:4.7}
\left| \calA_{T}(x) \right| \le \left|
x^{2} \int_{0}^{1} {\rm d} \tau \left(1 - \tau
\right) \partial_{1}^{2} \Val_{T}( \tau x) \right| \le
x^{2} C'(C'')^{k} {\rm e}^{-K_{2} 2^{m_{n}}} ,
\end{equation}
By summing over all possible self-energy values contributing to $M^{[n]}_{u,e}(x)$
the bound (\ref{eq:4.5}) follows.

In case (b), if $\al_{m_{\oln-1}}(\oo) \le 8|x|$ then one can bound
$|\calA_{T}(x) | \le |\Val_{T}(x) |+|\Val_{T}(0)| +|x
\, \partial_{x}\Val_{T}(0)|$ and use that $K(T)> 2^{m_{\oln-1}}$ to obtain
\begin{equation}\nonumber
{\rm e}^{-(K_{0}/2)K(T)} \le {\rm e}^{-(K_{0}/2)2^{m_{\oln-1}}} \le 
\al_{m_{\oln-1}}(\oo)^{2} \le 64 x^{2} .
\end{equation}
If $\al_{m_{\oln-1}}(\oo) > 8|x|$, for any line $\ell\in L(T)$ one has $ |\nn_{\ell}^{0}| \le
K(T) \le 2^{m_{\ol{n}}-1}$ and hence
$$
|\oo\cdot\nn_{\ell}^{0}| \ge
\al_{m_{\oln}-1}(\oo)> \frac{1}{2}\al_{m_{\oln-1}}(\oo) .
$$
Then one can reason as done in case (a) and obtain (\ref{eq:4.6}) for all
$t\in[0,1]$: in turn this yields the bound (\ref{eq:4.7}) and hence
the bound (\ref{eq:4.5}) follows once more.\EP

\begin{rmk}\label{rmk:4.12}
\emph{
 From (\ref{eq:4.1}) and Lemmas \ref{lem:4.10} and \ref{lem:4.11} it follows that
if property 1-$p$ (respectively property 1) is satisfied then
for all $n\le p$ (respectively for all $n\ge-1$) one has
$\MM^{[n]}(x)=\MM^{[n]}(0)+\partial_{x} \MM^{[n]}(0)x+O(\e x^{2})$,
where $\MM^{[n]}(0)$ is a real-valued matrix, while
$\partial_{x}\MM^{[n]}(0)$ is a purely imaginary one.
In particular this implies that if
$\Psi_{n+1}(x)\,|x^{2}-{\rm det}(\MM^{[n]}(0))|\ge\Psi_{n+1}(x)x^{2}/2$
for all $-1\le n<p$ (respectively for all $n\ge -1$) then property
1-$p$ (respectively property $1$) holds.
}
\end{rmk}

The following result will be crucial to check, in the forthcoming Section \ref{sec:5},
that property 1 is satisfied by $\MM$. The proof follows the lines of that
for Lemma 4.8 in \cite{CG2} and it is deferred to Appendix \ref{app:a}
(see (\ref{eq:3.17}) for the definition of $\Phi^{\RR,p}_{\vzero}$ and $\Ga^{\RR,p}_{\vzero}$).

\begin{lemma}\label{lem:4.13}
Assume $\MM$ to satisfy property 1-$p$. Then
\begin{equation} \label{eq:4.8}
\MM^{[p]}(0)=\begin{pmatrix}
\partial_{\be_{0}}\Phi^{\RR,{p}}_{\vzero}+e_{p,\be,\be} &
\partial_{B_{0}}\Phi^{\RR,{p}}_{\vzero}+e_{p,\be,B} \cr
& \cr
\partial_{\be_{0}}\Ga^{\RR,{p}}_{\vzero}+e_{p,B,\be} &
\partial_{B_{0}}\Ga^{\RR,{p}}_{\vzero}+e_{p,B,B}
\end{pmatrix},
\end{equation}
with $|e_{p,u,e}|\le |\e|\, A_{1} e^{-A_{2}2^{m_{p+1}}}$, $u,e=\be,B$,
for suitable positive constants $A_{1}$ and $A_{2}$.
\end{lemma}

\begin{lemma}\label{lem:4.5} 
Assume $\MM$ to satisfy property 1. Then 
the series (\ref{eq:3.15}) and (\ref{eq:3.16}) with the coefficients 
given by (\ref{eq:3.13}) and (\ref{eq:3.14}) respectively, converge for $\e$ small enough. 
\end{lemma} 

\prova 
Let $\theta\in\Theta^{\RR}_{k,\nn,h}$. By Remark \ref{rmk:4.4}
one can bound $|\calG^{[n]}_{e,u}(x)|\le \g_{0} \,
{\al_{m_{n}}(\oo)}^{-2}$ for all $n\ge 0$ and 
hence by Lemma \ref{lem:3.9} one can reason as in the proof of
the bound (\ref{eq:4.1a}) so as to obtain
$$ 
\sum_{\theta\in\Theta^{\RR}_{k,\nn,h}}|\Val(\theta)|\le C_{0} \ol{C}_{0}^{k} {\rm e}^{-\x |\nn|/2}, 
$$ 
for some constants $C_{0}$ and $\ol{C}_{0}$, which is enough to prove the assertion.\EP 

\begin{lemma}\label{lem:4.6} 
Assume $\MM$ to satisfy property 1. Then for $\e$ small enough the function (\ref{eq:3.15}), 
with the coefficients given by (\ref{eq:3.13}), solve the equations (\ref{eq:2.12}). 
\end{lemma} 

\prova 
We shall prove that, the functions $b^{\RR},B^{\RR}$ defined after (\ref{eq:3.16}) 
satisfy the equations of motion (\ref{eq:2.12}), i.e. we shall check that 
$f^{\RR}:=(b^{\RR},B^{\RR})=g\, \Xi(\oo t,f^{\RR})$,
where $g$ is the  pseudo-differential operator with kernel $g(\oo\cdot\nn)= 
1/(\ii\oo\cdot\nn)\uno$ and
$\Xi(\oo t,f^{\RR}):= \big(\om(B^{\RR})+\e F(\oo t,\be^{\RR},B^{\RR}),\e G(\oo t,\be^{\RR},B^{\RR})\big)$.
We can write the Fourier coefficients of $b^{\RR}$ and $B^{\RR}$  as 
\begin{equation*} 
\begin{aligned}
&b^{\RR}_{\nn}=\sum_{n\ge 0}b_{\nn}^{[n]},\qquad 
b_{\nn}^{[n]}=\sum_{k\ge1}\e^{k}\sum_{\theta \in \Theta^{\RR}_{k,\nn,\be}(n)} 
\Val(\theta), \\
&B^{\RR}_{\nn}=\sum_{n\ge 0}B_{\nn}^{[n]},\qquad 
B_{\nn}^{[n]}=\sum_{k\ge1}\e^{k}\sum_{\theta \in \Theta^{\RR}_{k,\nn,B}(n)} 
\Val(\theta),
\end{aligned}
\end{equation*} 
where $\Theta^{\RR}_{k,\nn,h}(n)$ is the subset of $\Theta^{\RR}_{k,\nn,h}$ 
such that $n_{\ell_{\theta}}=n$. Set also
$\ol{\Theta}^{\RR}_{k,\nn}(n):=\ol{\Theta}^{\RR}_{k,\nn,\be}(n)
\times\ol{\Theta}^{\RR}_{k,\nn,B}(n)$ and, for
$\tau=(\theta,\theta')\in\ol{\Theta}^{\RR}_{k,\nn}(n)$,
define $\Val(\tau):=(\Val(\theta),\Val(\theta'))$. 

Using Lemmas \ref{lem:3.4} and \ref{lem:4.5}, in Fourier space one can write
\begin{equation} \nonumber 
\begin{aligned} 
g(\oo\cdot\nn)&[\Xi(\oo t,f^{\RR})]_{\nn} 
= g(\oo\cdot\nn)\sum_{n\ge0}\Psi_{n}(\oo\cdot\nn) 
 [\Xi(\oo t,f^{\RR})]_{\nn} \\ 
&= g(\oo\cdot\nn)\sum_{n\ge0}\Psi_{n}(\oo\cdot\nn) 
(\calG^{[n]}(\oo\cdot\nn))^{-1}\calG^{[n]}(\oo\cdot\nn) 
 [\Xi(\oo t,f^{\RR})]_{\nn} \\ 
&= g(\oo\cdot\nn)\sum_{n\ge0}\left((\ii\oo\cdot\nn)\uno-
\MM^{[n-1]}(\oo\cdot\nn)\right)
\sum_{k\ge1}\e^{k}\sum_{\tau\in\ol{\Theta}^{\RR}_{k,\nn}(n)} 
 \Val(\tau), 
\end{aligned} 
\end{equation} 
where $\ol{\Theta}^{\RR}_{k,\nn}(n)$ differs from $\Theta^{\RR}_{k,\nn}(n)$ 
as it also includes couples of trees where the root line of one or both of
them is the exiting line of a renormalised self-energy cluster.
If we separate such couples from the others, we obtain 
\begin{equation} \nonumber 
\begin{aligned} 
& g(\oo\cdot\nn) [\Xi(\oo t,f^{\RR})]_{\nn} 
= g(\oo\cdot\nn) \Biggl[ \sum_{n\ge0}\left((\ii\oo\cdot\nn)\uno-
\MM^{[n-1]}(\oo\cdot\nn)\right)f_{\nn}^{[n]} \\ 
&\qquad\qquad  + \sum_{n\ge0}\Psi_{n}(\oo\cdot\nn) 
\sum_{p\ge n}\sum_{q=-1}^{n-1}M^{[q]}(\oo\cdot\nn)f_{\nn}^{[p]} + \sum_{n\ge1}\Psi_{n}(\oo\cdot\nn) 
\sum_{p=0}^{n-1}\sum_{q=-1}^{p-1}M^{[q]}(\oo\cdot\nn)  f_{\nn}^{[p]} \Biggr] \\ 
& \quad = g(\oo\cdot\nn) \Biggl[ \sum_{n\ge0}\left((\ii\oo\cdot\nn)\uno-
\MM^{[n-1]}(\oo\cdot\nn)\right)f_{\nn}^{[n]} +
\sum_{p\ge0}\Biggl(\sum_{q=-1}^{p-1} 
M^{[q]}(\oo\cdot\nn)\sum_{n\ge q+1}\Psi_{n}(\oo\cdot\nn) 
\Biggr) f_{\nn}^{[p]} \Biggr] \\ 
& \quad =  g(\oo\cdot\nn) \Biggl[ \sum_{n\ge0}\left((\ii\oo\cdot\nn)\uno-
\MM^{[n-1]}(\oo\cdot\nn\right)f_{\nn}^{[n]}  + \sum_{n\ge0}\Biggl(\sum_{q=-1}^{n-1} 
M^{[q]}(\oo\cdot\nn)\chi_{q}(\oo\cdot\nn)\Biggr) f_{\nn}^{[n]} \Biggr] \\ 
& \quad =  g(\oo\cdot\nn) \Biggl[ \sum_{n\ge0}\left((\ii\oo\cdot\nn)\uno-
\MM^{[n-1]}(\oo\cdot\nn)\right)f_{\nn}^{[n]} + \sum_{n\ge0}\MM^{[n-1]}(\oo\cdot\nn) 
f_{\nn}^{[n]} \Biggr] =\sum_{n\ge0}f_{\nn}^{[n]}=f^{\RR}_{\nn}, 
\end{aligned} 
\end{equation} 
so that the proof is complete.\EP

\begin{rmk}\label{rmk:4.14}
\emph{
 From Lemma \ref{lem:4.13} it follows that, if $\MM$ satisfies property $1$, one can define
\begin{equation*}
\MM^{[\io]}(x):=\lim_{n\to\io}\MM^{[n]}(x),
\end{equation*}
and one has
\begin{equation}\label{eq:4.9}
\MM^{[\io]}(0)=\begin{pmatrix}
\partial_{\be_{0}}\Phi^{\RR}_{\vzero} & \partial_{B_{0}}\Phi^{\RR}_{\vzero} \cr & \cr
\partial_{\be_{0}}\Ga^{\RR}_{\vzero} & \partial_{B_{0}}\Ga^{\RR}_{\vzero}
\end{pmatrix}.
\end{equation}
Note that (\ref{eq:4.9}) is pretty much the same equality provided
by Lemma 4.8 in \cite{CG2}, adapted to the present case.
}
\end{rmk}

\begin{rmk}\label{rmk:4.9}
\emph{If we take the formal expansion of the functions $\Phi_{\vzero}^{\RR}
(\e,\be_{0},B_{0})$, $\Ga_{\vzero}^{\RR}(\e,\be_{0},B_{0})$ and
$\MM^{[\io]}_{u,e}(0;\e,\be_{0},B_{0})$, $u,e\in\{\be,B\}$, we obtain
tree expansions where the self-energy clusters are allowed; see Section 6
where such a situation is discussed for the Hamiltonian case.
Then it is easy to prove the identity (\ref{eq:4.9})
to any perturbation order; in particular, if one expands
$$
\det\Biggl(\sum_{k=0}^{k_{0}-1}\e^{k}
\Biggl[ \MM^{[\io]}\Big(0;\e,\be_{0},\ol{B}_{0}+\sum_{h=1}^{k_{0}-1}
\e^{h}B^{(h)}_{0}+O(\e^{k_{0}})\Big)\Biggr]^{(k)}\Biggr)=
\sum_{k=0}^{k_{0}-1}\e^{k}\de^{(k)}+O(\e^{k_{0}}),
$$
one has $\de^{(k)}=\de^{(k)}(\be_{0})\equiv0$ for all $k=0,\ldots k_{0}-1$,
if the coefficients $B^{(h)}_{0}=B^{(h)}_{\vzero}
(\be_{0})$ are defined as in (\ref{eq:2.9}).
Moreover, for such an expansion, if one writes
\begin{equation}\nonumber
\det\Biggl(\sum_{k=0}^{k_{0}-1}\e^{k}\Biggl[\MM^{[n]}
(0,\be_{0},\ol{B}_{0}+\sum_{h=1}^{k_{0}-1}
\e^{h}B^{(h)}_{0}+O(\e^{k_{0}})\Big)\Biggr]^{(k)}\Biggr) =
\sum_{k=0}^{k_{0}-1}\e^{k}\de^{(k)}_{n}+O(\e^{k_{0}})
\end{equation}
one has
\begin{equation} \nonumber
\Biggl| \sum_{k=0}^{k_{0}-1}\e^{k}\de^{(k)}_{n} \Biggr|\le
A_{1} {\rm e}^{-A_{2}2^{m_{n}}} \end{equation}
for some positive constants $A_{1},A_{2}$. However, under Hypotheses
\ref{hyp1}, \ref{hyp2} and \ref{hyp4}, we are not able to
prove the convergence of the series and we need to introduce some
resummation procedure to give a meaning to the series.
}
\end{rmk}

\begin{lemma}\label{lem:4.7}
Assume $\MM$ to satisfy property $1$. Then there exists $B_{0}=B_{0}(\e,
\be_{0})$, $C^{\io}$ in both $\e, \be_{0}$, such that $B_{0}(\e,\be_{0})\to
\ol{B}_{0}$ for $\e\to0$ and $\Phi_{\vzero}^{\RR}(\e,\be_{0},
B_{0}(\e,\be_{0}))\equiv 0$ for any $\be_{0}$ and $\e$ small enough.
\end{lemma}

\prova
One has $\Phi^{\RR}_{\vzero}(\e;\be_{0},B_{0})=\om_{0}(B_{0})+O(\e)$ and it is
$C^{\io}$ in its arguments because of the assumption that $\MM$ satisfies
property $1$. Then, by Hypothesis \ref{hyp2} one can apply the
implicit function theorem to obtain the result. In particular one has
$$
B_{0}(\e,\be_{0})=\ol{B}_{0}+\sum_{h=1}^{k_{0}}\e^{h}B_{\vzero}^{(h)}(\be_{0})+
O(\e^{k_{0}+1}),
$$
where the coefficients $B^{(h)}_{\vzero}(\be_{0})$ coincide with those defined in (\ref{eq:2.9}).\EP

Given $x_{0}\in\RRR$ and an interval $(a,b) \subset \RRR$ such that
$x_{0}\in(a,b)$, we call half-neighbourhood of $x_{0}$ each of the two
intervals $(a,x_{0})$ and $(x_{0},b)$.

\begin{lemma}\label{lem:4.8}
Assume $\MM$ to satisfy property $1$ and set
$g(\e,\be_{0}):=\Ga^{\RR}_{\vzero}(\e;\be_{0},B_{0}(\e,\be_{0}))$, where
$B_{0}(\e,\be_{0})$ is the $C^{\io}$ function referred to in Lemma
\ref{lem:4.7}. Then there exists a continuous curve $\be_{0}=\be_{0}(\e)$
such that $g(\e,\be_{0}(\e))=0$ and moreover, at least in a suitable
half-neighbourhood
of $\e=0$, one has $\det\left(\MM^{[\io]}(0;\e,\be_{0}(\e),B_{0}(\e,\be_{0} (\e)))\right)\le 0$.
\end{lemma}

\prova
Using the same argument in the proof of Lemma 4.11 in \cite{CG2}, as
property $1$ and Hypothesis \ref{hyp4} imply
$g(\e,\be_{0})=\e^{k_{0}} \left(\Ga_{\vzero}^{(k_{0})}(\be_{0})+O(\e)\right)$
and $\om_{0}'(B_{0}(\e,\be_{0}))$ has the same sign of $\om_{0}'(\ol{B}_{0})$ for
$\e$ small enough, one can
find a continuous curve $\be_{0}=\be_{0} (\e)$ defined at least in a suitable
half-neighbourhood of $\e=0$ such that $g(\e,\be_{0}(\e))\equiv0$ and
$\partial_{\be_{0}}g(\e,\be_{0}(\e))\om_{0}'(B_{0}(\e,\be_{0}(\e)))\ge 0$. Moreover one has
\begin{equation}\nonumber
\partial_{\be_{0}}g(\e,\be_{0})=\partial_{2}\Ga^{\RR}_{\vzero}(\e;\be_{0},
B_{0}(\e,\be_{0}))+\partial_{3}\Ga^{\RR}_{\vzero}(\e;\be_{0},B_{0}
(\e,\be_{0}))\partial_{\be_{0}}B_{0}(\e,\be_{0}),
\end{equation}
where we denoted by $\partial_{j}$ the derivative with respect to the $j$-th argument, so that
\begin{equation}\nonumber
\det\left(\MM^{[\io]}(0;\e,\be_{0}(\e),B_{0}(\e,\be_{0}(\e)))\right)=
-\partial_{\be_{0}}g(\e,\be_{0}(\e))\left(\om_{0}'(B_{0}(\e,\be_{0}(\e))) +
O(\e)\right)
\end{equation}
and then the assertion follows. In particular if $k_{0}$ is even,
the curve $\be_{0}(\e)$ above can be defined in a whole neighbourood of $\e=0$.\EP

By the results above it follows that, if property $1$ is satisfied,
choosing $\be_{0}=\be_{0}(\e)$ and $B_{0}=B_{0}(\e,\be_{0}(\e))$ as above,
the series (\ref{eq:3.15}) solve the equation of motion (\ref{eq:2.1}).

In the forthcoming Section we shall prove that $\MM$ satisfies property 1, at least along
a continuous curve $C(\e):=(\be_{0}(\e),B_{0}(\e,\be_{0}(\e)))$ satisfying
$\Phi_{\vzero}^{\RR}(\e;C(\e))\equiv\Ga_{\vzero}^{\RR}(\e;C(\e))\equiv 0$,
adapting the analogous proof in \cite{CG2}.

\zerarcounters
\section{Convergence of the resummed series: fixing the initial phase}
\label{sec:5}

In this section, we shall complete the proof of Theorem \ref{thm:2.2}
by showing that, under Hypotheses \ref{hyp1}, \ref{hyp2} and \ref{hyp4},
by suitably choosing $\be_{0},B_{0}$, then $\MM$ turns out to satisfy property 1.

Define the $C^{\io}$ non-increasing function $\x$ such that
\begin{equation}\label{eq:5.1} 
\x(x)=\left\{ 
\begin{aligned} 
&1,\quad x\le 1/2,\\ 
&0,\quad x\ge 1, 
\end{aligned}\right. 
\end{equation} 
and set $\x_{-1}(x)=1$ and $\x_{n}(x)=\x(2^{8}x/\al_{m_{n+1}}^{2}(\oo))$ for all
$n\ge0$. Set also
\begin{equation}\label{eq:5.2}
B_{0}(\e,\be_{0},B_{0}'):=\ol{B}_{0}+
\sum_{h=1}^{k_{0}-1}\e^{h}B_{\vzero}^{(h)}(\be_{0})+\e^{k_{0}}B_{0}'
\end{equation}
where the coefficients $B^{(h)}_{\vzero}(\be_{0})$ are defined as in
(\ref{eq:2.9}) and $k_{0}$ is as in Hypothesis \ref{hyp4}. For all $n\ge 0$
we define recursively the \emph{regularised propagators} as
\begin{equation}\label{eq:5.3}
\ol{\calG}^{[n]}=\ol{\calG}^{[n]}(x;\e,\be_{0},B_{0}'):=
\Psi_{n}(x)\left((\ii x)\uno-
\ol{\MM}^{[n-1]}(x;\e,\be_{0},B_{0}')\x_{n-1}
(\De_{n-1})\right)^{-1},
\end{equation}
where
\begin{equation}\label{eq:5.4}
\ol{\MM}^{[n-1]}(x;\e,\be_{0},B_{0}'):
=\sum_{q=-1}^{n-1}\chi_{q}(x)
\ol{M}^{[q]}(x;\e,\be_{0},B_{0}'),
\end{equation}
with the $2\times 2$ matrix $\ol{M}^{[q]}(x;\e,\be_{0},B_{0}')$
defined so as
\begin{equation}\label{eq:5.5}
\ol{M}^{[q]}_{u,e}(x;\e,\be_{0},B_{0}'):=\sum_{T\in\gotR_{q,u,e}}
\e^{k(T)}\ol{\Val}_{T}(x;\e,\be_{0},B_{0}'),
\end{equation}
where
\begin{equation}\label{eq:5.6}
\ol{\Val}_{T}(x;\e,\be_{0},B_{0}'):=
\Biggl(\prod_{v\in N(T)}\widetilde{\calF}_{v}\Biggr)
\Biggl(\prod_{\ell\in L(T)}
\ol{\calG}^{[n_{\ell}]}_{e_{\ell},u_{\ell}} \Biggr),
\end{equation}
with $\widetilde{\calF}_{v}=\calF_{v}(\be_{0},B_{0}(\e,\be_{0},B_{0}'))$ and
$$
\De_{n-1}=\De_{n-1}(\e,\be_{0},B_{0}'):=
D_{n-1}(\e,\be_{0},B_{0}')-\sum_{k=0}^{k_{0}-1}\e^{k}
\left[D_{n-1}(\e,\be_{0},B_{0}')\right]^{(k)},
$$
with
$$
D_{n-1}(\e,\be_{0},B_{0}'):=\det\left(
\ol{\MM}^{[n-1]}(0;\e,\be_{0},B_{0}')\right).
$$

For any $\theta\in\Theta^{\RR}_{k,\nn,h}$, define also,
for all $k\ge0$, $\nn\in\ZZZ^{d}$, $h=\be,B$,
$$
\ol{\Val}(\theta;\e,\be_{0},B_{0}'):=
\Biggl(\prod_{v\in N(T)}\widetilde{\calF}_{v}\Biggr)
\Biggl(\prod_{\ell\in L(T)}\ol{\calG}^{[n_{\ell}]}_{e_{\ell},u_{\ell}} \Biggr) .
$$
Finally, set
$\ol{\MM}:=\{\ol{\MM}^{[n]}(x;\e,\be_{0},B_{0}')\}_{n\ge-1}$ and
$\ol{\MM}^{\x}:=\{\ol{\MM}^{[n]}(x;\e,\be_{0},B_{0}')\x_{n}(\De_{n})\}_{n\ge-1}$.

\begin{lemma}\label{lem:5.1}
For $\e$ small enough, $\ol{\MM}^{\x}$ satisfies property 1.
\end{lemma}

\prova
We shall prove that $\ol{\MM}^{\x}$ satisfies property 1-$p$ for
all $p\ge 0$, by induction on $p$.
For $p=0$ it is obvious if $\e$ is small enough. Assume then that
$\ol{\MM}^{\x}$ satisfies property 1-$p$. Then we can repeat almost word by
word the proofs of Lemmas \ref{lem:4.3} and \ref{lem:4.10}, so as to obtain
\begin{equation}\nonumber
\ol{\MM}^{[p]}(x;\e,\be_{0},B_{0}')=
\ol{\MM}^{[p]}(0;\e,\be_{0},B_{0}')
+ x \, \partial_{x}\ol{\MM}^{[p]}(0;\e,\be_{0},B_{0}')+ x^{2} \!\!
\int_{0}^{1} \!\! {\rm d} t
\left(1 - t\right) \partial_{x}^{2}
\ol{\MM}^{[n]}(tx;\e,\be_{0},B_{0}') ,
\end{equation}
with $\ol{\MM}^{[p]}(0;\e,\be_{0},B_{0}')$ a real-valued
matrix, $\partial_{x}\ol{\MM}^{[p]}(0;\e,\be_{0},B_{0}')$
a purely imaginary one and
\begin{equation} \nonumber
\left| x^{2} \int_{0}^{1} {\rm d} t \left(1 - t\right) \partial_{x}^{2}
\ol{\MM}^{[n]}(tx;\e,\be_{0},B_{0}') \right|  \le C \, |\e| x^{2}
\end{equation}
for some constant $C$, by Lemma \ref{lem:4.11}.
Then we only have to prove that -- see Remark \ref{rmk:4.12} --
$$
\Psi_{p+1}(x)\left|x^{2}-D_{p}(\e,\be_{0},B_{0}')\x_{p}(\De_{p})^{2}\right|
\ge \Psi_{p+1}(x)\frac{x^{2}}{2}.
$$
Note that, by the definition of $\De_{p}$, one has
$$
\sum_{k=0}^{k_{0}-1}\e^{k}
\left[D_{p}(\e,\be_{0},B_{0}')\right]^{(k)}=\sum_{k=0}^{k_{0}-1}
\e^{k}\de_{p}^{(k)}
$$
with the coefficients $\de_{p}^{(k)}$ as in Remark \ref{rmk:4.9},
and hence $\ol{\MM}^{\x}$ satisfies
property 1-$(p+1)$ by the definition of the function $\x_{p}$.\EP

Set
\begin{equation}\label{eq:5.7}
\ol{\MM}^{[\io]}(x;\e,\be_{0},B_{0}'):=\lim_{n\to\io}
\ol{\MM}^{[n]}(x;\e,\be_{0},B_{0}'),
\end{equation}
and define
\begin{equation}\label{eq:5.8}
\ol{\Phi}(\e,\be_{0},B_{0}') :=\sum_{k\ge 0}\e^{k}
\!\!\!\!\! \sum_{\theta\in\Theta^{\RR}_{k,\vzero,\be}}\ol{\Val}
(\theta,\e,\be_{0},B_{0}'),
\qquad \ol{\Ga}(\e,\be_{0},B_{0}'):=\sum_{k\ge 0}\e^{k}
\!\!\!\!\! \sum_{\theta\in\Theta^{\RR}_{k,\vzero,B}}\ol{\Val}
(\theta,\e,\be_{0},B_{0}').
\end{equation}
%

\begin{lemma}\label{lem:5.2}
One has
\begin{equation}\nonumber
[\ol{P}(\e,\be_{0},B_{0}')]^{(k)}=
[P^{\RR}_{\vzero}(\e,\be_{0},B_{0}(\e,\be_{0},B_{0}'))]^{(k)},
\qquad P=\Phi,\Ga,
\end{equation}
for all $k=0,\ldots,k_{0}$.
\end{lemma}

\proof
Set $\Theta_{k,\nn,h}^{\RR(n)}:=\{\theta\in\Theta_{k,\nn,h}^{\RR,n} :
\exists \ell\in L(\theta) \hbox{ such that }n_{\ell}=n\}$ and write
\begin{equation}\nonumber
\ol{P}(\e,\be_{0},B_{0}')=\sum_{k\ge 0}\e^{k}\sum_{n\ge0}
\sum_{\theta\in\Theta^{\RR(n)}_{k,\vzero,h}}\ol{\Val}
(\theta,\e,\be_{0},B_{0}'),
\end{equation}
with $h=\be,B$ for $P=\Phi,\Ga$, respectively, and note that if
$\theta\in\Theta_{k,\nn,h}^{\RR(n)}$ one has
$$
\prod_{v\in N(\theta)}|\widetilde{\calF}_{v}|\le E_{1}^{|N(\theta)|} {\rm e}^{-E_{2}2^{m_{n}}},
$$
for some constants $E_{1},E_{2}$. Moreover one can write formally
$$
\ol{\calG}_{n_{\ell}}(x) = \Psi_{n_{\ell}}(x) g_{n_{\ell}-1}(x) \Bigl( \uno +
\sum_{m\ge1}\Bigl( g_{n_{\ell}-1}(x) \widetilde{\MM}^{[n_{\ell}-1]}
(x;\e,\be_{0},B_{0}')\x_{n_{\ell}-1}(\De_{n_{\ell}-1}) \Bigr)^{m}\Bigr),
$$
with
$$
g_{n_{\ell}-1}(x) = \frac{1}{(\ii x)^2}\begin{pmatrix}
\ii x & \om_{0}'(\ol{B}_0)\xi_{n_{\ell}-1}(\De_{n_{\ell}-1}) \cr 0 & \ii x \end{pmatrix} 
$$
and
$$
\widetilde{\MM}^{[n_{\ell}-1]}(x;\e,\be_{0},B_{0}'):=
\ol{\MM}^{[n_{\ell}-1]}(x;\e,\be_{0},B_{0}')-\begin{pmatrix}
0 & \om_{0}'(\ol{B}_0) \cr
0 & 0 \end{pmatrix} 
= O(\e),
$$
and we can write
$\x_{n_{\ell}-1}(\De_{n_{\ell}-1}) = 1+\x_{n_{\ell}-1}'(\De^{*})
\De_{n_{\ell}-1}$ for some $\De^{*}$, where $\De_{n_{\ell}-1}=O(\e^{k_{0}})$
and
$$
|\x_{n_{\ell}-1}'(\De^{*})|\le \frac{E_3}{\al_{m_{n_{\ell}}}(\oo)^2}\le 
\frac{E_3}{\al_{m_{n}}(\oo)^2},
$$
for some positive constant $E_{3}$ independent of $n$. Hence the assertion
follows.\EP

Introduce the $C^{\io}$ functions $\hat{\Phi}(\e,\be_{0},B_{0})$,
$\widetilde{\Phi}(\e,\be_{0},B_{0})$,
$\hat{\Ga}(\e,\be_{0},B_{0})$ and
$\widetilde{\Ga}(\e,\be_{0},B_{0})$ such that\\
(1) the first $k_{0}$ coefficients of the Taylor expansion in $\e$ of both
$\hat{\Phi}(\e,\be_{0},B_{0}(\e,\be_{0},B_{0}'))$ and
$\widetilde{\Phi}(\e,\be_{0},B_{0}(\e,\be_{0},B_{0}'))$ coincide with those of
$\ol{\Phi}(\e,\be_{0},B_{0}') $,\\
(2) the first $k_{0}$ coefficients of the Taylor expansion in $\e$ of both
$\hat{\Ga}(\e,\be_{0},B_{0}(\e,\be_{0},B_{0}'))$ and
$\widetilde{\Ga}(\e,\be_{0},B_{0}(\e,\be_{0},B_{0}')) $ coincide with
those of $\ol{\Ga}(\e,\be_{0},B_{0}') $,\\
(3) one has
\begin{equation}\label{eq:5.9}
\ol{\MM}^{[\io]}(0;\e,\be_{0},B_{0}')=\begin{pmatrix}
\partial_{2}\hat{\Phi}(\e,\be_{0},B_{0}(\e,\be_{0},B_{0}')) &
\partial_{3}\widetilde{\Phi}(\e,\be_{0},B_{0}(\e,\be_{0},B_{0}')) \cr
& \cr
\partial_{2}\hat{\Ga}(\e,\be_{0},B_{0}(\e,\be_{0},B_{0}')) &
\partial_{3}\widetilde{\Ga}(\e,\be_{0},B_{0}(\e,\be_{0},B_{0}'))
\end{pmatrix},
\end{equation}
where we denoted by $\partial_{j}$ the derivative with respect to the $j$-th
argument.

Define also, for all $n\ge -1$ the $C^{\io}$ functions
$\hat{\Phi}_{n}(\e,\be_{0},B_{0})$,
$\widetilde{\Phi}_{n}(\e,\be_{0},B_{0})$,
$\hat{\Ga}_{n}(\e,\be_{0},B_{0})$ and
$\widetilde{\Ga}_{n}(\e,\be_{0},B_{0})$ such that
\begin{equation}\label{eq:5.10}
\ol{\MM}^{[n]}(0;\e,\be_{0},B_{0}')=\begin{pmatrix}
\partial_{2}\hat{\Phi}_{n}(\e,\be_{0},B_{0}(\e,\be_{0},B_{0}') ) &
\partial_{3}\widetilde{\Phi}_{n}(\e,\be_{0},B_{0}(\e,\be_{0},B_{0}') ) \cr
& \cr
\partial_{2}\hat{\Ga}_{n}(\e,\be_{0},B_{0}(\e,\be_{0},B_{0}') ) &
\partial_{3}\widetilde{\Ga}_{n}(\e,\be_{0},B_{0}(\e,\be_{0},B_{0}') )
\end{pmatrix},
\end{equation}
and
\begin{equation}\label{eq:5.11}
|P_{n}(\e,\be_{0},B_{0}(\e,\be_{0},B_{0}') )-
P(\e,\be_{0},B_{0}(\e,\be_{0},B_{0}') )|\le
P_{H}|\e| \, {\rm e}^{-B_{H}2^{m_{n}}},
\qquad P=\hat{\Phi},\widetilde{\Phi},\hat{\Ga},\widetilde{\Ga},
\end{equation}
for some constants $A_{P},B_{P}$.
Then, by reasoning as in the proofs of Lemmas \ref{lem:4.7} and
\ref{lem:4.8}, we can find $\widetilde{B}_{0}=\widetilde{B}_{0}(\e,\be_{0})$ and
$\widetilde{\be}_{0}=\widetilde{\be}_{0}(\e)$ such that\\
(i) $\hat{\Phi}(\e,\be_{0},{B}_{0}(\e,\be_{0},\widetilde{B}_{0}(\e,\be_{0})))
\equiv 0$ for all $\be_{0}$
and $\e$ small enough,\\
(ii) $\hat{\Ga}(\e,\widetilde{\be}_{0}(\e),{B}_{0}(\e,
\widetilde{\be}_{0}(\e),\widetilde{B}_{0}
(\e,\widetilde{\be}_{0}(\e))))\equiv 0$
for all $\e$ small enough and\\
(iii) $\partial_{3}\widetilde{\Phi}(\e,\be_{0},{B}_{0}
(\e,\be_{0},\widetilde{B}_{0}(\e,\be_{0})))
\partial_{\be_{0}}\hat{\Ga}(\e,\be_{0},{B}_{0}
(\e,\be_{0},\widetilde{B}_{0}
(\e,\be_{0}))) \big|_{\be_{0}=\widetilde{\be}_{0}(\e)}\ge0$,
at least in a suitable half-neighbourhood of $\e=0$.

\begin{lemma}\label{lem:5.3}
Set $\widetilde{C}(\e)=(\widetilde{\be}_{0}(\e),\widetilde{B}_{0}
(\e,\widetilde{\be}_{0}(\e)))$ with $\widetilde{B}_{0}
(\e,\be_{0})$ and $\widetilde{\be}_{0}(\e)$ as above. Then, along
$\widetilde{C}(\e)$ one has $\x_{n}(\De_{n})\equiv1$ for all $n\ge-1$.
\end{lemma}

\prova
We shall prove the result by induction on $n$. For $n=-1$ it is obvious.
Assume then $\x_{p}(\De_{p})\equiv1$ for all
$p=-1,\ldots,n-1$ along $\widetilde{C}(\e)$ and set
$C(\e)=(\widetilde{\be}_{0}(\e),B_{0}(\e,\widetilde{C}(\e)))$. Hence $\ol{\calG}^{[p]}
(x;\e,\widetilde{C}(\e))\equiv\calG^{[p]}(x;\e,C(\e))$ for all $p=0,\ldots,n$
and thence $\ol{\MM}^{[n]}(x;\e,\widetilde{C}(\e))\equiv\MM^{[n]}(x;\e,C(\e))$.
In particular $\MM$ satisfies property 1-$n$ so that, using Lemma
\ref{lem:4.13} one has
\begin{equation}\nonumber
\ol{\MM}^{[n]}(0;\e,\widetilde{C}(\e))=
\begin{pmatrix}
\partial_{2}\Phi^{\RR,n}_{\vzero}(\e,C(\e))+e_{n,\be,\be} &
\partial_{3}\Phi^{\RR,n}_{\vzero}(\e,C(\e))+e_{n,\be,B} \cr
& \cr
\partial_{2}\Ga^{\RR,n}_{\vzero}(\e,C(\e))+e_{n,B,\be} &
\partial_{3}\Ga^{\RR,n}_{\vzero}(\e,C(\e))+e_{n,B,B}
\end{pmatrix},
\end{equation}
with $|e_{n,u,e}|\le |\e| A_{1}e^{-A_{2}2^{m_{n+1}}}$, $u,e=\be,B$.
On the other hand one has
\begin{equation}\nonumber
\begin{aligned}
&\partial_{2}\hat{\Phi}(\e,C(\e))=-\partial_{3}\hat{\Phi}(\e,C(\e))
\left. \partial_{\be_{0}}{B}_{0}(\e,\be_{0},\widetilde{B}_{0}
(\e,\be_{0})) \right|_{\be_{0}=\widetilde{\be}_{0}(\e)}, \\
&\partial_{2}\hat{\Ga}(\e,C(\e))=
\left.\partial_{\be_{0}}\hat{\Ga}(\e,\be_{0},{B}_{0}
(\e,\be_{0},\widetilde{B}_{0}(\e,\be_{0})))\right|_{\be_{0}=\widetilde{\be}_{0}(\e)}
\!\!\!\!\!\!\!\! -\partial_{3}\hat{\Ga}(\e,C(\e)) \left. \partial_{\be_{0}}
{B}_{0}(\e,\be_{0},\widetilde{B}_{0}(\e,
\be_{0})) \right|_{\be_{0}=\widetilde{\be}_{0}(\e)} ,
\end{aligned}
\end{equation}
so that, without writing explicitly the dependence on $(\e,C(\e))$, one has
\begin{equation}\nonumber
\ol{\MM}^{[n]}(0;\e,\widetilde{C}(\e))=
\begin{pmatrix}
-\partial_{3}\Phi^{\RR,n}_{\vzero}
\partial_{\be_{0}}{B}_{0}
+\gamma_{n} &
\partial_{3}\Phi^{\RR,n}_{\vzero}+e_{n,\be,B} \cr
& \cr
\partial_{\be_{0}}\Ga^{\RR,n}_{\vzero}
-\partial_{3}\Ga^{\RR,n}_{\vzero}\partial_{\be_{0}}
{B}_{0}+\gamma'_{n} &
\partial_{3}\Ga^{\RR,n}_{\vzero}+e_{n,B,B}
\end{pmatrix},
\end{equation}
with $|\gamma_{n}|,|\gamma_{n}'|\le |\e|\, C_{1}e^{-C_{2}2^{m_{n+1}}}$
for some $C_{1},C_{2}$. Hence
\begin{equation}\nonumber
\De_{n}=-\partial_{\be_{0}}\Ga^{\RR,n}_{\vzero}
\partial_{3}\Phi^{\RR,n}_{\vzero}+c_{n}
=-\partial_{\be_{0}}\hat{\Ga}_{n}\partial_{3}\widetilde{\Phi}_{n}+c'_{n}
=-\partial_{\be_{0}}\Ga\partial_{3}\widetilde{\Phi}+c''_{n}\le c''_{n},
\end{equation}
with  $|c_{n}|,|c'_{n}|,|c''_{n}|\le |\e|\, D_{1}e^{-D_{2}2^{m_{n+1}}}$
for some constants $D_{1}$ and $D_{2}$, so that
the assertion follows by the definition of $\xi_{n}$.\EP

\begin{lemma}\label{lem:5.4}
Let $\widetilde{C}(\e)$ be as in Lemma \ref{lem:5.3} and set
$C(\e)=(\widetilde{\be}_{0}(\e),B_{0}(\e,\widetilde{C}(\e)))$.
One can choose
the functions $\hat{\Phi},\widetilde{\Phi},\hat{\Ga},\widetilde{\Ga}$ such that
$\hat{\Phi}(\e,C(\e))=\widetilde{\Phi}(\e,C(\e))=\Phi^{\RR}_{\vzero}(\e,C(\e))
\equiv0$ and $\hat{\Ga}(\e,C(\e))=\widetilde{\Ga}(\e,C(\e))=\Ga^{\RR}_{\vzero}
(\e,C(\e))\equiv 0$. In particular $(\be(t,\e),B(t,\e))=C(\e)+
(b^{\RR}(t;\e,C(\e)),B^{\RR}(t;\e,C(\e)))$ defined in (\ref{eq:3.15}) solves
the equation of motion (\ref{eq:2.1})
\end{lemma}

\prova
It follows from the results above. Indeed, for any $\hat{\Phi}$, $\hat{\Ga}$
there is a curve $C(\e)$ along which $\MM=\ol{\MM}=\ol{\MM}^{\xi}$ (hence
$\MM$ satisfies property $1$) and $\hat{\Phi}(\e,C(\e))=\hat{\Ga}
(\e,C(\e))\equiv 0$. By Remark \ref{rmk:4.14} also
$\Phi^{\RR}_{\vzero}$ and $\Ga^{\RR}_{\vzero}$ are among the primitives of
$\MM^{[\io]}_{\be,\be}$ and $\MM^{[\io]}_{B,\be}$ respectively, and then the
assertion follows.
\EP

\begin{rmk}\label{rmk:5.5}
\emph{Note that without  Remark \ref{rmk:4.14} we were able to prove only
the existence of curves on which the solution of the range equations
is well-defined. On the other hand Remark \ref{rmk:4.14} (which follows from
Lemma \ref{lem:4.13}) guarantees that the solution of the bifurcation]
equations is one of such curves.}
\end{rmk}

Lemma \ref{lem:5.4} completes the proof of Theorem \ref{thm:2.2}: indeed the
function $(\be(t,\e),B(t,\e))$ is a quasi-periodic solution to (\ref{eq:2.1}) with
frequency vector $\oo$ and, by construction, it reduces to $(\ol{\be}_{0},\ol{B}_{0})$
as $\e$ tends to $0$.

\zerarcounters
\section{The Hamiltonian case}
\label{sec:6}

In this section we prove Theorem \ref{thm:2.3}. Consider (\ref{eq:2.1})
of the form (\ref{eq:2.15}), i.e. with $F=\partial_{B}f$ and $G=-\partial_{\be}f$,
where $f$ is the function appearing in (\ref{eq:2.14}).
We look for a solution $(\be(t),B(t))$ which can be formally written as in (\ref{eq:2.8}),
with the coefficients given by (\ref{eq:2.9}). If there exists $k_{0}\ge1$ such that all the
coefficients $\Ga_{\vzero}^{(k)}(\be_{0})$ vanish identically for all
$0\le k \le k_{0}-1$, while $\Ga_{\vzero}^{(k_{0})}(\be_{0})$ is not identically zero,
we can solve the equations of motion up to order $k_{0}$ without fixing
the parameter $\be_{0}$. Moreover one has $\Ga^{(k_{0})}(\be_{0})=\partial_{\be_{0}}g^{(k_{0})}(\be_{0})$, with
$$
g^{(k_{0})}(\be_{0}):=[\ol{B}\,\dot{\ol{b}}
]^{(k_{0})}_{\vzero}-[h_{0}(\ol{B}_{0}+\ol{B}+B^{(k_{0})})]^{(k_{0})}_{\vzero}
-[f(\oo t,\be_{0}+\ol{b},\ol{B}_{0}+\ol{B})]^{(k_{0}-1)}_{\vzero},
$$
because, if we denote
\begin{equation*}
\ol{b}=\sum_{k=1}^{k_{0}-1}b^{(k)},\qquad
\ol{B}=\sum_{k=1}^{k_{0}-1}B^{(k)} ,
\end{equation*}
one has
\begin{equation*}
\begin{aligned}
\partial_{\be_{0}}
[f(\oo t,\be_{0}+\ol{b},\ol{B}_{0}+\ol{B})]^{(k_{0}-1)}_{\vzero} &=
[\partial_{\be}f(\oo t,\be_{0}+\ol{b},\ol{B}_{0}+\ol{B})(1+
\partial_{\be_{0}}\ol{b})]^{(k_{0}-1)}_{\vzero}\\
&\qquad+
[\partial_{B}f(\oo t,\be_{0}+\ol{b},\ol{B}_{0}+\ol{B})
\partial_{\be_{0}}\ol{B}]^{(k_{0}-1)}_{\vzero} \\
&=-\Ga^{(k_{0})}_{\vzero}-[\dot{\ol{B}}\partial_{\be_{0}}\ol{b}]^{(k_{0})}_{\vzero}
+[\dot{\ol{b}}\partial_{\be_{0}}\ol{B}]^{(k_{0})}_{\vzero}\\
&\qquad-
[\om_{0}(\ol{B}_{0}+\ol{B}+B^{(k_{0})})\partial_{\be_{0}}(\ol{B}+
B^{(k_{0})})]^{(k_{0})}_{\vzero}\\
&=-\Ga^{(k_{0})}_{\vzero}+\partial_{\be_{0}}[\ol{B}\,\dot{\ol{b}}
]^{(k_{0})}_{\vzero}
-\partial_{\be_{0}}[h_{0}(\ol{B}_{0}+\ol{B}+B^{(k_{0})})]^{(k_{0})}_{\vzero}.
\end{aligned}
\end{equation*}
Since $g^{(k_{0})}$
is analytic and periodic, it has at least one maximum $\be_{0}'$ and
one minimum $\be_{0}''$.  Then Hypothesis \ref{hyp4} automatically holds.
Indeed, if $\e^{k_{0}}\om'_{0}(\ol{B}_{0})>0$ one can choose $\ol{\be}_{0}=
\be_{0}''$, while if $\e^{k_{0}}\om'_{0}(\ol{B}_{0})<0$ one can choose $\ol{\be}_{0}=
\be_{0}'$ and hence in both cases Hypothesis \ref{hyp4} is satisfied.
Therefore the existence of a quasi-periodic solution with frequency vector $\oo$
follows from Theorem \ref{thm:2.2}.

Assume now that $\Ga_{\vzero}^{(k)}\equiv 0$ for all $k\ge 0$. We
shall prove the following result, which, together with the argument
given above, implies Theorem \ref{thm:2.3}.

\begin{prop}\label{prop:6.1}
Consider the system (\ref{eq:2.15}) and assume Hypotheses \ref{hyp1}
and \ref{hyp2} to be satisfied. Assume also that
$\Ga_{\vzero}^{(k)}\equiv 0$ for all $k\ge 0$.
Then for $\e$ small enough there exists a resonant maximal torus run
with frequency vector $\oo$, which can be parameterised as
$\be=\be_{0}+\be_{\e}(\pps,\be_{0})$,
$B=\ol{B}_{0}+B_{\e}(\pps,\be_{0})$, $\aaa=\pps$, $\AAA=\AAA_{\e}(\pps,\be_{0})$,
with $(\pps,\be_{0})\in\TTT^{d+1}$, where the functions $\be_{\e}$,
$B_{\e}$ and $\AAA_{\e}$ are analytic in $\e$, as well as in $\pps$ and $\be_{0}$,
and are at least of order $\e$.
The time evolution along the torus is given by $(\pps,\be_{0})\to(\pps +\oo t,\be_{0})$.
\end{prop}

As we shall see, no resummation is needed in such a case and the formal
expansion (\ref{eq:2.8}) turns out to converge, i.e. the solution
is analytic in the perturbation parameter. However we still need a multiscale decomposition
of the propagators to show that the small divisors `do not accumulate too much'.
To this aim, we slightly change some of the definitions in Section \ref{sec:3} as follows.

First of all, when defining the labelled trees, the following changes are
made. We associate with each node $v$ a \emph{mode label} $\nn_{v}\in\ZZZ^{d}$, a
\emph{component label} $h_{v}\in\{\be,B\}$ and an \emph{order label} $k_{v}\in\{0,1\}$
with the constraint that $k_{v}=1$ if $\nn_{v}\ne\vzero$.
With each line $\ell=\ell_{v}$, $\ell\ne\ell_{\theta}$, we associate a
\emph{component label} $h_{\ell}\in\{\be,B\}$, with the constraint that $h_{\ell_{v}} =h_{v}$,
and a \emph{momentum label} $\nn_{\ell}\in \ZZZ^{d}$,
with the constraint that $\nn_{\ell}\ne \vzero$ if $h_{\ell}=\be$.
We associate with the root line $\ell_{\theta}$ a component label
$h_{\ell_{\theta}}\in\{\be,B,\Phi,\Ga\}$ and a
momentum label $\nn_{\ell_{\theta}}\in\ZZZ^{d}$ with the following constraints.
Call $v_{0}$ the node which $\ell_{\theta}$ exists: then (i)
$h_{\ell_{\theta}}=B,\Ga$ if $h_{v_{0}}=B$, while $h_{\ell_{\theta}}=\be,\Phi$
if $h_{v_{0}}=\be$ and (ii) $\nn_{\ell_{\theta}}\ne\vzero$ for
$h_{\ell_{\theta}}=\be$, while $\nn_{\ell_{\theta}}=\vzero$ for
$h_{\ell_{\theta}}=\Phi,\Ga$. We require $k_{v_{0}}=1$ if
$\ell=\ell_{v_{0}}$ is such that either $\ell=\ell_{\theta}$ and
$h_{\ell}=\Ga$ or $h_{\ell}=B$ and $\nn_{\ell}\ne\vzero$.
Denote by $p_{v}$ and $q_{v}$ the numbers of lines with component label
$\be$ and $B$, respectively, entering the node $v$ and set $s_{v}=p_{v}+q_{v}$: 
if $k_{v}=0$ for some $v \in N(\theta)$, then we impose $p_{v}=0$ and
$q_{v}\ge1$ if $h_{v}=\be$, while we impose $p_{v}=0$ and $q_{v}\ge2$ if $h_{v}=B$.
Finally with each line $\ell\in L(\theta)$ we associate
a \emph{scale label} $n_{\ell}\in\ZZZ_{+}\cup\{-1\}$, with the constraint that
$n_{\ell}\ge 0$ if $\nn_{\ell}\neq\vzero$ and $n_{\ell}=-1$ if $\nn_{\ell}=\vzero$. 

We do not change the definition of cluster, while
it is more convenient to change slightly the definition of self-energy cluster:
a cluster $T$ on scale $n$ is a \emph{self-energy cluster} if
(i) it has only one entering line $\ell_{T}'$ and one exiting line $\ell_{T}$,
(ii) one has $\nn_{\ell_{T}}=\nn_{\ell_{T}'}$,
(iii) either $n=-1$ and $\calP_{T}=\emptyset$ or $n\ge0$ and one has
$n_{\ell}\ge0$ and $\nn_{\ell}^{0} \neq \vzero$ for all $\ell\in\calP_{T}$.
We shall say that a subgraph $T$ constituted by only one node
$v$ with $\nn_{v}=\vzero$ and $s_{v}=1$,
is also a self-energy cluster on scale $-1$.

We denote by $\Theta_{k,\nn,h}$ the set of trees with order $k$, total momentum
$\nn$ and total component $h$ and by $\gotS^{k}_{n,u,e}$ the set of
self-energy clusters with order $k$, scale $n$ and such that
$h_{\ell'_{T}}=e$ and $h_{\ell_{T}}=u$, with $e,u\in\{\be,B\}$. Note that self-energy clusters
are allowed both in $\Theta_{k,\nn,h}$ and in $\gotS^{k}_{n,u,e}$;
in other words we are considering trees and subgraphs which are not renormalised.

For any tree $\theta$ or any subgraph $S$ of $\theta$
we define their (non-renormalised)
value as in (\ref{eq:3.12}), but with the node factors defined as
\begin{equation}\label{eq:6.1}
\calF_{v}= \left\{ \begin{aligned}
& \frac{1}{q_{v}!}\partial_{1}^{q_{v}}\om_{0}(\ol{B}_{0}) ,
& \qquad k_{v}=0, & & & \\
& \frac{1}{p_{v}!q_{v}!}\partial_{\be}^{p_{v}}\partial_{B}^{q_{v}+1} f_{\nn_{v}}(\be_{0},\ol{B}_{0}) ,
& \qquad k_{v}=1, & \quad h_{{v}}=\be, & & \\
& -\frac{1}{p_{v}!q_{v}!}\partial_{\be}^{p_{v}+1}\partial_{B}^{q_{v}} f_{\nn_{v}}(\be_{0},\ol{B}_{0}) ,
& \qquad k_{v}=1, & \quad h_{{v}}=B, & \; \nn_{\ell_{v}}\ne\vzero, & \\
& \frac{1}{p_{v}!q_{v}!}\partial_{\be}^{p_{v}}\partial_{B}^{q_{v}+1} f_{\nn_{v}}(\be_{0},\ol{B}_{0}) ,
& \qquad k_{v}=1, & \quad h_{{v}}=B, & \; \nn_{\ell_{v}}=\vzero, & \quad h_{\ell_{v}}=B,\\
& -\frac{1}{p_{v}!q_{v}!}\partial_{\be}^{p_{v}+1}\partial_{B}^{q_{v}} f_{\nn_{v}}(\be_{0},\ol{B}_{0}) ,
& \qquad k_{v}=1, & \quad h_{{v}}=B, & \; \nn_{\ell_{v}}=\vzero, & \quad h_{\ell_{v}}=\Ga,
\end{aligned} \right.
\end{equation}
and the propagators defined as
\begin{equation}\label{eq:6.2}
\calG_{n_{\ell}}(\oo\cdot\nn_{\ell}):=\left\{
\begin{aligned}
&\frac{\Psi_{n_{\ell}}(\oo\cdot\nn_{\ell})}{\ii\oo\cdot\nn_{\ell}},
& \qquad n_{\ell}\ge0, \hskip.3truecm  & \\
& - \frac{1}{\om'_{0}(\ol{B}_{0})},
& \qquad n_{\ell}=-1, & \quad h_{\ell}=B, \\
&1,
& \qquad n_{\ell}=-1, & \quad h_{\ell}=\Ga,\Phi.
\end{aligned}
\right.
\end{equation}
%

\begin{lemma}\label{lem:6.1b}
Let $T$ be a subgraph of any tree $\theta$. Then one has $|N(T)|\le 4k(T)-2$.
\end{lemma}

\prova
We shall prove the result by induction on $k=k(T)$.
For $k=1$ the bound is trivially satisfied as a direct check shows.
Assume then the bound to hold for all $k'<k$. Call $v$ the node
which $\ell_{T}$ (possibly $\ell_{\theta}$) exits, $\ell_{1},\ldots,
\ell_{s_{v}}$ the lines
entering $v$ and $T_{1},\ldots,T_{s_{v}}$ the subgraphs of $T$ with
exiting lines $\ell_{1},\ldots,\ell_{s_{v}}$.
If $k_{v}=1$ then by the inductive hypothesis one has
$$
|N(T)|=1+\sum_{i=1}^{s_{v}}|N(T_{i})|\le 1+4(k-1)-2s_{v}\le 4k-2.
$$
If $k_{v}=0$ and $h_{v}=B$ then one has $s_{v}=q_{v}\ge2$ and hence
$$
|N(T)|=1+\sum_{i=1}^{q_{v}}|N(T_{i})|\le 1+4k-2q_{v}\le 4k-2 .
$$
If $k_{v}=0$ and $h_{v}=\be$, then if $q_{v}\ge2$ one can reason as in the
previous case. Otherwise, the line $\ell=\ell_{w}$ entering $v$ is such that
$h_{\ell}=B$, so that either $k_{w}=1$ or $k_{w}=0$ and $q_{w}\ge2$.
call $\ell'_{1},\ldots,\ell'_{s_{w}}$ the lines
entering $w$ and $T'_{1},\ldots,T'_{s_{v}}$ the subgraphs of $T$ with
exiting lines $\ell'_{1},\ldots,\ell'_{s_{w}}$.
In the first case one has
$$
|N(T)|=2+\sum_{i=1}^{s_{w}}|N(T'_{i})|\le 2+4(k-1)-2s_{w}\le 4k-2,
$$
while in the second case one has
$$
|N(T)|=2+\sum_{i=1}^{q_{w}}|N(T'_{i})|\le 2+4k-2q_{w}\le 4k-2 .
$$
Therefore the bound follows.\EP

 From now on we shall not write explicitly the dependence on $\be_{0}$
and $\ol{B}_{0}$ to lighten the notations.
If $T$ is a self-energy cluster we can (and shall) write
$\Val(T)=\Val_{T}(\oo\cdot\nn_{\ell'_{T}})$ to stress the dependence on
$\oo\cdot\nn_{\ell'_{T}}$; see also Remark \ref{rmk:3.7}.
For all $k\ge1$, define
\begin{subequations}
\begin{align}
&b_{\nn}^{(k)}:=\sum_{\theta\in\Theta_{k,\nn,\be}} \Val(\theta),
\quad \nn\in\ZZZ^{d}_{*}, \qquad
B_{\nn}^{(k)}:= \sum_{\theta\in\Theta_{k,\nn,B}} \Val(\theta),
\quad \nn\in\ZZZ^{d}, 
\label{eq:6.3a} \\
&M^{(k)}_{u,e}(x,n):= \!\!\!\! \sum_{T\in\gotS^{k}_{n,u,e}} \!\!\!\! \Val_{T}(x),
\quad \MM^{(k)}_{u,e}(x,n):= \!\! \sum_{p=-1}^{n}M^{(k)}_{u,e}(x,p) , \quad
\MM^{(k)}_{u,e}(x):=\lim_{n\to\io}\MM^{(k)}_{u,e}(x,n).
\label{eq:6.3b} 
\end{align}
\label{eq:6.3}
\end{subequations}
\vskip-.3truecm
\noindent  
Note that the coefficients (\ref{eq:6.3a}) coincide with those in
(\ref{eq:2.9}). Set also
\begin{equation}\label{eq:6.3bis}
\Phi^{(k)}_{\vzero}:= \sum_{\theta\in\Theta_{k,\vzero,\Phi}} \Val(\theta), \qquad
\Ga^{(k)}_{\vzero}:= \sum_{\theta\in\Theta_{k,\vzero,\Ga}} \Val(\theta), 
\end{equation}
and note that $\Phi^{(k)}_{\vzero}=[\om_{0}(B(t))+\e\partial_{B}f(\oo t,\be(t),B(t))]^{(k)}_{\vzero}$
and $\Ga^{(k)}_{\vzero}=[-\e\partial_{\be}f(\oo t,\be(t),B(t))]^{(k)}_{\vzero}$.

\begin{rmk}\label{rmk:6.2}
\emph{
One has
\begin{equation}\nonumber
\begin{aligned}
&\MM_{\be,\be}^{(0)}(x,n)=\MM_{B,B}^{(0)}(x,n)=\MM_{B,\be}^{(0)}(x,n)=0, 
\qquad \MM_{\be,B}^{(0)}(x,n)=\om_{0}'(\ol{B}_{0}), \\
& \MM_{\be,\be}^{(1)}(x,n)=\partial_{\be}\partial_{B}f_{\vzero}=-\MM_{B,B}^{(1)}(x,n),
\end{aligned}
\end{equation}
for all $n\ge-1$ and all $x\in\RRR$.}
\end{rmk}

We shall say that a line $\ell$ is \emph{resonant} if there exist two
self-energy clusters $T$, $T'$, such that $\ell_{T'}=\ell=\ell'_{T}$,
otherwise $\ell$ is \emph{non-resonant}. Given any subgraph $S$ of any tree
$\theta$, we denote by $\gotN^{*}_{n}(S)$ the number of non-resonant lines
on scale $\ge n$ in $S$. Define also, for any line $\ell \in \theta$, the
\emph{minimum scale} of $\ell$ as
\begin{equation*} 
\ze_{\ell} := \min\{ n \in \ZZZ_{+} : \Psi_{n} (\oo\cdot\nn_{\ell}) \neq 0 \} 
\end{equation*}
and denote by $\gotN^{\bullet}_{n}(S)$ the number of non-resonant lines
$\ell\in L(S)$ such that $\ze_{\ell} \ge n$. If $\Val(S)\neq0$,
for each line $\ell\in L(S)$ either $n_{\ell}=\ze_{\ell}$ or $n_{\ell}=\ze_{\ell}+1$.
Then one can prove the following results.

\begin{lemma}\label{lem:6.3}
For all $h\in\{\be,B,\Phi,\Ga\}$, $\nn\in\ZZZ^{d}$, $k\ge1$ and for any
$\theta\in\Theta_{k,\nn,h}$ with $\Val(\theta)\ne 0$, one has
$\gotN_{n}^{\bullet}(\theta)\le 2^{-(m_{n}-3)}K(\theta)$ for all $n\ge0$.
\end{lemma}

\begin{lemma}\label{lem:6.4}
For all $e,u\in\{\be,B\}$, $n \ge 0$, $k\ge1$ and for any $T\in\gotS^{k}_{n,u,e}$
with $\Val_{T}(x)\ne 0$, one has $K(T)>2^{m_{n}-1}$ and
$\gotN_{p}^{\bullet}(T)\le 2^{-(m_{p}-3)}K(T)$ for all $0\le p\le n$.
\end{lemma}

The proofs of the two results above follow the lines of those for
Lemmas \ref{lem:3.9} and \ref{lem:3.10}, and are given in Appendix \ref{app:b}.

\begin{lemma}\label{lem:6.6b}
Let $k\ge 1$, $\nn\in\ZZZ^{d}$, $h\in\{\be,B,\Phi,\Ga\}$, $u,e\in\{\be,B\}$ and $n\ge 0$ arbitrarily fixed.
For any tree $\theta\in\Theta_{k,\nn,h}$ and any self-energy cluster
$T\in\gotS^{k}_{n,u,e}$ denote by $L_{NR}(\theta)$ and $L_{NR}(T)$ the sets
of non-resonant lines in $\theta$ and $T$, respectively, and set
\begin{equation}\nonumber
\Val_{NR}(\theta):=\Biggl(\prod_{v\in N(\theta)} \!\!\!\! \calF_{v}\Biggr)
\Biggl(\prod_{\ell\in L_{NR}(\theta)} \!\!\!\!\!\!\calG_{n_{\ell}}(\oo\cdot\nn_{\ell}) \Biggr),
\quad
\Val_{T,NR}(\oo\cdot\nn_{\ell_{T}'}):=\Biggl(\prod_{v\in N(T)} \!\!\!\! \calF_{v}\Biggr)
\Biggl(\prod_{\ell\in L_{NR}(T)} \!\!\!\!\!\! \calG_{n_{\ell}}(\oo\cdot\nn_{\ell}) \Biggr),
\end{equation}
Then
\begin{subequations}
\begin{align}
&|\Val_{NR}(\theta)|
\le c_{1}^{k}e^{-\xi|\nn|/2},
\label{eq:6.4bis.a}\\
&|\Val_{T,NR}(\oo\cdot\nn_{\ell_{T}'})|\le c_{2}^{k}e^{-\xi K(T)/2},
\label{eq:6.4bis.b}
\end{align}
\label{eq:6.4bis}
\end{subequations}
\vskip-.3truecm
\noindent
for some positive constants $c_{1},c_{2}$.
\end{lemma}

\prova
The proof follows the lines of the proof of the bound
\eqref{eq:4.1a} in Lemma \ref{lem:4.3}.
Indeed by Lemma \ref{lem:6.1b} and the analyticity of $f$ and $\om_{0}$ one
has
$$
\prod_{v\in N(\theta)}|\calF_{v}|\le AB^{k}e^{-\xi K(\theta)},
$$
for some positive constants $A,B$, while
\begin{equation}\nonumber
\begin{aligned}
\prod_{\ell\in L_{NR}(\theta)}\left| \calG_{n_{\ell}}(\oo\cdot\nn_{\ell}) \right|&\le
\prod_{n\ge0}\left(\frac{16}{\al_{m_{n}}(\oo)}\right)^{\gotN^{\bullet}_{n}(\theta)}
\le \left(\frac{16}{\al_{m_{n_{0}}}(\oo)}\right)^{4k-2} 
\!\!\!\!\!\!\!\!\!\! 
\prod_{n\ge n_{0}+1}
\left(\frac{16}{\al_{m_{n}}(\oo)}\right)^{\gotN_{n}^{\bullet}(\theta)}\\ 
 &\le \left(\frac{16}{\al_{m_{n_{0}}}(\oo)}\right)^{4k-2} \!\!\!\!\!\!\!\!\!\! 
\prod_{n\ge  n_{0}+1}\left(\frac{16}{\al_{m_{q}}(\oo)}
\right)^{2^{-(m_{n}-3)} K(\theta)} \!\!\!\!\!\!\!\!\!\! \\
& \le D(n_{0})^{4k-2}{\rm exp}(\x(n_{0})K(\theta)), 
\end{aligned} 
\end{equation}
with 
$$ 
D(n_{0})=\frac{16}{\al_{m_{n_{0}}}(\oo)},\qquad 
\x(n_{0})=8\sum_{n\ge n_{0}+1}\frac{1}{2^{m_{n}}}\log
\frac{16}{\al_{m_{n}}(\oo)}. 
$$ 
Then, by Hypothesis \ref{hyp2}, one can choose $n_{0}$ such that 
$\x(n_{0})\le \x/2$, so that \eqref{eq:6.4bis.a} follows, recalling
$K(\theta)\ge |\nn|$.
To obtain \eqref{eq:6.4bis.b} one can reason in the same way, simply
with $T$ playing the role of $\theta$.
\EP

\begin{rmk} \label{rmk:6.6c}
\emph{
By using Remark \ref{rmk:3.11} one can show that also
$\partial_{x}^{j}\Val_{T,NR}(\tau x)$ admits the same bound
as $\Val_{T,NR}(x)$ in (\ref{eq:6.4bis.b}) for $j=0,1,2$ and $\tau\in[0,1]$,
possibly with a different constant $c_{2}$.
This will be used later on (in Appendix \ref{app:c}).
}
\end{rmk}

In the light of the results above, although the propagators are bounded
proportionally to $1/|x|$ (since no resummation is performed),
in principle one can have accumulation of small divisors
because of the presence of resonant lines. Therefore one needs a `gain factor'
proportional to $\oo\cdot\nn_{\ell}$ for each resonant line $\ell$, in order to prove
the convergence of the power series (\ref{eq:2.8}). One could also envisage
performing the resummation and exploiting cancellations between the
self-energies, but in practice this would make the analysis more complicated.

\begin{lemma}\label{lem:6.5}
Assume that $\Ga_{\vzero}^{(k)}\equiv 0$ for all $k\ge 0$. Then
for all  $k\ge 1$ one has
\begin{subequations}
\begin{align}
&\MM^{(k)}_{B,\be}(0)+\sum_{k_{1}+k_{2}=k}\MM^{(k_{1})}_{B,B}(0)\partial_{\be_{0}} B_{\vzero}^{(k_{2})}=0,
\label{eq:6.4a} \\
&\MM^{(k)}_{\be,\be}(0)+\sum_{k_{1}+k_{2}=k}\MM^{(k_{1})}_{\be,B}(0)\partial_{\be_{0}} B_{\vzero}^{(k_{2})}=0,
\label{eq:6.4b} \\
&\MM^{(k)}_{\be,\be}(0)=-\MM^{(k)}_{B,B}(0).
\label{eq:6.4c}
\end{align}
\label{eq:6.4}
\end{subequations}
\vskip-.3truecm
\end{lemma}

\vskip-0.5truecm
\prova
Both (\ref{eq:6.4a}) and (\ref{eq:6.4b}) follow from the fact that
(see also Remark \ref{rmk:4.9})
\begin{equation*}
\begin{aligned}
&\partial_{\be_{0}}\Ga_{\vzero}^{(k)}=
\MM^{(k)}_{B,\be}(0)+\sum_{k_{1}+k_{2}=k}\MM^{(k_{1})}_{B,B}(0)\partial_{\be_{0}} B_{\vzero}^{(k_{2})}, \\
& \partial_{\be_{0}}\Phi_{\vzero}^{(k)}=
\MM^{(k)}_{\be,\be}(0)+\sum_{k_{1}+k_{2}=k}\MM^{(k_{1})}_{\be,B}(0)\partial_{\be_{0}} B_{\vzero}^{(k_{2})} ,
\end{aligned}
\end{equation*}
where we have used that, for any function $g=g(\be_{0},B_{\vzero}(\be_{0}))$, one has
$\partial_{\be_{0}}g=\partial_1 g + \partial_{\be_0}B_{\vzero} \, \partial_2 g$.

To obtain (\ref{eq:6.4c}) one can reason as follows.
For any self-energy cluster $T$, denote by $\ol{N}(T)$ the set of nodes
$v\in N(T)$ such that $\ell_{v}\in\calP_{T}\cup\{\ell_{T}\}$ and the line $\ell'_{v}\in\calP_{T}\cup
\{\ell'_{T}\}$ entering $v$ has component $h_{\ell'_{v}}=h_{\ell_{v}}$.
Let $T\in\gotS^{k}_{n,\be,\be}$ and consider the self-energy cluster
$T'\in\gotS^{k}_{n,B,B}$ obtained from $T$ by changing all the component
labels of the lines in $\calP_{T}\cup\{\ell_{T}\}\cup\{\ell_{T}'\}$ and reversing their orientation;
in particular the entering line $\ell_{T}'$ of $T$ becomes the exiting
line $\ell_{T'}$ of $T'$ and, vice versa, the exiting line $\ell_{T}$ of $T$ becomes the entering
line $\ell_{T'}'$ of $T'$. Both $\MM_{\be,\be}^{(k)}(0)$ and $\MM_{B,B}^{(k)}(0)$
are given by (\ref{eq:6.3b}) with $x=0$, so that
the external lines $\ell_{T}'$ and $\ell_{T}$ carry momenta $\nn_{\ell_{T}'}=\nn_{\ell_{T}}=\vzero$. 
Hence any line $\ell'\in\calP_{T'}$ has momentum $\nn'=-\nn$ if $\nn$ is the
momentum of the corresponding line in $\calP_{T}$ and the corresponding
propagator changes sign (see (\ref{eq:6.2}) and recall that
$n_{\ell}\ge0$ for all $\ell\in\calP_{T}$). 
Moreover, for any $v\in\ol{N}(T)$ the node factor $\calF_{v}$ changes
sign when regarded as a node in $\ol{N}(T')$, see (\ref{eq:6.1}).
All the other factors remains the same i.e. we can write
$\Val_{T}(0)=\gotA(T)\Val(\calP_{T})$ and
$\Val_{T'}(0)=\gotA(T')\Val(\calP_{T'})$, where
$$
\Val(\calP_{T}):= \Biggl(\prod_{v\in\ol{N}(T)}\calF_{v}\Biggr)
\Biggl(\prod_{\ell\in \calP_{T}}\calG_{n_{\ell}}(\oo\cdot\nn_{\ell})\Biggr)
$$
and  $\Val(\calP_{T'})$ is analogously defined with $T'$ instead of $T$,
while $\gotA(T)=\gotA(T')$. Now, one has
$$
\prod_{v\in \ol{N}(T)}\calF_{v}= \s \!\!\!\!\prod_{v\in \ol{N}(T')}\calF_{v}
$$
with $\s=\pm 1$. If $\s=1$, then $|\ol{N}(T)|=|\ol{N}(T')|$ is even and hence there is an odd number of
lines in $\calP_{T}$. If on the contrary $\s=-1$,
then there is an even number of lines in $\calP_{T}$. In both cases the assertion follows.\EP

\begin{lemma}\label{lem:6.6}
Assume that $\Ga_{\vzero}^{(k)}\equiv 0$ for all $k\ge 0$. Then
for all  $k\ge 1$ one has
\begin{subequations}
\begin{align}
& \partial_{x} \MM^{(k)}_{B,\be}(0) = \partial_{x} \MM^{(k)}_{\be,B}(0) =0,
\label{eq:6.5a} \\
& \partial_{x} \MM^{(k)}_{\be,\be}(0) = \partial_{x} \MM^{(k)}_{B,B}(0)  .
\label{eq:6.5b}
\end{align}
\label{eq:6.5}
\end{subequations}
\vskip-.3truecm
\end{lemma}

\vskip-0.5truecm
\prova
One reason along the same lines as the proof of (\ref{eq:6.4c}) in Lemma \ref{lem:6.5}.
Let $T\in\gotS^{k}_{n,\be,B}$ and consider the self-energy cluster
$T'\in\gotS^{k}_{n,\be,B}$ obtained from $T$ by changing all the component
labels of the lines in $\calP_{T}\cup\{\ell_{T}\}\cup\{\ell_{T}'\}$
and reversing their orientation. The derivative $\partial_{x}$ acts on
the propagator of some line $\ell\in\calP_{T}$.
After differentiation, when computing the propagators at $x=0$,
any line $\ell'\in\calP_{T'}$ turns out to have momentum $\nn'=-\nn$, if $\nn$ is the
momentum of the corresponding line in $\calP_{T}$ and the corresponding
propagator changes sign, except the differentiated propagator
$\left. \partial_{x}\calG_{n_{\ell}}(\oo\cdot\nn_{\ell}^0+x)
\right|_{x=0}$, which is even in its argument.
Moreover, for any $v\in\ol{N}(T)$ (we use the same notations as in the
proof of Lemma \ref{lem:6.5}) the node factor $\calF_{v}$ changes
sign when regarded as a node in $\ol{N}(T')$, while all the other factors remains the same.
If $|\ol{N}(T)|=|\ol{N}(T')|$ is even (resp. odd) then there is an even (resp. odd) number of
lines in $\calP_{T}$, but, as we said, the differentiated propagator
does not change sign. Therefore the two contributions have the same
modulus but different sign, so that, once summed together, they gives zero.
Therefore (\ref{eq:6.5a}) is proved.

To prove (\ref{eq:6.5b}) reason as in proving (\ref{eq:6.4c}): the
only difference is that, as in the previous case, the differentiated
propagator does not change sign, so that the two contributions are
equal to (and not the opposite of) each other.\EP

\begin{rmk}\label{rmk:6.7}
\emph{
Note that the Hamiltonian structure is fundamental in order to prove
both the identity (\ref{eq:6.4c}) and the identities (\ref{eq:6.5}).
}
\end{rmk}

Given $p\ge 2$ self-energy clusters $T_{1},\ldots,T_{p}$ of any tree $\theta$, with
$\ell'_{T_{i}}=\ell_{T_{i+1}}$ for $i=1,\ldots,p-1$ and $\ell_{T_{1}},\ell'_{T_{p}}$ being non-resonant,
we say that $C=\{T_{1},\ldots,T_{p}\}$ is a \emph{chain}.
Define $\ell_{0}(C):=\ell_{T_{1}}$ and $\ell_{i}(C):=\ell_{T_{i}}'$
for $i=1,\ldots,p$ and set $n_{i}(C)=n_{\ell_{i}(C)}$ for
$i=0,\ldots,p$; we also  call $k(C):=k(T_{1})+\ldots+k(T_{p})$ the \emph{total order}
of the chain $C$ and $p(C)=p$ the \emph{length} of $C$.
Given a chain $C=\{T_{1},\ldots,T_{p}\}$ we define the \emph{value} of $C$ as
\begin{equation}\label{eq:6.6}
\Val_{C}(x)=\prod_{i=1}^{p}\Val_{T_{i}}(x).
\end{equation}
We denote by $\gotC(k;h,h';n_{0},\ldots,n_{p})$ the set of all chains $C=\{T_{1},\ldots,T_{p}\}$
with total order $k$ and with fixed labels $h_{\ell_{0}(C)}=h$,
$h_{\ell_{p}(C)}=h'$ and $n_{i}(C)=n_{i}$ for $i=0,\ldots,p$.

\begin{rmk}\label{rmk:6.8}
\emph{
Let $\ell$ be a resonant line. Then there exists a chain $C$ such that
$\ell=\ell_{i}(C)$ for some $i=1,\ldots,p(C)-1$. If there exists a minimal
self-energy cluster $T$ containing $\ell$, then $T$ contains the whole chain $C$ and
all lines $\ell_{0}(C),\ldots,\ell_{p}(C)$
(this follows from the fact that, by definition of self-energy cluster,
$\nn_{\ell'}^{0} \neq \vzero$ for all $\ell'\in\calP_{T}$).
In particular $L(T)$ contains the two non-resonant lines $\ell_{0}(C)$ and
$\ell_{p}(C)$, with $\ze_{\ell_{0}(C)}=\ze_{\ell_{p}(C)}=\ze_{\ell}$.
}
\end{rmk}  

\begin{lemma}\label{lem:6.9}
Assume that $\Ga_{\vzero}^{(k')}\equiv 0$ for all $k'\ge 0$. Then
for all $p\ge2$, all $k\ge 1$, all $h,h'\in\{\be,B\}$ and all $\ol{n}_{0},\ldots,\ol{n}_{p}\in\ZZZ_{+}$
such that $\Psi_{\ol{n}_{i}}(x) \neq 0$, $i=0,\ldots,p$, one has
\begin{equation}\label{eq:6.7}
\Biggl| \sum_{C\in \gotC(k;h,h';\ol{n}_{0},\ldots,\ol{n}_{p})} \Val_{C}(x) \Biggr| \le B^{k}|x|^{p-1} ,
\end{equation}
for some constant $B>0$.
\end{lemma}

The proof is deferred to Appendix \ref{app:c}.

The bound (\ref{eq:6.7}) provides exactly the gain
factor which is needed in order to prove the convergence of
the power series. Indeed given a tree $\theta$, sum together
the values of the trees obtained from $\theta$ by replacing
each maximal chain $C$ (i.e. each chain which is not contained inside any other chain)
with any other chain which has the same total
order, the same length, the same scale labels associated with the
lines $\ell_{0}(C),\ldots,\ell_{p}(C)$ and the same component labels
associated with the lines $\ell_{0}(C)$ and $\ell_{p}(C)$; in other words,
if $C\in \gotC(k;h,h';\ol{n}_{0},\ldots,\ol{n}_{p})$ for some values of the labels, sum over
all possible chains belonging to the set $\gotC(k;h,h';\ol{n}_{0},\ldots,\ol{n}_{p})$.
Then we can bound the product of the propagators of the non-resonant lines outside
the maximal chains thanks to Lemma \ref{lem:6.3}, while the product of the
propagators of the lines $\ell_{1}(C),\ldots,\ell_{p-1}(C)$ of any chain $C$
times the sum of the corresponding chain values is bounded through Lemma \ref{lem:6.9}.
Then we can reason as in the proof of Lemma \ref{lem:4.5}, by using
that the propagator of any line $\ell$ is bounded proportionally to $\al_{m_{\ze_{\ell}}}(\oo)^{-1}$.

\begin{rmk}\label{rmk:6.11}
\emph{
We obtained the convergence of the power series (\ref{eq:2.8})
for any $\be_{0}$ and any $\e$ small enough. Thus the solution
turns out to be analytic in both $\e$ and $\be_{0}$. Moreover, since the
solution is parameterised by $\be_0\in\TTT$, in that case the full
resonant torus survives. Of course, such a situation
is highly non-generic and hence very unlikely.
}
\end{rmk}

\appendix 
 
\zerarcounters 
\section{Proof of Lemma \ref{lem:4.13}} 
\label{app:a} 

\begin{figure}[ht] 
\centering 
\ins{135pt}{-16pt}{$n$} 
\ins{237pt}{-24pt}{$n$} 
\ins{330pt}{-16pt}{$n+1$} 
\ins{135pt}{-35pt}{$\ell_{T}$}
\ins{290pt}{-50pt}{$T$}
\ins{330pt}{-35pt}{$\ell_{T}'$}
\ins{125pt}{-096pt}{$n+1$} 
\ins{237pt}{-104pt}{$n$} 
\ins{340pt}{-096pt}{$n$} 
\ins{135pt}{-115pt}{$\ell_{T}$}
\ins{290pt}{-130pt}{$T$}
\ins{330pt}{-115pt}{$\ell_{T}'$}
\includegraphics[width=4in]{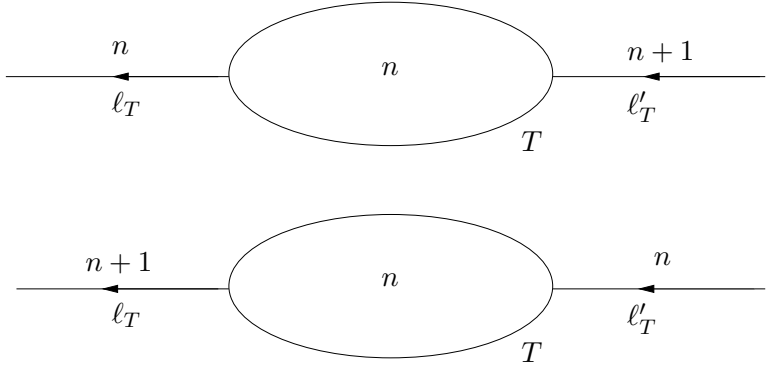} 
\vskip.2truecm 
\caption{Graphical representation of the left-fake (top) and right-fake (bottom) clusters $T$ 
on scale $n$; by construction $\nn_{\ell_{T}}=\nn_{\ell_{T}'}$. As for
the self-energy clusters in Figure \ref{fig:34} one has $\ell_{T},\ell_{T}'\notin L(T)$.}
\label{fig:35} 
\end{figure}

A \emph{left-fake cluster} $T$ on scale $n$ is a connected 
subgraph of a tree $\theta$ with only one entering line $\ell'_{T}$ 
and one exiting line $\ell_{T}$ such that (i) all the lines in $T$ have 
scale $\le n$ and there is in $T$ at least one line on scale $n$, 
(ii) $\ell'_{T}$ is on scale $n+1$ 
and $\ell_{T}$ is on scale $n$ and (iii) one has 
$\nn_{\ell_{T}}=\nn_{\ell'_{T}}$. 
Analogously a \emph{right-fake cluster} $T$ on scale $n$ is 
a connected subgraph of a tree $\theta$ with only one entering 
line $\ell'_{T}$ and one exiting line $\ell_{T}$ such that (i) all the 
lines in $T$ have scale $\le n$ and there is in $T$ at least one line 
on scale $n$, (ii)  $\ell'_{T}$ is on scale $n$ 
and $\ell_{T}$ is on scale $n+1$ and (iii) one has 
$\nn_{\ell_{T}}=\nn_{\ell'_{T}}$. 
Roughly speaking, a left-fake (respectively right-fake) cluster $T$ fails to 
be a self-energy cluster (or even a cluster) only because the exiting
(respectively the entering) line is on scale equal to the scale of $T$. 
Left-fake and right-fake clusters will be represented graphically as in Figure \ref{fig:35}.

The sets of renormalised left-fake clusters and renormalised right-fake clusters $T$
on scale $n$ such that $u_{\ell_{T}}=u$ and $e_{\ell'_{T}}=e$ will be denoted by
$\gotL\gotF_{n,u,e}$ and $\gotR\gotF_{n,u,e}$, respectively.

\begin{rmk}\label{rmk:a.-1} 
\emph{ 
If $T$ is a renormalised left-fake (respectively right-fake) cluster, 
we can (and shall) write $\Val(T;\e,\be_{0},B_{0})=
\Val_{T}(\oo\cdot\nn_{\ell'_{T}};\e,\be_{0},B_{0})$ to stress that the
propagators of the lines in $\calP_{T}$ depend on 
$\oo\cdot\nn_{\ell'_{T}}$. In particular one has 
$$ 
\sum_{T\in \gotL\gotF_{n,u,e}}\e^{k(T)}\Val_{T}(x;\e,\be_{0},B_{0})= 
\sum_{T\in \gotR\gotF_{n,u,e}}\e^{k(T)}\Val_{T}(x;\e,\be_{0},B_{0})= 
M^{[n]}_{u,e}(x;\e,\be_{0},B_{0}). 
$$ 
} 
\end{rmk} 

Throughout this appendix, for sake of simplicity, we shall omit the adjective ``renormalised''
referred to trees, self-energy clusters, left-fake clusters and right-fake clusters.

We shall prove explicitly only the bound
\begin{equation} \label{eq:a.0}
|\MM_{\be,\be}^{[p]}(0;\e,\be_{0},B_{0})-
\partial_{\be_{0}}\Phi_{\vzero}^{\RR,p}(\e,\be_{0},B_{0})|
\le |\e| \, A_{1} e^{-A_{2}2^{m_{p+1}}},
\end{equation}
as the others relations in (\ref{eq:4.8}) can be proved exactly in the same way.

We want to compute $\partial_{\be_{0}}\Phi_{\vzero}^{\RR,p}(\e,\be_{0},B_{0})$,
with $\Phi_{\vzero}^{\RR,p}(\e,\be_{0},B_{0})$ given by the first line
of (\ref{eq:3.15}). We start by considering trees $\theta\in
\Theta^{\RR,{p}}_{k,\vzero,\be}$ such that 
\begin{equation} \label{eq:a.0bis} 
\max_{\ell\in \Theta^{\RR,{p}}_{k,\vzero,\be}}
\{ n \in \ZZZ_{+} : \Psi_{n}(\oo\cdot\nn_{\ell}) \neq 0 \} \le p ,
\end{equation}
and shall see later how to deal with trees in $\Theta^{\RR,{p}}_{k,\vzero,\be}$
for which the condition (\ref{eq:a.0bis}) is not satisfied (see case {\bf 7}
at the end).

First of all, for any tree $\theta$ set 
\begin{equation} \label{eq:a.1} 
\partial_{v}\Val(\theta;\e,\be_{0},B_{0}) := 
\partial_{\be_{0}}\calF_{v}\Biggl(\prod_{w \in N(\theta)\setminus\{v\}} 
\calF_{w}\Biggr) 
\Biggl(\prod_{\ell\in L(\theta)}\calG^{[n_{\ell}]}_{e_{\ell},u_{\ell}}
(\oo\cdot\nn_{\ell}; \e,\be_{0},B_{0})\Biggr),
\end{equation} 
and 
\begin{equation}\label{eq:a.2} 
\begin{aligned} 
\partial_{\ell}\Val(\theta;\e,\be_{0},B_{0})&:= 
\partial_{\be_{0}}\calG^{[n_{\ell}]}_{e_{\ell},u_{\ell}}(x_{\ell};\e,\be_{0},B_{0}) 
\Biggl(\prod_{v\in N(\theta)}\calF_{v} \Biggr)
\Biggl(\prod_{\la\in L(\theta)\setminus\{\ell\}} 
\calG^{[n_{\la}]}_{e_{\la},u_{\la}}(x_{\la};\e,\be_{0},B_{0})\Biggr)\\ 
 &=\calA_{\ell} (\theta,x_{\ell};\e,\be_{0},B_{0}) \, 
\partial_{\be_{0}}\calG^{[n_{\ell}]}_{,e_{\ell},u_{\ell}}
(x_{\ell};\e,\be_{0},B_{0})\, 
\BB_{\ell}(\theta;\e,\be_{0},B_{0}), 
\end{aligned} 
\end{equation} 
where $x_{\ell}:=\oo\cdot\nn_{\ell}$, 
$\partial_{\be_{0}}\calG^{[n_{\ell}]}_{e_{\ell},u_{\ell}}(x_{\ell};\e,\be_{0},B_{0})$ 
is written according to Remark \ref{rmk:3.5} and
\begin{subequations}
\begin{align}
& \calA_{\ell}(\theta,x_{\ell};\e,\be_{0},B_{0}):= 
\Biggl(\prod_{\substack{v\in N(\theta)\\ v \not\prec \ell}}\calF_{v}\Biggr) 
\Biggl(\prod_{\substack{\ell'\in L(\theta)\\ \ell' \not\preceq\ell}}
\calG^{[n_{\ell'}]}_{e_{\ell'},u_{\ell'}}  (x_{\ell'};\e,\be_{0},B_{0})\Biggr), 
\label{eq:a.3} \\
& \BB_{\ell}(\theta;\e,\be_{0},B_{0}):= 
\Biggl(\prod_{\substack{v\in N(\theta)\\ v\prec \ell}}\calF_{v}\Biggr) 
\Biggl(\prod_{\substack{\ell'\in L(\theta)\\ \ell' \prec\ell}}
\calG^{[n_{\ell'}]}_{e_{\ell'},u_{\ell'}} (x_{\ell'};\e,\be_{0},B_{0})\Biggr) . 
\label{eq:a.4} 
\end{align}
\label{eq:a.3-4}
\end{subequations}
\vskip-.3truecm
%
Let us define in the analogous way $\partial_{v}\Val_{T}(x;\e,\be_{0},B_{0})$ 
and $\partial_{\ell}\Val_{T}(x;\e,\be_{0}B_{0})$ for any self-energy cluster
$T$ and let us write 
\begin{equation} \label{eq:a.5} 
\partial_{\be_{0}}\Val(\theta;\e,\be_{0},B_{0})= 
\partial_{N}\Val(\theta;\e,\be_{0},B_{0}) + 
\partial_{L}\Val(\theta;\e,\be_{0},B_{0}), 
\end{equation} 
where 
\begin{equation} \label{eq:a.6} 
\partial_{N}\Val(\theta;\e,\be_{0},B_{0}) :=\sum_{v \in N(\theta)} 
\partial_{v}\Val(\theta;\e,\be_{0},B_{0}), 
\end{equation} 
and 
\begin{equation}\label{eq:a.7} 
\begin{aligned} 
\partial_{L}\Val(\theta;\e,\be_{0},B_{0})&:=\sum_{\ell \in L(\theta)} 
\partial_{\ell}\Val(\theta;\e,\be_{0},B_{0}). \\ 
\end{aligned} 
\end{equation} 
Let us also write 
\begin{equation}\label{eq:a.8} 
\partial_{\be_{0}}\Val_{T}(x;\e,\be_{0},B_{0})= 
\partial_{N}\Val_{T}(x;\e,\be_{0},B_{0})+\partial_{L}\Val_{T}(x;\e,\be_{0},
B_{0}), 
\end{equation} 
for any $T\in\gotR_{n,u,e}$, $n\ge0$ and $u,e\in\{\be,B\}$, where the
derivatives $\partial_{N}$ 
and $\partial_{L}$ are defined analogously with the previous cases 
(\ref{eq:a.6}) and (\ref{eq:a.7}), with $N(T)$ and $L(T)$ replacing 
$N(\theta)$ and $L(\theta)$, respectively, so that we can split 
\begin{equation}\label{eq:a.9} 
\begin{aligned} 
\partial_{\be_{0}}\Phi^{\RR,p}_{\vzero} (x;\e,\be_{0},B_{0})&= 
\partial_{N}\Phi^{\RR,p}_{\vzero} (x;\e,\be_{0},B_{0})+\partial_{L}\Phi^{\RR,p}_{\vzero} (x;\e,\be_{0},B_{0}), \\ 
\partial_{\be_{0}}M^{[n]}(x;\e,\be_{0},B_{0})&= 
\partial_{N}M^{[n]}(x;\e,\be_{0},B_{0})+\partial_{L}M^{[n]} (x;\e,\be_{0},B_{0}), \\ 
\partial_{\be_{0}}\MM^{[n]}(x;\e,\be_{0},B_{0})&= 
\partial_{N}\MM^{[n]}(x;\e,\be_{0},B_{0})+\partial_{L}\MM^{[n]}(x;\e,\be_{0},
B_{0}) , 
\end{aligned} 
\end{equation} 
again with obvious meaning of the symbols. 

\begin{rmk}\label{rmk:a.1} 
\emph{ 
We can interpret the derivative $\partial_{v}$ as all the possible 
ways to attach an extra line $\ell$ (with $\nn_{\ell}=\vzero$
and $u_{\ell}=\be$) to 
the node $v$, so that 
$$ 
\sum_{k\ge0}\e^{k+1}\sum_{\theta\in\Theta^{\RR,p}_{k+1,\vzero}}\partial_{N} 
\Val(\theta;\e,\be_{0},B_{0}), 
$$ 
produces contributions to $\MM^{[p]}_{\be,\be}(0;\e,\be_{0},B_{0})$. 
} 
\end{rmk} 
 
In order to compute $\partial_{\be_{0}}\Phi_{\vzero}^{\RR,p}(\e,\be_{0},B_{0})$,
we have to study the derivative (\ref{eq:a.5}) for any $\theta \in \Theta^{\RR,p}_{k,\vzero,\be}$.
The terms (\ref{eq:a.6}) produce immediately 
contributions to $\MM^{[p]}_{\be,\be}(0;\e,\be_{0},B_{0})$ by Remark \ref{rmk:a.1}. 
Thus, we have to study the derivatives $\partial_{\ell}\Val(\theta;\e,\be_{0},
B_{0})$ appearing in the sum (\ref{eq:a.7}). From now on, we shall not write
any longer explicitly the dependence on $\e,\be_{0}$ and $B_{0}$, 
in order not to overwhelm the notation. 

For any $\theta \in \Theta^{\RR,{p}}_{k,\vzero,\be}$
satisfying the condition (\ref{eq:a.0bis}) and for any line $\ell\in L(\theta)$, 
either there is only one scale $n$ such that $\Psi_{n}(x_{\ell}) \neq 0$ 
(and in that case $\Psi_{n}(x_{\ell})=1$ and $\Psi_{n'}(x_{\ell})=0$ 
for all $n'\neq n$) or there exists only one 
$0\le n\le p-1$ such that $\Psi_{n}(x_{\ell})\Psi_{n+1}(x_{\ell})\ne 0$. 
To help following the argument below, we divide the discussion into several steps
(cases {\bf 1} to {\bf 7}), marking the end of each step with a white box ($\square$).

\noindent 
\textbf{1.} If $\Psi_{n}(x_{\ell})=1$ one has 
\begin{equation}\label{eq:a.10} 
\begin{aligned} 
\partial_{\ell}\Val(\theta)&=\calA_{\ell}(\theta,x_{\ell}) \, 
\left(\calG^{[n]}(x_{\ell})
\partial_{\be_{0}}\MM^{[n-1]}(x_{\ell}) 
\left((\ii x_{\ell})\uno-\MM^{[n-1]}
(x_{\ell})\right)^{-1}\right)_ {e_ {\ell},u_{\ell}} \BB_{\ell}(\theta)\\ 
&=\calA_{\ell}(\theta,x_{\ell}) \, 
\Big(\calG^{[n]}(x_{\ell}) 
\partial_{\be_{0}}\MM^{[n-1]}(x_{\ell}) 
\calG^{[n]}(x_{\ell})\Big)_{e_{\ell},u_{\ell}} \BB_{\ell}(\theta),
\end{aligned} 
\end{equation} 
with $\calA_{\ell}(\theta,x_{\ell})$ and $\BB_{\ell}(\theta)$ defined in (\ref{eq:a.3-4}).\qed

\begin{rmk}\label{rmk:a.2} 
\emph{ 
Note that if we split $\partial_{\be_{0}}=\partial_{N}+ 
\partial_{L}$ in (\ref{eq:a.10}), the term with $\partial_{N}\MM^{[n-1]} 
(x_{\ell})$ is a contribution to $\MM^{[p-1]}_{\be,\be}(0)$ and hence to
$\MM^{[p]}_{\be,\be}(0)$. 
} 
\end{rmk} 
  
If there is only one $0\le n\le p-1$ such that $\Psi_{n}(x_{\ell})
\Psi_{n+1}(x_{\ell}) 
\ne 0$, then $\Psi_{n}(x_{\ell})+\Psi_{n+1}(x_{\ell})=1$ and
$\chi_{q}(x_{\ell})=1$ for all 
$q=-1,\ldots,n-1$, so that $\psi_{n+1}(x_{\ell})=1$ and hence 
$\Psi_{n+1}(x_{\ell})=\chi_{n}(x_{\ell})$. Moreover
it can happen only (see Remark \ref{rmk:3.8}) $n_{\ell}=n$ or $n_{\ell}=n+1$. 

\begin{figure}[ht] 
\centering 
\ins{050pt}{-40pt}{$n\!+\!1$} 
\ins{106pt}{-47pt}{$\le n$} 
\ins{154pt}{-40pt}{$n\!+\!1$} 
\ins{234pt}{-48pt}{$+$} 
\ins{298pt}{-40pt}{$n\!+\!1$} 
\ins{342pt}{-47pt}{$\le n\!-\!1$} 
\ins{416pt}{-41pt}{$n$} 
\ins{090pt}{-113pt}{$+$} 
\ins{152pt}{-105pt}{$n\!+\!1$} 
\ins{200pt}{-112pt}{$\le n\!-\!1$} 
\ins{260pt}{-106pt}{$n$} 
\ins{318pt}{-113pt}{$n$} 
\ins{363pt}{-105pt}{$n\!+\!1$} 
\includegraphics[width=6in]{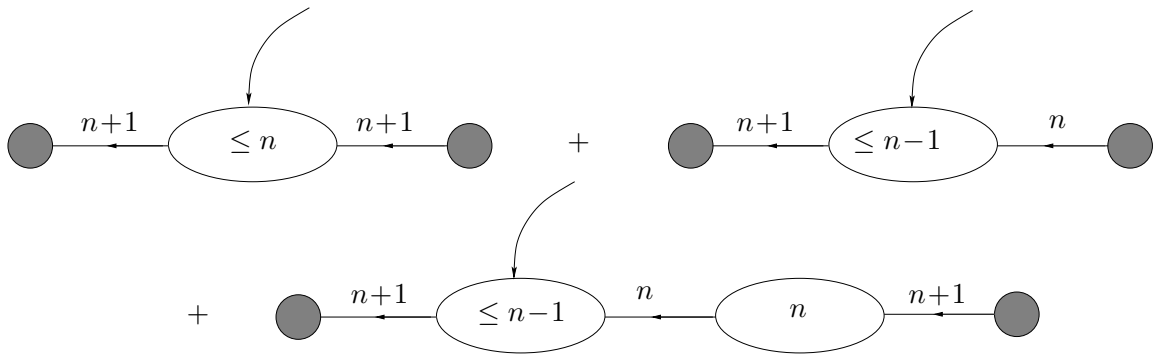} 
\caption{Graphical representation of the derivative 
$\partial_{\ell}\Val(\theta)$ according to (\ref{eq:a.11}).} 
\label{fig:a1} 
\end{figure}

\noindent 
\textbf{2.} Consider first the case $n_{\ell}=n+1$. One has
\begin{equation} \label{eq:a.11a}
\partial_{\ell}\Val(\theta) = \calA_{\ell}(\theta,x_{\ell}) \, 
\left(\calG^{[n+1]}(x_{\ell}) \partial_{\be_{0}}\MM^{[n]}(x_{\ell}) 
((\ii x_{\ell})\uno-\MM^{[n]}(x_{\ell}))^{-1}\right)_{e_{\ell},u_{\ell}} \BB_{\ell}(\theta) , 
\end{equation}
with
\begin{equation} \label{eq:a.11}
\begin{aligned} 
\calG^{[n+1]} & (x_{\ell}) 
\partial_{\be_{0}}\MM^{[n]}(x_{\ell}) 
\big((\ii x_{\ell})\uno- \MM^{[n]}(x_{\ell}) \big)^{-1} \\ 
& = \calG^{[n+1]}(x_{\ell}) 
  \partial_{\be_{0}}\MM^{[n-1]}(x_{\ell}) 
\big(\Psi_{n}(x_{\ell})+\Psi_{n+1}(x_{\ell}) \big)
\big((\ii x_{\ell})\uno-\MM^{[n]}(x_{\ell}) \big)^{-1} \\ 
& \quad + \;  \calG^{[n+1]}(x_{\ell}) 
  \partial_{\be_{0}}M^{[n]}(x_{\ell}) \, \chi_{n}(x_{\ell})
\big( (\ii x_{\ell})\uno-\MM^{[n]}(x_{\ell}) \big)^{-1} \\ 
&  = \calG^{[n+1]}(x_{\ell}) 
  \Biggl(\sum_{q=-1}^{n}\partial_{\be_{0}}M^{[q]}(x_{\ell})\Biggr) 
  \calG^{[n+1]}(x_{\ell}) + \calG^{[n+1]}(x_{\ell}) 
  \Biggl(\sum_{q=-1}^{n-1}\partial_{\be_{0}}M^{[q]}(x_{\ell})\Biggr) 
  \calG^{[n]}(x_{\ell}) \\ 
& \quad+ \; \calG^{[n+1]}(x_{\ell}) 
  \Biggl(\sum_{q=-1}^{n-1}\partial_{\be_{0}}M^{[q]}(x_{\ell}) \Biggr)
\calG^{[n]}(x_{\ell}) M^{[n]}(x_{\ell}) 
  \calG^{[n+1]}(x_{\ell}) ,
\end{aligned} 
\end{equation} 
where we have used that $\chi_{n}(x_{\ell})=\Psi_{n+1}(x_{\ell})$ and
$$ \big( (\ii x_{\ell})\uno- \MM^{[n]}(x_{\ell}) \big)^{-1} 
\big( \uno + M^{[n]}(x_{\ell}) \Psi_{n+1}(x_{\ell})
\big( (ix_{\ell})\uno-\MM^{[n]}(x_{\ell}) \big)^{-1} \big)^{-1}=
\big( (\ii x_{\ell})\uno- \MM^{[n-1]}(x_{\ell}) \big)^{-1} . $$
We represent graphically the three contributions in (\ref{eq:a.11}) as in Figure \ref{fig:a1}: 
we represent the derivative $\partial_{\be_{0}}$ as an arrow pointing
toward the graphical representation of the differentiated quantity.\qed
 
\begin{rmk}\label{rmk:a.3} 
\emph{ 
Note that the $M^{[n]}(x_{\ell})$ appearing in the latter line of 
(\ref{eq:a.11}) has 
to be interpreted (see Remark \ref{rmk:3.7}) as the matrix with
components
$$ 
\sum_{T\in\gotL\gotF_{n,u,e}}\e^{k(T)}\Val_{T}(x_{\ell}). 
$$ 
Note also that, again, if we split $\partial_{\be_{0}}=\partial_{N}+ 
\partial_{L}$ in (\ref{eq:a.11}), all the terms with $\partial_{N}M^{[q]} 
(x_{\ell})$ are contributions to $\MM^{[p]}_{\be,\be}(0)$. 
} 
\end{rmk} 

Now consider the case $n_{\ell}=n$. We distinguish among several cases
(see Remark \ref{rmk:a.4bis} 
below for the meaning of ``removal'' and ``insertion'' of left-fake clusters):\\
(a) $\ell$ \textit{does not} exit any left-fake cluster and one \textit{can}
insert a left-fake cluster, together its entering line, between $\ell$ and the
node $\ell$ exists without creating any self-energy cluster (case {\bf 3}
below);\\
(b) $\ell$ \textit{does not} exit any left-fake cluster and one
\textit{cannot} insert any left-fake cluster between $\ell$ and the node
$\ell$ exists because
this way a self-energy cluster would appear (case {\bf 4} below);\\
(c) $\ell$ \textit{does} exit a left-fake cluster and one \textit{can} remove
the left-fake cluster, together its entering line, without creating
a self-energy cluster (case {\bf 3} below);\\
(d) $\ell$ \textit{does} exit a left-fake cluster and one \textit{cannot}
remove the left-fake cluster because a self-energy cluster would be produced
(case {\bf 5} below).

\begin{rmk}\label{rmk:a.4bis} 
\emph{
Here and henceforth, if $S$ is a subgraph with only one entering line
$\ell'_{S}=\ell_{v}$ and one exiting line $\ell_{S}$, by saying that we
``remove'' $S$
together with $\ell'_{S}$, we mean that we change $u_{\ell_{S}}$ into
$h_{v}$ and we also reattach the line $\ell_{S}$ to the node $v$ (so
that $\ell_{S}$ becomes the line exiting $v$).
Analogously, whenever we ``insert'' a subgraph $S$ with only one entering line
$\ell'$ between a line $\ell$ and the node $v$ which $\ell$ exits,
we mean that we set $u_{\ell'}=h_{v}$ and change  $u_{\ell}$ into $h_{w}$ if
$w\in N(S)$ is the node to which we reattach $\ell$ (and $\ell$ becomes the
line $\ell_{w}$ exiting $S$). 
}
\end{rmk}

\begin{figure}[ht] 
\centering 
\ins{028pt}{-11pt}{$\bar{\theta}=$} 
\ins{086pt}{-04pt}{$n$} 
\ins{086pt}{-20pt}{$\ell$} 
\ins{203pt}{-11pt}{$\theta'=$} 
\ins{270pt}{-04pt}{$n$} 
\ins{270pt}{-20pt}{$\ell$} 
\ins{322pt}{-13pt}{$n$} 
\ins{378pt}{-03pt}{$n\!+\!1$} 
\includegraphics[width=5.2in]{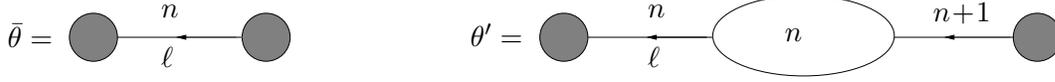} 
\caption{The renormalised tree $\bar{\theta}$ and the renormalised trees 
$\theta'$ of the set $\tau_{1}(\bar{\theta},\ell)$ associated with $\bar{\theta}$.} 
\label{fig:a2} 
\end{figure}

\noindent 
\textbf{3.} If $\ell$ is not the exiting line of a 
left-fake cluster, set $\bar{\theta}=\theta$; otherwise, 
if $\ell$ is the exiting line of a left-fake cluster $T$, 
define -- if possible -- $\bar{\theta}$ as the tree 
obtained from $\theta$ by removing $T$ and $\ell_{T}'$ 
In both cases, define -- if possible -- $\tau_{1}(\bar{\theta}, 
\ell)$ as the set constituted by all the renormalised trees $\theta'$ 
obtained from $\bar{\theta}$ by inserting a left-fake cluster, together 
with its entering line, between $\ell$ and the node $v$ which $\ell$ exits; see Figure \ref{fig:a2}.

\begin{rmk}\label{rmk:a.4} 
\emph{ 
The construction of the set $\tau_{1}(\bar{\theta},\ell)$ could be impossible 
if the removal or the insertion of a left-fake cluster $T$, together with its 
entering line $\ell'_{T}$, would produce a self-energy cluster. We shall see 
later (see cases {\bf 4} and {\bf 5} below) how to deal with these cases. 
} 
\end{rmk} 

\begin{figure}[ht] 
\centering 
\ins{057pt}{-41pt}{$n$} 
\ins{098pt}{-47pt}{$\le n-1$} 
\ins{162pt}{-41pt}{$n$} 
\ins{234pt}{-48pt}{$+$} 
\ins{310pt}{-41pt}{$n$} 
\ins{342pt}{-47pt}{$\le n\!-\!1$} 
\ins{400pt}{-40pt}{$n\!+\!1$} 
\ins{090pt}{-113pt}{$+$} 
\ins{160pt}{-106pt}{$n$} 
\ins{200pt}{-112pt}{$\le n\!-\!1$} 
\ins{265pt}{-106pt}{$n$} 
\ins{318pt}{-113pt}{$n$} 
\ins{363pt}{-105pt}{$n\!+\!1$} 
\includegraphics[width=6in]{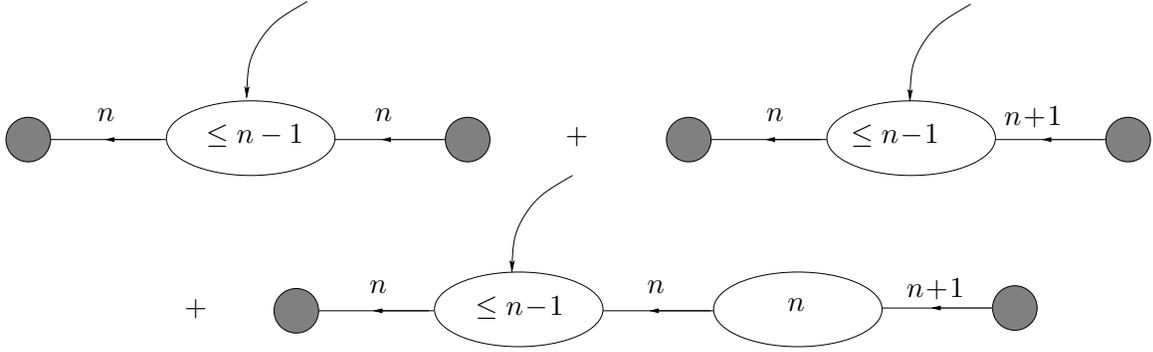} 
\caption{Graphical representation of the three contributions in 
the last two lines of (\ref{eq:a.13}).} 
\label{fig:a3} 
\end{figure}

Then one has 
\begin{equation}\label{eq:a.12} 
\partial_{\ell}\Val(\bar{\theta})+ 
\partial_{\ell}\!\!\!\sum_{\theta'\in\tau_{1}(\bar{\theta},\ell)}\!\!\!\Val(\theta') 
=\calA_{\ell}(\bar{\theta},x_{\ell}) \, \Big(\partial_{\be_{0}}\calG^{[n]}
(x_{\ell}) 
\left(\uno+M^{[n]}(x_{\ell})\calG^{[n+1]}(x_{\ell})\right)\Big)_{e_{\ell},u_{\ell}}
\, \BB_{\ell}(\bar{\theta}), 
\end{equation} 
with
\begin{equation} \nonumber 
\begin{aligned} 
\partial_{\be_{0}}&\calG^{[n]}(x_{\ell})\left(\uno+ 
  M^{[n]}(x_{\ell})\calG^{[n+1]}(x_{\ell})\right)\\ 
&=\calG^{[n]}(x_{\ell}) 
  \partial_{\be_{0}}\MM^{[n-1]}(x_{\ell}) 
  \calG^{[n]}(x_{\ell})\\ 
&\quad+\calG^{[n]}(x_{\ell}) 
  \partial_{\be_{0}}\MM^{[n-1]}(x_{\ell}) 
  \Psi_{n+1}(x_{\ell})
  \left((\ii x_{\ell})\uno-\MM^{[n-1]}(x_{\ell})\right)^{-1}\\ 
&\quad+\calG^{[n]}(x_{\ell}) 
  \partial_{\be_{0}}\MM^{[n-1]}(x_{\ell}) 
  \calG^{[n]}(x_{\ell}) 
  M^{[n]}(x_{\ell}) 
  \calG^{[n+1]}(x_{\ell})\\  
&\quad+\calG^{[n]}(x_{\ell}) 
  \partial_{\be_{0}}\MM^{[n-1]}(x_{\ell}) 
  \Psi_{n+1}(x_{\ell})
  \left((\ii x_{\ell})\uno-\MM^{[n-1]}(x_{\ell})\right)^{-1} 
  M^{[n]}(x_{\ell}) 
  \calG^{[n+1]}(x_{\ell})
\end{aligned} 
\end{equation} 
and hence
\begin{equation} \nonumber \label{eq:a.13} 
\begin{aligned} 
\partial_{\be_{0}}&\calG^{[n]}(x_{\ell})\left(\uno+ M^{[n]}(x_{\ell})\calG^{[n+1]}(x_{\ell})\right)\\ 
&=\calG^{[n]}(x_{\ell}) 
  \partial_{\be_{0}}\MM^{[n-1]}(x_{\ell}) 
  \calG^{[n]}(x_{\ell}) +\calG^{[n]}(x_{\ell}) 
  \partial_{\be_{0}}\MM^{[n-1]}(x_{\ell}) 
  \calG^{[n+1]}(x_{\ell})\\ 
&\quad-\calG^{[n]}(x_{\ell}) 
  \partial_{\be_{0}}\MM^{[n-1]}(x_{\ell}) 
  \chi_{n}(x_{\ell})
  \left((\ii x_{\ell})\uno-\MM^{[n-1]}(x_{\ell})\right)^{-1} 
  M^{[n]}(x_{\ell}) 
  \calG^{[n+1]}(x_{\ell})\\  
&\quad+\calG^{[n]}(x_{\ell}) 
  \partial_{\be_{0}}\MM^{[n-1]}(x_{\ell}) 
  \calG^{[n]}(x_{\ell}) 
  M^{[n]}(x_{\ell}) 
  \calG^{[n+1]}(x_{\ell})\\ 
&\quad+\calG^{[n]}(x_{\ell}) 
  \partial_{\be_{0}}\MM^{[n-1]}(x_{\ell}) 
  \Psi_{n+1}(x_{\ell})
  \left((\ii x_{\ell})\uno-\MM^{[n-1]}(x_{\ell})\right)^{-1}
  M^{[n]}(x_{\ell}) 
  \calG^{[n+1]}(x_{\ell})\\ 
&=\calG^{[n]}(x_{\ell}) 
  \partial_{\be_{0}}\MM^{[n-1]}(x_{\ell}) 
  \calG^{[n]}(x_{\ell}) +\calG^{[n]}(x_{\ell}) 
  \partial_{\be_{0}}\MM^{[n-1]}(x_{\ell}) 
  \calG^{[n+1]}(x_{\ell})\\ 
&\quad+\calG^{[n]}(x_{\ell}) 
  \partial_{\be_{0}}\MM^{[n-1]}(x_{\ell}) 
  \calG^{[n]}(x_{\ell}) 
  M^{[n]}(x_{\ell}) 
  \calG^{[n+1]}(x_{\ell}), 
\end{aligned} 
\end{equation} 
where we have used that
$$ \big( \uno - \Psi_{n+1}(x_{\ell})
\big( (ix_{\ell})\uno-\MM^{[n-1]}(x_{\ell}) \big)^{-1} M^{[n]}(x_{\ell}) \big)^{-1}
\big( (\ii x_{\ell})\uno- \MM^{[n-1]}(x_{\ell}) \big)^{-1} =
\big( (\ii x_{\ell})\uno- \MM^{[n]}(x_{\ell}) \big)^{-1} . $$
Also in this case, if we split $\partial_{\be_{0}}=\partial_{N}+ 
\partial_{L}$, all the terms with $\partial_{N}\MM^{[n-1]}$ are 
contributions to $\MM^{[p]}_{\be,\be}(0)$ -- see Remark \ref{rmk:a.2}. 
Again, we can represent graphically the three contributions obtained 
inserting (\ref{eq:a.13}) in (\ref{eq:a.12}); see Figure \ref{fig:a3}.\qed
 
\begin{figure}[ht] 
\centering 
\ins{028pt}{-11pt}{$\theta'=$} 
\ins{082pt}{-03pt}{$n\!+\!1$} 
\ins{090pt}{-19pt}{$\ol{\ell}$} 
\ins{144pt}{-12pt}{$n$} 
\ins{202pt}{-04pt}{$n$} 
\ins{202pt}{-20pt}{$\ell$} 
\ins{320pt}{-11pt}{$\ol{\theta}=$} 
\ins{384pt}{-20pt}{$\ol{\ell}$} 
\ins{378pt}{-03pt}{$n\!+\!1$} 
\includegraphics[width=5.2in]{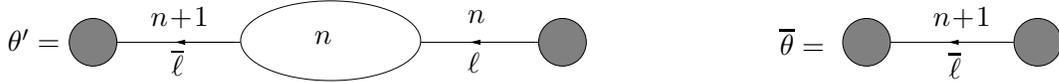} 
\caption{The  trees $\theta'$ of the set $\tau_{2}(\ol{\theta},\ol\ell)$ 
obtained from $\ol{\theta}$ when $\ell\in L(\theta)$ enters a right-fake 
cluster.} 
\label{fig:a4} 
\end{figure}

\noindent 
\textbf{4.} Assume now that $\ell$ is not the exiting line of a left-fake cluster 
and the insertion of a left-fake cluster, together 
with its entering line, produces a self-energy cluster. Note that this can 
happen only if $\ell$ is the entering line 
of a right-fake cluster $T$. Let $\ol{\ell}$ be 
the exiting line (on scale $n+1$) of the 
right-fake cluster $T$, call $\ol{\theta}$ the tree obtained 
from $\theta$ by removing $T$ and $\ell$ and call 
$\tau_{2}(\ol{\theta},\ol\ell)$ the set of trees $\theta'$ 
obtained from $\ol{\theta}$ by inserting a right-fake cluster, 
together with its entering line, before $\ol{\ell}$; see Figure \ref{fig:a4}. 
 
By construction one has 
\begin{equation} \nonumber 
\begin{aligned} 
\Val(\ol{\theta}) & = 
 \calA_{\ol{\ell}}(\ol{\theta},x_{\ell}) \, 
  \calG^{[n+1]}_{e_{\ol{\ell}},u_{\ol{\ell}}}(x_{\ol{\ell}}) 
  \, \BB_{\ol{\ell}}(\ol{\theta}) \\ 
\sum_{\theta'\in\tau_{2}(\ol{\theta},\ol\ell)} 
\!\!\!\Val(\theta') & = 
 \calA_{\ol{\ell}}(\ol{\theta},x_{\ell}) \, 
  \Big(\calG^{[n+1]}(x_{\ol{\ell}}) \, M^{[n]}(x_{\ol{\ell}})\, 
  \calG^{[n]}(x_{\ol{\ell}}) \Big)_{e_{\ol{\ell}},u_{\ol{\ell}}}
  \, \BB_{\ol{\ell}}(\ol{\theta}), 
\end{aligned} 
\end{equation} 
where $u_{\ol{\ell}}$ denotes the $u$-component of $\ol{\ell}$ as line in
$\ol{\theta}$ and we have used that $x_{\ell}=x_{\bar{\ell}}$.

\begin{figure}[ht] 
\centering 
\ins{193pt}{-20pt}{$T$} 
\ins{271pt}{-20pt}{$T''$} 
\ins{330pt}{-80pt}{$T'$} 
\ins{104pt}{-38pt}{$n\!+\!1$} 
\ins{178pt}{-45pt}{$n$} 
\ins{234pt}{-39pt}{$n$} 
\ins{292pt}{-45pt}{$n$} 
\ins{358pt}{-38pt}{$n\!+\!1$} 
\ins{236pt}{-54pt}{$\ell$} 
\ins{110pt}{-54pt}{$\ell_{T'}$} 
\ins{360pt}{-54pt}{$\ell_{T'}'$} 
\includegraphics[width=4.0in]{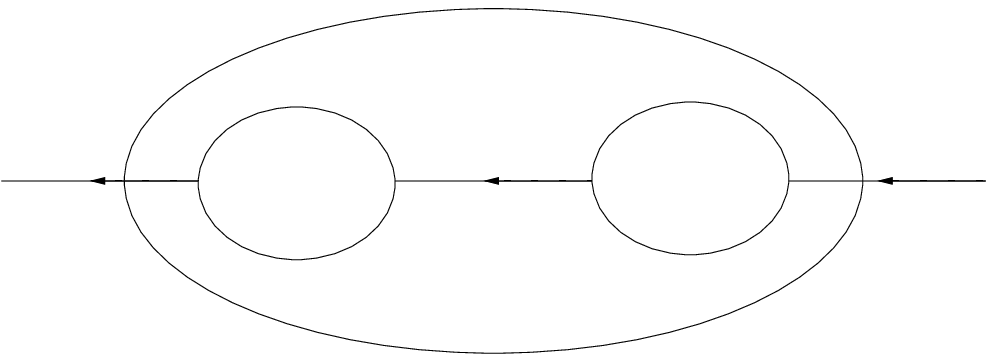} 
\caption{A self-energy cluster $T'\in \gotR_{n}(T)$.} 
\label{fig:a5} 
\end{figure}

Consider the contribution to 
$\partial_{\ol{\ell}}\Val(\ol{\theta})$ given by 
\begin{equation} \label{eq:a.14} 
\calA_{\ol{\ell}}(\ol{\theta},x_{\ol{\ell}}) 
  \Big(\calG^{[n+1]}(x_{\ol{\ell}}) 
  \partial_{L}M^{[n]}(x_{\ol{\ell}}) 
  \calG^{[n+1]}(x_{\ol{\ell}}) \Big)_{e_{\ol{\ell}},u_{\ol{\ell}}}
  \BB_{\ol{\ell}}(\ol{\theta}) ,
\end{equation} 
arising from (\ref{eq:a.11}). For $u,e,e'\in\{\be,B\}$ and $T\in\gotR\gotF_{n,u,e'}$ call
$\gotR_{n,u,e}(T)$ the subset of $\gotR_{n,u,e}$ such that 
if $T'\in\gotR_{n,u,e}(T)$ the exiting line $\ell_{T'}$ exits also 
the renormalised right-fake cluster $T$; note that the entering line 
$\ell$ of $T$ must be also the exiting line of some renormalised 
left-fake cluster $T''$ contained in $T'$; see Figure \ref{fig:a5}. 

\begin{figure}[ht] 
\centering 
\ins{040pt}{-28pt}{$n\!+\!1$} 
\ins{073pt}{-38pt}{$n$} 
\ins{105pt}{-29pt}{$n$} 
\ins{130pt}{-36pt}{$\le \! n\!-\!1$} 
\ins{185pt}{-29pt}{$n$} 
\ins{232pt}{-38pt}{$+$} 
\ins{272pt}{-28pt}{$n\!+\!1$} 
\ins{309pt}{-38pt}{$n$} 
\ins{338pt}{-29pt}{$n$} 
\ins{366pt}{-36pt}{$\le n\!-\!1$} 
\ins{418pt}{-28pt}{$n\!+\!1$} 
\ins{076pt}{-088pt}{$+$} 
\ins{119pt}{-080pt}{$n\!+\!1$} 
\ins{153pt}{-090pt}{$n$} 
\ins{185pt}{-081pt}{$n$} 
\ins{211pt}{-087pt}{$\le n\!-\!1$} 
\ins{266pt}{-081pt}{$n$} 
\ins{307pt}{-089pt}{$n$} 
\ins{340pt}{-080pt}{$n\!+\!1$} 
\includegraphics[width=6.0in]{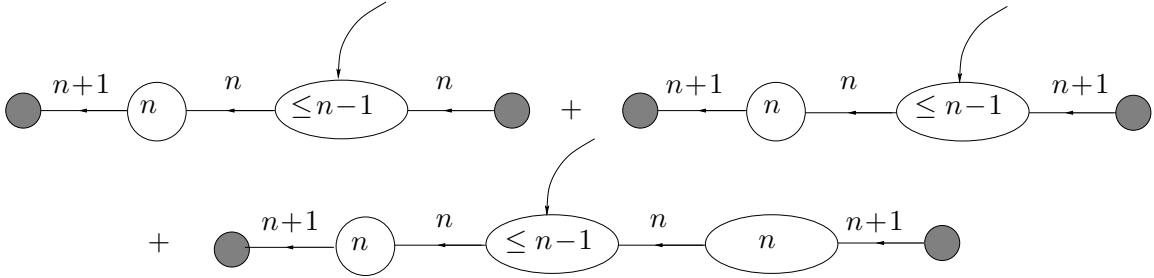} 
\caption{Graphical representation of the three contributions arising from (\ref{eq:a.16}).} 
\label{fig:a6} 
\end{figure}

Define $\gotM^{[n]}(x_{\ol{\ell}})$ as the $2\times 2$ matrix with components
\begin{equation} \label{eq:a.15} 
\gotM^{[n]}_{u,e}(x_{\ol{\ell}})= \sum_{e'=\be,B}\sum_{T\in\gotR\gotF_{n,u,e'}}
\sum_{T'\in\gotR_{n,u,e}(T)}\e^{k(T')}\Val_{T'}(x_{\ol{\ell}}) 
\end{equation} 
and consider the contribution $\gotM^{[n]}(x_{\ol{\ell}})$
to $M^{[n]}(x_{\ol{\ell}})$ in (\ref{eq:a.14}). Let us pick up the
term with the derivative acting on the line $\ell$: one has 
\begin{equation}\label{eq:a.16} 
\begin{aligned} 
& \partial_{\ell} \!\!\!\sum_{\theta'\in\tau_{2}(\ol{\theta},\ol{\ell})} 
\!\!\!\Val(\theta') + 
 \calA_{\ol{\ell}}(\ol{\theta},x_{\ell}) \, 
  \Big(\calG^{[n+1]}(x_{\ell}) \, \partial_{\ell} 
\gotM^{[n]}(x_{\ol{\ell}}) 
  \calG^{[n+1]}(x_{\ell})\Big)_{e_{\ol{\ell}},u_{\ol{\ell}}}
  \, \BB_{\ol{\ell}}(\ol{\theta})\\ 
&=\calA_{\ol{\ell}}(\ol{\theta},x_{\ell}) \, 
\Big(\calG^{[n+1]}(x_{\ell})  M^{[n]}(x_{\ell}) 
\partial_{\be_{0}}\calG^{[n]}(x_{\ell})\left(\uno+ 
M^{[n]}(x_{\ell})\calG^{[n+1]}(x_{\ell})\right)\Big)_{e_{\ol{\ell}},u_{\ol{\ell}}}
\BB_{\ol{\ell}}(\ol{\theta}) , 
\end{aligned} 
\end{equation} 
where we have used again that $x_{\ell}=x_{\ol{\ell}}$. Thus, one can reason 
as in (\ref{eq:a.13}), so as to obtain the sum 
of three contributions, as represented in Figure \ref{fig:a6}. \qed
 
\noindent 
\textbf{5.} Finally, consider the case in which $\ell$ is the 
exiting line of a left-fake cluster, 
$T_{0}$ and the removal of $T_{0}$ and $\ell'_{T_{0}}$ 
(see Remark \ref{rmk:a.4bis}) creates a self-energy cluster. 
 
Set (for a reason that will become clear later) 
$\theta_{0}=\theta$ and $\ell_{0}=\ell$. 
Then there is a maximal $m\ge1$ such that 
there are $2m$ lines $\ell_{1},\ldots,\ell_{m}$ and 
$\ell'_{1},\ldots\ell'_{m}$, with the following properties: 
 
\noindent 
(i) $\ell_{i}\in \calP(\ell_{\theta_{0}},\ell_{i-1})$, for $i=1,\ldots,m$, \\ 
(ii) $n_{\ell_{i}}=n+i<\max\{p:\Psi_{p}(x_{\ell_{i}})\ne0\}=n+i+1$, for 
$i=0,\ldots,m-1$, while $n_{m}:=n_{\ell_{m}}= 
n+m+\s$, with $\s\in\{0,1\}$,\\ 
(iii) $\nn_{\ell_{i}}\ne\nn_{\ell_{i-1}}$ and the lines preceding 
$\ell_{i}$ but not $\ell_{i-1}$ are on scale $\le n+i-1$, for 
$i=1,\ldots,m$, \\ 
(iv) $\nn_{\ell'_{i}}=\nn_{\ell_{i}}$, for $i=1,\ldots,m$,\\ 
(v) if $m\ge2$, $\ell'_{i}$ is the exiting line of a left-fake cluster 
$T_{i}$, for $i=1,\ldots,m-1$, \\ 
(vi) $\ell'_{i}\prec\ell'_{T_{i-1}}$ and all the lines 
preceding $\ell'_{T_{i-1}}$ but not $\ell'_{i}$ are on scale 
$\le n+i-1$,  for $i=1,\ldots,m$, \\ 
(vii) $n'_{m}:=n_{\ell'_{m}}=n+m+\s'$ with $\s'\in\{0,1\}$.

Note that one cannot have $\s=\s'=1$, otherwise the subgraph between
$\ell_{m}$ and $\ell'_{m}$ would be a self-energy cluster. 
Note also that (ii), (iv) and (v) imply $n_{\ell'_{i}}=n+i$ for $i=1,\ldots,m-1$ 
if $m\ge2$.
Call $S_{i}$ the subgraph between $\ell_{i+1}$ and $\ell_{i}$ and $S'_{i}$ 
the cluster between $\ell'_{T_{i}}$ and $\ell'_{i+1}$, for all $i=0,\ldots,m-1$. 
For $i=1,\ldots,m$, call $\theta_{i}$ the tree obtained from 
$\theta_{0}$ by removing everything between $\ell_{i}$ and the part of  
$\theta_{0}$ preceding $\ell'_{i}$ Note that, if $m\ge2$,
properties (i)--(vii) hold for $\theta_{i}$ but with $m-i$ instead of $m$, for all $i=1,\ldots,m-1$.

\begin{figure}[ht] 
\vskip.3truecm 
\centering 
\ins{028pt}{-23pt}{$\theta_{0}=$} 
\ins{082pt}{-16pt}{$n\!+\!1$} 
\ins{090pt}{-32pt}{$\ell_{1}$} 
\ins{135pt}{-24pt}{$\le n$} 
\ins{163pt}{-46pt}{$S_{0}$} 
\ins{187pt}{-17pt}{$n$} 
\ins{187pt}{-32pt}{$\ell_{0}$} 
\ins{236pt}{-25pt}{$n$} 
\ins{261pt}{-46pt}{$T_{0}$} 
\ins{278pt}{-16pt}{$n\!+\!1$} 
\ins{284pt}{-32pt}{$\ell_{T_{0}}'$} 
\ins{324pt}{-24pt}{$\le n$} 
\ins{354pt}{-46pt}{$S'_{0}$} 
\ins{380pt}{-32pt}{$\ell'_{1}$} 
\ins{374pt}{-16pt}{$n\!+\!1$} 
\ins{028pt}{-82pt}{$\theta_{1}=$} 
\ins{082pt}{-75pt}{$n\!+\!1$} 
\ins{090pt}{-92pt}{$\ell_{1}$} 
\ins{236pt}{-88pt}{$R_{1}$} 
\ins{236pt}{-129pt}{$n$} 
\ins{236pt}{-143pt}{$\ell_{0}$} 
\ins{187pt}{-133pt}{$\le n$} 
\ins{284pt}{-133pt}{$\le n$} 
\ins{263pt}{-155pt}{$S'_{0}$} 
\ins{213pt}{-155pt}{$S_{0}$} 
\includegraphics[width=5.2in]{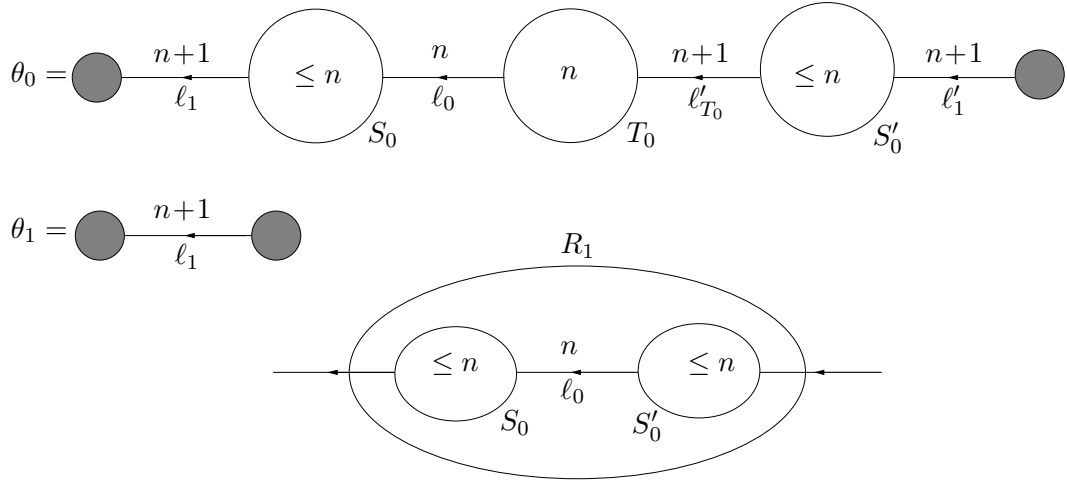} 
\caption{The renormalised trees $\theta_{0}$ and $\theta_{1}$ and the 
self-energy cluster $R_{1}$ in case {\bf 5} 
with $m=1$ and $\s=\s'=0$. Note that the set $S_{0}'$ is a 
cluster, but not a self-energy cluster.} 
\label{fig:a7} 
\end{figure}

For $i=1,\ldots,m$, call $R_{i}$ the self-energy 
cluster obtained from the subgraph of $\theta_{i-1}$ between $\ell_{i}$ and
$\ell'_{i}$, by removing the left-fake cluster $T_{i-1}$ together with
$\ell'_{T_{i}}$. Note that $L(R_{i})=L(S_{i-1})\cup\{\ell_{i-1}\}\cup 
L(S'_{i-1})$ and $N(R_{i})=N(S_{i-1})\cup N(S'_{i-1})$; see Figure  \ref{fig:a7}.
 
For $i=0,\ldots,m-1$, given $\ell',\ell\in L(\theta_{i})$, with 
$\ell'\prec\ell$, call $\calP^{(i)}(\ell,\ell')$ the path of lines in 
$\theta_{i}$ connecting 
$\ell'$ to $\ell$ (hence $\calP^{(i)}(\ell,\ell')=\calP(\ell,\ell') 
\cap L(\theta_{i})$). For any $i=0,\ldots,m-1$ and any $\ell\in\calP^{(i)} 
(\ell_{i},\ell'_{m})$, let 
$\tau_{3}(\theta_{i},\ell)$ be the set of all renormalised trees 
which can be obtained from $\theta_{i}$ by replacing each left-fake cluster 
preceding $\ell$ but not $\ell'_{m}$ with all possible left-fake 
clusters. Set also $\tau_{3}(\theta_{m-1},\ell_{m}')=\theta_{m-1}$. 

Note that, by construction,
\begin{equation}\label{eq:a.18} 
\begin{aligned} 
&\calA_{\ell_{m}}(\theta_{m},x_{\ell_{m}}) 
  \calG^{[n_{m}]}_{e_{\ell_{m}},u_{\ell_{m}}}(x_{\ell_{m}}) \Val(S_{m-1})= 
\calA_{\ell_{m-1}}(\theta_{m-1},x_{\ell_{m-1}}),\\ 
&\Val(S'_{m-1}) \calG^{[n'_{m}]}_{e_{\ell'_{m}},u_{\ell'_{m}}}(x_{\ell_{m}}) 
\BB_{\ell_{m}}(\theta_{m}) 
=\BB_{\ell'_{T_{m-1}}}(\theta_{m-1}) .
\end{aligned} 
\end{equation} 
One among cases {\bf 1}--{\bf 4} holds for $\ell_{m}\in L(\theta_{m})$,
so that we can consider the contribution to $\partial_{\ell_{m}}\Val 
(\theta_{m})$ (together with other contributions as in
{\bf 3} and {\bf 4}, if necessary) given by -- see (\ref{eq:a.10}), (\ref{eq:a.11}) and
(\ref{eq:a.13}) --
$$ 
\calA_{\ell_{m}}(\theta_{m},x_{\ell_{m}})\Big(\calG^{[n_{m}]}(x_{\ell_{m}}) 
\partial_{\ell_{m-1}}\Val_{R_{m}}(x_{\ell_{m}})\calG^{[n'_{m}]}(x_{\ell_{m}})
\Big)_{e_{\ell_{m}},u_{\ell'_{m}}}
\BB_{\ell_{m}} (\theta_{m}). 
$$ 
Then one has 
\begin{eqnarray} 
& & \calA_{\ell_{m}}(\theta_{m},x_{\ell_{m}})\Big(\calG^{[n_{m}]}(x_{\ell_{m}}) 
\partial_{\ell_{m-1}}\Val_{R_{m}}(x_{\ell_{m}})\calG^{[n'_{m}]}(x_{\ell_{m}}) 
\Big)_{e_{\ell_{m}},u_{\ell'_{m}}} \!\!\!\!\!\!\!\!
\BB_{\ell_{m}}(\theta_{m}) + \partial_{\ell_{m-1}}\!\!\!\!\!\!\!\!\!\!\! 
\sum_{\theta'\in \tau_{3}(\theta_{m-1},\ell_{m-1})} 
\!\!\!\!\!\!\!\!\!\!\! 
\Val(\theta';\e,\be_{0}) \nonumber \\ 
& & =\calA_{\ell_{m-1}}(\theta_{m-1},x_{\ell_{m-1}})
\Big(\partial_{\be_{0}} \calG^{[n+m-1]}(x_{\ell_{m-1}}) \times
\label{eq:a.19} \\
& & \qquad\times
\left(\uno+M^{[n+m-1]}(x_{\ell_{m-1}})
\calG^{[n+m]}(x_{\ell_{m-1}})\right) \Big)_{e,u'} \;
\BB_{\ell'_{T_{m-1}}}\!\!\!\!\!\!\!\!  (\theta_{m-1}), \nonumber
\end{eqnarray} 
where we have shortened $e,u'=e_{\ell_{m-1}},u_{\ell'_{T_{m-1}}}$
to simplify notation. By reasoning as in (\ref{eq:a.13}), this gives
\begin{eqnarray}
& & \calA_{\ell_{m-1}}(\theta_{m-1},x_{\ell_{m-1}}) 
 \Big(\calG^{[n+m-1]}(x_{\ell_{m-1}})\partial_{\be_{0}}\MM^{[n+m-2]}(x_{\ell_{m-1}})
\calG^{[n+m-1]}(x_{\ell_{m-1}})\Big)_{e,u'} \; \BB_{\ell'_{T_{m-1}}}\!\!\!\!\!\!\!\! (\theta_{m-1}) \nonumber \\ 
& & \quad + \; \calA_{\ell_{m-1}}(\theta_{m-1},x_{\ell_{m-1}}) 
\Big(\calG^{[n+m-1]}(x_{\ell_{m-1}})\partial_{\be_{0}}\MM^{[n+m-2]} (x_{\ell_{m-1}}) 
\calG^{[n+m]}(x_{\ell_{m-1}})\Big)_{e,u'} \; \BB_{\ell'_{T_{m-1}}}\!\!\!\!\!\!\!\! (\theta_{m-1}) \nonumber \\ 
& & \quad + \; \calA_{\ell_{m-1}}(\theta_{m-1},x_{\ell_{m-1}}) 
\Big(\calG^{[n+m-1]}(x_{\ell_{m-1}})\partial_{\be_{0}}\MM^{[n+m-2]}(x_{\ell_{m-1}}) 
\calG^{[n+m-1]}(x_{\ell_{m-1}}) \times
\label{eq:a.20}  \\ 
& & \qquad\qquad\times \;
M^{[n+m-1]}(x_{\ell_{m-1}}) \calG^{[n+m]}(x_{\ell_{m-1}})
\Big)_{e,u'} \; \BB_{\ell'_{T_{m-1}}}\!\!\!\!\!\!\!\! (\theta_{m-1}). \nonumber
\end{eqnarray} 
where again $e,u'=e_{\ell_{m-1}},u_{\ell'_{T_{m-1}}}$.

\begin{figure}[ht] 
\vskip-.3truecm 
\centering 
\ins{065pt}{-49pt}{$n\!+\!i$} 
\ins{070pt}{-64pt}{$\ell_{i}$} 
\ins{100pt}{-56pt}{$\le\!\! n\!\!+\!\!i\!\!-\!\!1$} 
\ins{128pt}{-80pt}{$S_{i-1}$} 
\ins{140pt}{-26pt}{$R_{i}$} 
\ins{145pt}{-64pt}{$\ell_{i-1}$} 
\ins{180pt}{-56pt}{$\le\!\! n\!\!+\!\!i\!\!-\!\!1$} 
\ins{170pt}{-80pt}{$S'_{i-1}$} 
\ins{236pt}{-64pt}{$\ell'_{i}$} 
\ins{230pt}{-49pt}{$n\!+\!i$} 
\ins{266pt}{-54pt}{$n\!+\!i$} 
\ins{300pt}{-49pt}{$n\!+\!i\!+\!1$} 
\ins{310pt}{-64pt}{$\ell'_{_{T_{i}}}$} 
\ins{384pt}{-56pt}{$+$} 
\ins{065pt}{-153pt}{$n\!+\!i$} 
\ins{070pt}{-167pt}{$\ell_{i}$} 
\ins{101pt}{-158pt}{$\le\!\! n\!\!+\!\!i\!\!-\!\!1$} 
\ins{128pt}{-184pt}{$S_{i-1}$} 
\ins{145pt}{-167pt}{$\ell_{i-1}$} 
\ins{180pt}{-161pt}{$n\!+\!i\!-\!1$} 
\ins{230pt}{-153pt}{$n\!+\!i$} 
\ins{230pt}{-167pt}{$\ell'_{T_{i-1}}$} 
\ins{260pt}{-158pt}{$\le\!\! n\!\!+\!\!i\!\!-\!\!1$} 
\ins{287pt}{-184pt}{$S'_{i-1}$} 
\ins{305pt}{-153pt}{$n\!+\!i$} 
\ins{310pt}{-167pt}{$\ell'_{i}$} 
\ins{344pt}{-158pt}{$n+i$} 
\ins{382pt}{-153pt}{$n\!+\!i\!+\!1$} 
\ins{390pt}{-167pt}{$\ell'_{_{T_{i}}}$} 
\includegraphics[width=5.4in]{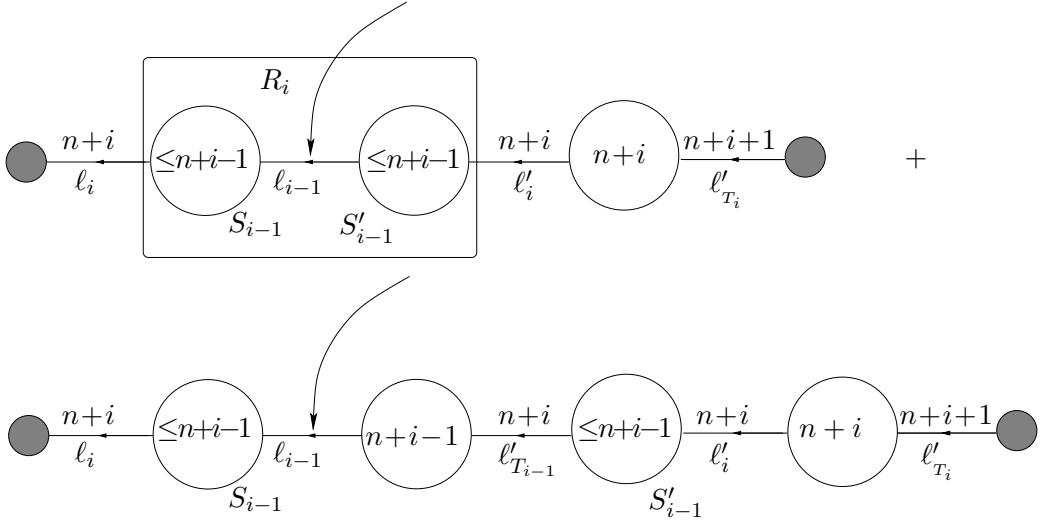} 
\vskip.3truecm 
\caption{Graphical representation of the left hand side of (\ref{eq:a.23}).} 
\label{fig:a8} 
\end{figure}

Then, for $i=m-1,\ldots,1$ we recursively reason as follows. Set 
$$ 
\BB_{\ell'_{T_{i}}}(\tau_{3}(\theta_{i},\ell'_{i+1})):= \!\!\!\!\!\!\!\!
\sum_{\theta'\in\tau_{3}(\theta_{i},\ell'_{i+1})} 
\BB_{\ell'_{T_{i}}}(\theta')
$$ 
and note that 
\begin{eqnarray}
& & \calA_{\ell_{i}}(\theta_{i},x_{\ell_{i}}) 
  \calG^{[n+i]}_{e_{\ell_{i}},u_{\ell_{i}}}(x_{\ell_{i}}) \Val(S_{i-1})= 
\calA_{\ell_{i-1}}(\theta_{i-1},x_{\ell_{i-1}}), 
\label{eq:a.21} \\ 
& & \Val(S'_{i-1}) \Big(\calG^{[n+i]}(x_{\ell_{i}}) 
M^{[n+i]}(x_{\ell_{i}}) \calG^{[n+i+1]}(x_{\ell_{i}})
\Big)_{e_{\ell'_{i-1}},u_{\ell'_{T_{i}}}}
\BB_{\ell'_{T_{i}}}(\tau_{3}(\theta_{i},\ell'_{i+1}) )
=\BB_{\ell'_{T_{i-1}}}(\tau_{3}(\theta_{i-1},\ell'_{i})). \nonumber
\end{eqnarray} 
Consider the contribution 
\begin{equation}\label{eq:a.22} 
\begin{aligned} 
&\calA_{\ell_{i}}(\theta_{i},x_{\ell_{i}}) 
  \Big(\calG^{[n+i]}(x_{\ell_{i}}) 
  \partial_{\ell_{i-1}}\Val_{R_{i}}(x_{\ell_{i}}) 
  \calG^{[n+i]}(x_{\ell_{i}}) \times \\
&\qquad\times
M^{[n+i]}(x_{\ell_{i}})  \calG^{[n+i+1]}(x_{\ell_{i}}) 
\Big)_{e_{\ell_{i}},u_{\ell'_{T_{i}}}}
\BB_{\ell'_{T_{i}}}(\tau_{3}(\theta_{i},\ell'_{i+1}) ),
\end{aligned} 
\end{equation}
obtained at the $(i+1)$-th step of the recursion.
By (\ref{eq:a.21}) one has (see Figure \ref{fig:a8}) 
\begin{equation}\label{eq:a.23} 
\begin{aligned} 
&\calA_{\ell_{i}}(\theta_{i},x_{\ell_{i}}) 
  \Big(\calG^{[n+i]}(x_{\ell_{i}}) 
  \partial_{\ell_{i-1}}\Val_{R_{i}}(x_{\ell_{i}}) 
  \calG^{[n+i]}(x_{\ell_{i}}) \times \\
&\qquad\qquad\times
M^{[n+i]}(x_{\ell_{i}})  \calG^{[n+i+1]}(x_{\ell_{i}}) 
\Big)_{e_{\ell_{i}},u_{\ell'_{T_{i}}}}
\BB_{\ell'_{T_{i}}}(\tau_{3}(\theta_{i},\ell'_{i+1}) )
+ 
\partial_{\ell_{i-1}}\!\!\!\!\!\!\!\!\!\!\! 
\sum_{\theta'\in \tau_{3}(\theta_{i-1},\ell_{i-1})} 
\!\!\!\!\!\!\!\!\!\!\!  \Val(\theta')\\
&\qquad
=\calA_{\ell_{i-1}}(\theta_{i-1},x_{\ell_{i-1}}) \, 
\Big(\partial_{\be_{0}}\calG^{[n+i-1]}(x_{\ell_{i-1}}) \times
\\ &\qquad\qquad \times \left(\uno+ 
M^{[n+i-1]}(x_{\ell_{i-1}})\calG^{[n+i]}(x_{\ell_{i-1}})\right) 
\Big)_{e_{\ell_{i-1}},u_{\ell'_{T_{i-1}}}}
\!\!\!\!\!\!\!\! \BB_{\ell'_{T_{i-1}}}(\tau_{3}(\theta_{i-1},\ell'_{i})), 
\end{aligned} 
\end{equation} 
which produces, as in (\ref{eq:a.20}), the contribution 
\begin{equation}\label{eq:a.24} 
\begin{aligned} 
\calA_{\ell_{i-1}}(\theta_{i-1},x_{\ell_{i-1}}) 
  & \Big(\calG^{[n+i-1]}(x_{\ell_{i-1}}) 
  \partial_{\ell_{i-2}}\Val_{R_{i-1}}(x_{\ell_{i-1}}) 
  \calG^{[n+i-1]}(x_{\ell_{i-1}}) \; \times
\\ &\times 
  M^{[n+i-1]}(x_{\ell_{i-1}})\calG^{[n+i]}(x_{\ell_{i-1}})
\Big)_{e_{\ell_{i-1}},u_{\ell'_{T_{i-1}}}}
\!\!\!\!\!\!\!\!  \BB_{\ell'_{T_{i-1}}}(\tau_{3}(\theta_{i-1},\ell'_{i})). 
\end{aligned} 
\end{equation} 
%

Hence we can proceed recursively from $\theta_{m}$ up to $\theta_{0}$, 
until we obtain 
\begin{eqnarray}
& & \calA_{\ell_{0}}(\theta_{0},x_{\ell_{0}}) 
 \Big(\calG^{[n]}(x_{\ell_{0}})\partial_{\be_{0}}\MM^{[n-1]}(x_{\ell_{0}}) 
 \calG^{[n]}(x_{\ell_{0}})  \Big)_{e_{\ell_{0}},u_{\ell'_{T_{0}}}}
\!\!\!\!\!\!
  \BB_{\ell'_{T_{0}}}(\tau_{3}(\theta_{0},\ell'_{1})) \nonumber \\ 
& & \qquad+\;\calA_{\ell_{0}}(\theta_{0},x_{\ell_{0}}) \Big(
 \calG^{[n]}(x_{\ell_{0}})\partial_{\be_{0}}\MM^{[n-1]}(x_{\ell-{0}}) 
 \calG^{[n+1]}(x_{\ell_{0}})  \Big)_{e_{\ell_{0}},u_{\ell'_{T_{0}}}}
\!\!\!\!\!\!
 \BB_{\ell'_{T_{0}}}(\tau_{3}(\theta_{0},\ell'_{1}))
\label{eq:a.25} \\ 
& & \qquad+\;\calA_{\ell_{0}}(\theta_{0},x_{\ell_{0}}) \Big(
 \calG^{[n]}(x_{\ell_{0}})\partial_{\be_{0}}\MM^{[n-1]}(x_{\ell_{0}}) 
 \calG^{[n]}(x_{\ell_{0}})M^{[n]}(x_{\ell_{0}}) \calG^{[n+1]}(x_{\ell_{0}}) 
 \Big)_{e_{\ell_{0}},u_{\ell'_{T_{0}}}}
\!\!\!\!\!\!\!\! \BB_{\ell'_{T_{0}}}(\tau_{3}(\theta_{0},\ell'_{1})). \nonumber
\end{eqnarray} 
Once again, if we split $\partial_{\be_{0}}=\partial_{N}+\partial_{L}$, 
all the terms with $\partial_{N}\MM^{[n-1]}$ are contributions to 
$\MM^{[p]}_{\be,\be}(0)$.\qed

\noindent
{\bf 6.} We are left with the derivatives
$\partial_{L}M^{[q]}(x)$, $q\le n$, 
when the differentiated propagator is not one of those used along 
the cases {\bf 4} or {\bf 5}; see for instance (\ref{eq:a.16}),
(\ref{eq:a.19}) and (\ref{eq:a.23}). One can reason 
as in the case $\partial_{L}\Val(\theta)$, by studying the derivatives 
$\partial_{\ell}\Val_{T}(x_{\ell})$ and proceed iteratively 
along the lines of cases {\bf 1} to {\bf 5} above, until 
only lines on scales $0$ are left. In that case the derivatives 
$\partial_{\be_{0}}\calG^{[0]}(x_{\ell})$ 
produce derivatives
$$\partial_{\be_{0}}M^{[-1]}(x)=\begin{pmatrix}
\e\partial_{\be_{0}}^2F_{\vzero}(\beta_{0},B_{0}) &
\e\partial_{\be_{0},B_{0}}^2F_{\vzero}(\beta_{0},B_{0}) \cr
& \cr
\e\partial_{\be_{0}}^2G_{\vzero}(\beta_{0},B_{0}) &
\e\partial_{\be_{0},B_{0}}^2G_{\vzero}(\beta_{0},B_{0}) \cr
\end{pmatrix}
$$
(see Remarks \ref{rmk:3.5} and \ref{rmk:3.6}). Therefore, for $n=-1$, 
in the splitting (\ref{eq:a.9}), there are no terms with the 
derivatives $\partial_{\ell}$ and the derivatives $\partial_{v}$ can 
be interpreted as said in Remark \ref{rmk:a.1}.\qed

\noindent
{\bf 7.} By construction, each contribution to $\MM^{[p-1]}_{\be,\be}(0)$
appears as one term among those considered in the discussion above,
that is among the contributions to $\partial_{\be_{0}}\Phi_{\vzero}^{\RR,p}(\e,\be_{0},B_{0})$
arising from the trees $\theta\in\Theta^{\RR,{p}}_{k,\vzero,\be}$ satisfying the condition (\ref{eq:a.0bis}).
Of course, when computing $\partial_{\be_{0}}\Val(\theta)$ for such trees,
also some contributions to $M^{[p]}_{\be,\be}(0)$ have been produced.
Call $W^{[p]}$ the contributions to $M^{[p]}_{\be,\be}(0)$
which are not obtained in the previous steps.
Define also $R^{[p]}$ as the sum of the contributions to $\partial_{\be_{0}}\Phi^{\RR,p}_{\vzero}$ such that
\begin{equation} \label{eq:a.27}
\partial_{\be_{0}}\Phi^{\RR,{p-1}}_{\vzero}+
R^{[p]}=\MM^{[p-1]}_{\be,\be}(0)+ \left( M^{[p]}_{\be,\be}(0)-W^{[p]} \right) ,
\end{equation}
where we have used that $\MM^{[p]}_{\be,\be}(0)=\MM^{[p-1]}_{\be,\be}(0)+
M^{[p]}_{\be,\be}(0)$ -- see definition (\ref{eq:3.8}) and use that
$\chi_{q}(0)=1$ for all $q\ge-1$. Hence $\partial_{\be_{0}}\Phi^{\RR,{p-1}}_{\vzero}+R^{[p]}$ represents
the sum of all contributions to $\partial_{\be_{0}}\Phi^{\RR,p}_{\vzero}$ used in {\bf 1}--{\bf 6}.
One can write
\begin{equation} \label{eq:a.28}
\partial_{\be_{0}}\Phi^{\RR,{p}}_{\vzero}=
\partial_{\be_{0}}\Phi^{\RR,{p-1}}_{\vzero}+R^{[p]}+S^{[p]} ,
\end{equation}
for a suitable $S^{[p]}$: by construction $S^{[p]}$ takes into account all contributions
arising from the trees $\theta\in\Theta^{\RR,{p}}_{k,\vzero,\be}$
which do not satisfy the condition (\ref{eq:a.0bis}), i.e. such that
\begin{equation} \label{eq:a.29} 
\max_{\ell\in \Theta^{\RR,{p}}_{k,\vzero,\be}}
\{ n \in \ZZZ_{+} : \Psi_{n}(\oo\cdot\nn_{\ell}) \neq 0 \} = p+1 .
\end{equation}
Such trees have been excluded in the discussion above, because on
the one hand they would produce the remaining contributions to $M^{[p]}_{\be,\be}(0)$,
on the other hand  they would equally produce contributions to $M^{[p+1]}_{\be,\be}(0)$.
Therefore, by combining (\ref{eq:a.27}) and (\ref{eq:a.28}), we obtain
$\partial_{\be_{0}}\Phi^{\RR,{p}}_{\vzero}$ $=$ 
$\MM^{[p]}_{\be,\be}(0)$ $+$ $\left( S^{[p]}-W^{[p]} \right)$,
where both $W^{[p]}$ and $S^{[p]}$ arise from trees containing at
least one line $\ell$ on scale $p$ and such that
$\Psi_{p+1}(\oo\cdot\nn_{\ell})\neq0$:
for such a line $\ell$ one has $|\nn_{\ell}| > 2^{m_{p+1}-1}$ by Remark \ref{rmk:3.8}.
Therefore, one has
$\max\big\{ |S^{[p]}|, |W^{[p]}| \big\} \le |\e|\, D_{1}e^{-D_{2}2^{m_{p+1}}}$,
for some constants $D_{1},D_{2}$ and this is enough to prove the bound (\ref{eq:a.0}).\qed

\zerarcounters \section{Proof of Lemmas \ref{lem:6.3} and \ref{lem:6.4}} 
\label{app:b} 

We want to prove by induction that
\begin{equation} \label{eq:b.1}
\gotN^{\bullet}_{n}(\theta) \le \max\{ 2^{-(m_{n}-3)} K(\theta) - 2 , 0 \} .
\end{equation}
First of all note that $\gotN^{\bullet}_{n}(\theta) \ge 1$ implies
$\gotN_{n}^{*}(\theta)\ge 1$ and hence $K(\theta) \ge 2^{m_{n}-1}$.

Set $\ze_{0}:=\ze_{\ell_{\theta}}$, $n_{0}:=n_{\ell_{\theta}}$
and $\nn:=\nn_{\ell_{\theta}}$,
and note that either $n_{0}=\ze_{0}$ or $n_{0}=\ze_{0}+1$.
If $\ze_{0}<n$ the bound (\ref{eq:b.1}) follows from the inductive hypothesis.
If $\ze_{0} \ge n$, call $\ell_{1},\ldots,\ell_{r}$ the lines with minimum scale
$\ge n$ closest to $\ell_{\theta}$ and $\theta_{1},\ldots,\theta_{r}$
the subtrees with root lines $\ell_{1},\ldots,\ell_{r}$, respectively.
If $r=0$ the bound trivially holds. If $r\ge 2$, by the inductive hypothesis
one has
$\gotN_{n}^{\bullet}(\theta)=1+\gotN_{n}^{\bullet}(\theta_{1})+\ldots+
\gotN_{n}^{\bullet}(\theta_{r})
\le 1 + 2^{-(m_{n}-3)}K(\theta)-2r \le 2^{-(m_{n}-3)}K(\theta)-3$, so that
the bound follows once more.

If $r=1$ call $T$ the subgraph with exiting line $\ell_{\theta}$ and
entering line $\ell_{1}$. Then either $T$ is a
self-energy cluster or $K(T) \ge 2^{m_{n}-1}$. This can be proved as follows.
Set $\nn_{1}:=\nn_{\ell_{1}}$. If $T$ is not a cluster, then it must contain
at least one line
$\ell'$ on scale $n_{\ell'}=n$, so that if $\ell'\notin\calP_{T}$ and
$\theta'$ is the subtree with root line $\ell'$ one has 
$K(T) \ge K(\theta') \ge 2^{m_{n}-1}$, while if $\ell\in\calP_{T}$ then
$\nn_{\ell'} \neq \nn_{1}$ (because $\ze_{\ell'}=n-1$ and
$\ze_{\ell_{1}} \ge n$), so that $|\oo\cdot(\nn_{\ell'}-\nn_{1})| < \al_{m_{n}-1}(\oo)$
implies $K(T) \ge |\nn_{\ell'}-\nn_{1}| \ge 2^{m_{n}-1}$. If $T$ is a cluster
then either
(i) $\nn_{1}\ne \nn$ so that
$K(T) \ge |\nn-\nn_{1}| \ge 2^{m_{n}-1}$ or (ii) $\nn_{1}=\nn$ and there is a
line $\ell\in\calP_{T}$ with $n_{\ell}=-1$ so that $K(T)\ge|\nn_{\ell}^{0}|=
|\nn_{1}|\ge 2^{m_{n}-1}$, or (iii) $\nn_{1}=\nn$ and $T$ is a self-energy
cluster, otherwise there would be a line
$\ell'\in\calP_{T}$ with $\nn_{\ell'}=\nn_{1}$, which is incompatible
with $\ze_{\ell'}\le n-1$ and $\ze_{\ell_{1}}\ge n$.

Therefore, if $K(T) \ge 2^{m_{n}-1}$, the inductive
hypothesis yields the bound (\ref{eq:b.1}).
If $K(T) < 2^{m_{n}-1}$ then $T$ is a self-energy cluster (and hence $\nn_{1}=
\nn$).
In such a case call $\theta_{1}$ the tree with root line $\ell_{1}$;
by construction $\gotN^{\bullet}_{n}(\theta)=1+\gotN^{\bullet}(\theta_{1})$.
We can repeat the argument above: call $\ell'_{1},\ldots,\ell'_{r'}$
the lines with minimum scale
$\ge n$ closest to $\ell_{1}$ and $\theta'_{1},\ldots,\theta'_{r'}$
the subtrees with root lines $\ell'_{1},\ldots,\ell'_{r'}$, respectively.
Again the case $r'=0$ is trivial.
If $r'\ge2$ then
$\gotN_{n}^{\bullet}(\theta)=2+\gotN_{n}^{\bullet}(\theta'_{1})+
\ldots+\gotN_{n}^{\bullet}(\theta'_{r'})
\le 2 + 2^{-(m_{n}-3)}K(\theta)-2r \le 2^{-(m_{n}-3)}K(\theta)-2$, so
yielding the bound. Therefore the only case which does not
imply immediately the bound (\ref{eq:b.1}) through the inductive hypothesis
is when $\ell_{1}$ exits a subgraph $T'$ with only one entering
line $\ell_{1}'$ on minimum scale $\ge n$. Set $\nn_{1}':=\nn_{\ell_{1}'}$ and
call $\theta_{1}'$ the tree with root line $\ell_{1}'$. As before we have that
either $K(T') \ge 2^{m_{n}-1}$ or $T'$ is a self-energy cluster.
If $K(T') \ge 2^{m_{n}-1}$ then $\gotN^{\bullet}_{n}(\theta)=2+\gotN^{\bullet}_{n}(\theta_{1}')$
and one can reason as before to obtain the bound
by relying on the inductive hypothesis. If $T'$ is a self-energy cluster
then $\ell_{1}$ is a resonant line and $\gotN^{\bullet}_{n}(\theta)=1+\gotN^{\bullet}_{n}(\theta_{1}')$.

One can iterate again the argument until either one reaches a case which
can be dealt with through the inductive hypothesis or one obtains
$\gotN^{\bullet}_{n}(\theta)=1+\gotN^{\bullet}_{n}(\theta'')$, for some
tree $\theta''$ which has no line $\ell$ with $\ze_{\ell}\ge n$. Thus
$\gotN^{\bullet}_{n}(\theta'')=0$ and the bound (\ref{eq:b.1}) follows.
This completes the proof of Lemma \ref{lem:6.3}. 

Now, we pass to the proof of Lemma \ref{lem:6.4}.
Consider a self-energy cluster $T\in\gotS^{k}_{n,u,e}$. First of
all we prove that $K(T) \ge 2^{m_{n}-1}$: the argument proceeds as in
the proof of Lemma \ref{lem:3.10}, by noting that if there is a line
$\ell\in\calP_{T}$ on scale $n$ then $\nn_{\ell} \neq \nn_{\ell_{T}'}$
(otherwise $\nn_{\ell}^{0}=0$ and hence $T$ would not be a self-energy cluster).

Define $\CCCC(n,p)$ as the set of subgraphs $T$ of $\theta$ with
only one entering line $\ell_{T}'$ and one exiting line $\ell_{T}$
both on minimum scale $\ge p$, such that $L(T)\neq \emptyset$ and $n_{\ell}
\le n$ for any line $\ell\in L(T)$. We prove by induction on the order the
bound
\begin{equation} \label{eq:b.2}
\gotN^{\bullet}_{p}(T)  \le 2^{-(m_{p}-3)}K(T) 
\end{equation}
for all $T\in \CCCC(n,p)$ and all $0\le p \le n$.
Consider $T\in\CCCC(n,p)$, $p\le n$: 
call $\ell_{1},\ldots,\ell_{r}$ the lines with minimum scale
$\ge p$ closest to $\ell_{T}$. The case $r=0$ is trivial.
If $r\ge 1$ and none of such lines is along the path $\calP_{T}$ then the
bound follows
from (\ref{eq:b.1}). If one of such lines, say $\ell_{1}$, is along
the path $\calP_{T}$, then denote by $\theta_{2},\ldots,\theta_{r}$ the subtrees
with root lines $\ell_{2},\ldots,\ell_{r}$, respectively, and by
$T_{1}$ the subgraph with exiting line $\ell_{1}$ and entering line $\ell_{T}'$.
One has $\gotN_{p}^{\bullet}(T) \le 1 + \gotN_{p}^{\bullet}(T_{1}) +
\gotN_{p}^{\bullet}(\theta_{2})+\ldots+\gotN_{p}^{\bullet}(\theta_{r})$.
By construction $T_{1} \in \CCCC(n,p)$, so that 
the bound (\ref{eq:b.2}) follows by the inductive hypothesis for $r\ge 2$.

If $r=1$ then call $T_{0}$ the subgraph with exiting line $\ell_{T}$
and entering line $\ell_{1}$. By reasoning as in the proof of Lemma
\ref{lem:6.3} we find that either $K(T_{0})\ge 2^{m_{p}-1}$ or $T_{0}$ is a
self-energy cluster.
Since $\gotN_{p}^{\bullet}(T) \le 1 + \gotN_{p}^{\bullet}(T_{1})$, if $K(T_{0})\ge 2^{m_{p}-1}$
the bound follows once more. If on the contrary $T_{0}$ is a self-energy cluster
we can iterate the construction: call $\ell'_{1},\ldots,\ell'_{r'}$ the lines with minimum scale
$\ge p$ closest to $\ell_{1}$. If either $r'=0$ or no line among
$\ell'_{1},\ldots,\ell'_{r'}$ is along the path $\calP_{T_{1}}$, the bound follows easily.
Otherwise if a line, say $\ell'_{1}$ is along the path $\calP_{T_{1}}$ and $r'\ge 2$
one has $\gotN_{p}^{\bullet}(T) \le 2 + \gotN_{p}^{\bullet}(T'_{1}) +
\gotN_{p}^{\bullet}(\theta'_{2})+\ldots+\gotN_{p}^{\bullet}(\theta'_{r'})$, where
$T_{1}'$ is the subgraph with exiting line $\ell_{1}'$ and entering
line $\ell'_{T}$ and hence $\gotN_{p}^{\bullet}(T) \le 2 + 2^{-(m_{p}-3)}K(T)-2$,
by the inductive hypothesis, so that (\ref{eq:b.2}) follows.

If $r'=1$ let $T_{0}'$ be the subgraph with exiting line $\ell_{1}$
and entering line $\ell'_{1}$.
If $K(T_{0}') \ge 2^{m_{p}-1}$ then the inductive hypothesis implies
once more the bound (\ref{eq:b.2}), while if $K(T_{0}') <
2^{m_{p}-1}$ then, by the same argument as above, $T_{0}'$ must be a
self-energy cluster, so that $\ell_{1}$ does not contribute to
$\gotN_{p}^{\bullet}(T)$,
i.e. $\gotN_{p}^{\bullet}(T)\le1+\gotN_{p}(T''_{1})$ where $T''_{1}$ is the
subgraph with exiting line $\ell_{1}'$ and entering line $\ell_{T}'$.
Again we can
iterate the argument until either one finds a subgraph $T''$ with
$K(T'')\ge 2^{m_{p}-1}$, so that the inductive hypothesis compels the
bound (\ref{eq:b.2}) for $T$, or one obtains
$\gotN_{p}^{\bullet}(T)\le1+\gotN_{p}(T'')$
for some subgraph $T''$ which has no line on minimum scale $\ge p$, so that
$\gotN_{p}^{\bullet}(T)\le1$.

\zerarcounters \section{Proof of Lemma \ref{lem:6.9}} 
\label{app:c} 

We say that a self-energy cluster $T$ is \emph{isolated} if both its
external lines
are non-resonant and that is \emph{relevant} if it is not isolated.
As will emerge from the proof, it is convenient to introduce a further label
$\gotd_{T}\in\{0,1\}$ to be associated with each \resonance $T$.
We shall see later how to fix such a label: for the time being we consider it
as an abstract label and we define the \emph{subchains} as follows.
Given $p\ge2$ \resonances $T_{1},\ldots,T_{p}$ of a tree $\theta$,
with $\ell'_{T_{i}}=\ell_{T_{i+1}}$ for all $i=1,\ldots,p$, we say
that $C=\{T_{1},\ldots,T_{p}\}$ is a subchain if $\gotd_{T_{i}}=1$ for $i=1,
\ldots,p$, the line $\ell_{T_{1}}$ either is non-resonant or enters a
\resonance $T_{0}$ with $\gotd_{T_{0}}=0$ and the line
$\ell_{T_{p}}'$ either is non-resonant or exits
a \resonance $T_{p+1}$ with $\gotd_{T_{p+1}}=0$.
We say that a \resonance $T$ is a \emph{link} if $\gotd_{T}=1$.

Given a subchain $C=\{T_{1},\ldots,T_{p}\}$ of a tree $\theta$,
the \resonances $T_{i}$ are called the \emph{links} of $C$.
Define $\ell_{0}(C):=\ell_{T_{1}}$ and $\ell_{i}(C):=\ell_{T_{i}}'$
for $i=1,\ldots,p$ and set $n_{i}(C)=n_{\ell_{i}(C)}$ for $i=0,\ldots,p$.
The lines $\ell_{0}(C),\ldots,\ell_{p}(C)$ are the \emph{chain-lines} of $C$:
we call $\ell_{1}(C),\ldots,\ell_{p-1}(C)$ the \emph{internal chain-lines} of $C$
and $\ell_{0}(C),\ell_{p}(C)$ the \emph{external chain-lines} of $C$.
For future convenience we also set $\ell_{C}=\ell_{0}(C)$ and $\ell_{C}'=\ell_{p}(C)$.
We also  call $k(C):=k(T_{1})+\ldots+k(T_{p})$ the \emph{total order}
of the subchain $C$ and $p(C)=p$ the \emph{length} of $C$.
The value of a subchain $C$ is defined as in (\ref{eq:6.6}). 
Note that for all $i=1,\ldots,p-1$ one has
$\ze_{\ell_{i}(C)}=\ze_{\ell_{C}}=\ze_{\ell'_{C}}$ if $\Val(\theta)\neq0$.

We denote by $\gotC_{1}(k;h,h';n_{0},\ldots,n_{p})$ the set of all subchains $C=\{T_{1},\ldots,T_{p}\}$
with total order $k$ and with fixed labels $h_{\ell_{0}(C)}=h$,
$h_{\ell_{p}(C)}=h'$ and $n_{i}(C)=n_{i}$ for $i=0,\ldots,p$.

\begin{figure}[ht] 
\centering 
\ins{032pt}{-13pt}{$\ol{n}_{0},h_{0}$} 
\ins{042pt}{-30pt}{$\ell_{0}$} 
\ins{116pt}{-46pt}{$T_{1}$} 
\ins{084pt}{-22pt}{$n_{1},k_{1}$} 
\ins{130pt}{-13pt}{$\ol{n}_{1},h_{1}$} 
\ins{140pt}{-30pt}{$\ell_{1}$} 
\ins{214pt}{-46pt}{$T_{2}$} 
\ins{178pt}{-22pt}{$n_{2},k_{2}$} 
\ins{222pt}{-13pt}{$\ol{n}_{2},h_{2}$} 
\ins{234pt}{-30pt}{$\ell_{2}$} 
\ins{266pt}{-26pt}{$\ldots\ldots\ldots$}
\ins{316pt}{-13pt}{$\ol{n}_{p\!-\!1},h_{p\!-\!1}$} 
\ins{332pt}{-30pt}{$\ell_{p\!-\!1}$} 
\ins{407pt}{-46pt}{$T_{p}$} 
\ins{372pt}{-22pt}{$n_{p},k_{p}$} 
\ins{420pt}{-13pt}{$\ol{n}_{p},h_{p}$} 
\ins{430pt}{-30pt}{$\ell_{p}$} 
\includegraphics[width=6.0in]{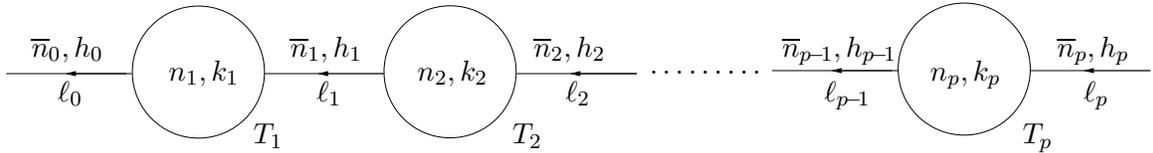} 
\vskip.3truecm 
\caption{A subchain $C$ of length $p$ with links $T_{1},\ldots,T_{p}$ and
chain-lines $\ell_{0},\ldots,\ell_{p}$; summing over all
possible $C$ with $h_{0}=h$, $h_{p}=h'$, $k_{1}+\ldots+k_{p}=k$
and $\ol{n}_{1},\ldots,\ol{n}_{p}$ fixed, one obtains
a graphical representation of (\ref{eq:c.1}).}
\label{fig:c1} 
\end{figure}

If all \resonances $T$ of $\theta$ carried a label $\gotd_{T}=1$ the definition of subchain 
would reduce to that of chain in Section \ref{sec:6}.  We want to
prove the bound (\ref{eq:6.7}). The sum is over all chains
$C=\{T_{1},\ldots,T_{p}\}$ in $\gotC(k;h,h';n_{0},\ldots,n_{p})$;
then we set $\gotd_{T_{i}}=1$ for $i=1,\ldots,p$, so that we can replace
$\gotC(k;h,h';n_{0},\ldots,n_{p})$ with $\gotC_{1}(k;h,h';n_{0},\ldots,n_{p})$.
Thus in (\ref{eq:6.7}) we can write
\begin{equation}\label{eq:c.1}
\sum_{C\in \gotC(k;h,h';\ol{n}_{0},\ldots,\ol{n}_{p})} \Val_{C}(x) =
\sum_{C\in \gotC_{1}(k;h,h';\ol{n}_{0},\ldots,\ol{n}_{p})} \Val_{C}(x) =
\sum_{\substack{h_{1},\ldots,h_{p-1}\in\{\be,B\} \\ k_{1}+\ldots+k_{p}=k}}\prod_{i=1}^{p}
\MM_{h_{i-1},h_{i}}^{(k_{i})}(x,n_{i}),
\end{equation}
where $h_{0}=h$, $h_{p}=h'$ and $n_{i}=\min\{\ol{n}_{i-1},\ol{n}_{i}\}-1$ for $i=1,\ldots,p$;
of course $|n_{i}-n_{j}|\le1$ for all $i,j=1,\ldots,p$ and $k_{i}\ge
0$ for $i=1,\ldots,p$; see Figure \ref{fig:c1}.

For all $k\ge1$, all $n\ge-1$ and all $h,h'\in\{\be,B\}$ let us write
\begin{equation}\label{eq:c.2}
\MM^{(k)}_{h,h'}(x,n)=\sum_{\de\in\De}\MM^{(k)}_{h,h'}(x,n,\de)
\end{equation}
where $\Delta:=\{\calL,\partial,\partial^{2},\RR\}$ is a set of labels and
\begin{equation}\label{eq:c.3}
\begin{aligned}
&\MM^{(k)}_{h,h'}(x,n,\calL):=\MM^{(k)}_{h,h'}(0)\qquad
\MM^{(k)}_{h,h'}(x,n,\partial):=x\partial\MM^{(k)}_{h,h'}(0), \\
&\MM^{(k)}_{h,h'}(x,n,\partial^{2}):=x^{2}\int_{0}^{1}{\rm d}\tau(1-\tau)
\partial^{2}\MM^{(k)}_{h,h'}(\tau x),\\
&\MM^{(k)}_{h,h'}(x,n,\RR):=\MM^{(k)}_{h,h'}(x,n)-\MM^{(k)}_{h,h'}(x),
\end{aligned}
\end{equation}
so that we can decompose the sum in (\ref{eq:c.1}) as
\begin{equation}\label{eq:c.4}
\sum_{\de_{1},\ldots,\de_{p}\in\De}
\sum_{\substack{h_{1},\ldots,h_{p-1}\in\{\be,B\} \\ k_{1}+\ldots+k_{p}=k}}
\prod_{i=1}^{p} \MM_{h_{i-1},h_{i}}^{(k_{i})}(x,n_{i},\de_{i}).
\end{equation}

There are several contributions to (\ref{eq:c.4}) which sum up to zero. 
This holds for all contributions with $\de_{j}=\de_{j+1}=\calL$ for
some $j=1,\ldots,p-1$. Indeed one can write such contributions as
\begin{eqnarray}
& & \sum_{\de_{1,\ldots,\de_{j-1},\de_{j+1},\ldots,\de_{p}}\in\Delta}
\sum_{\substack{h_{1},\ldots,h_{j-1},h_{j+1},\ldots h_{p-1} \in \{\be,B\} \\ k_{1}+\ldots+k_{j-1}+\ol{k}+k_{j+2}+\ldots+k_{p}=k}}
\prod_{\substack{i=1 \\ i\ne j}}^{p}
\MM^{(k_{i})}_{h_{i-1},h_{i}}(x,n_{i},\de_{i}) \; \times \nonumber \\
& & \times \; \Biggl(\sum_{k_{j}+k_{j+1}=\ol{k}}\MM^{(k_{j})}_{h_{j-1},\be}(0)
\MM_{\be,h_{j+1}}^{(k_{j+1})}(0)
+\sum_{k_{j}+k_{j+1}=\ol{k}}\MM^{(k_{j})}_{h_{j-1},B}(0)\MM_{B,h_{j+1}}^{(k_{j+1})}(0) \Biggr),
\label{eq:c.5} 
\end{eqnarray}
and by Lemma \ref{lem:6.5} one has (for instance)
\begin{equation*}
\begin{aligned}
\sum_{k_{j}+k_{j+1}=\ol{k}}&\MM^{(k_{j})}_{\be,\be}(0)\MM_{\be,\be}^{(k_{j+1})}(0)
+\sum_{k_{j}+k_{j+1}=\ol{k}}\MM^{(k_{j})}_{\be,B}(0)\MM_{B,\be}^{(k_{j+1})}(0) \\
&=\sum_{k_{j}+k_{j+1}=\ol{k}}\MM^{(k_{j})}_{\be,\be}(0)\MM_{\be,\be}^{(k_{j+1})}(0)
+\sum_{k_{j}+k_{j+1}=\ol{k}}\MM^{(k_{j})}_{\be,B}(0)
\Biggl(-\sum_{k'+k''=k_{j+1}}\MM_{BB}^{(k')}(0)\partial_{\be_{0}}
B_{0}^{(k'')}\Biggr)\\
&=\sum_{k_{j}+k_{j+1}=\ol{k}}\MM^{(k_{j})}_{\be,\be}(0)\MM_{\be,\be}^{(k_{j+1})}(0)
+\sum_{k_{j}+k_{j+1}=\ol{k}}\MM_{BB}^{(k_{j})}(0)\MM_{\be,\be}^{(k_{j+1})}(0)=0 ;
\end{aligned}
\end{equation*}
one can reason in the same way  also for the cases $(h_{j-1},h_{j+1})\neq(\be,\be)$.
By (\ref{eq:6.5a}) of Lemma \ref{lem:6.6}, also the contributions with $\de_{j}=\partial$
and $h_{j-1}\ne h_{j}$ for some $j=1,\ldots,p$ sum up zero. Finally
we obtain zero also when we sum together all contributions with $\de_{j}=\ldots=\de_{j+q}
=\partial$, $h_{j-1}=\ldots=h_{j+q}=\ol{h}$ for some $\ol{h}\in\{\be,B\}$ and
$\de_{j-1}=\de_{j+q+1}=\calL$ for some $j=2,\ldots,p-1$ and $q=0,\ldots,p-1-j$.
Indeed we can write the sum of such contributions as
\begin{equation*}
\begin{aligned}
&\sum_{\substack{h_{1},\ldots,h_{j-2},\ol{h},h_{j+q+1},\ldots,h_{p-1}=\be,B \\ k_{1}+\ldots+k_{p}=k}}
\Biggl( \prod_{i=1}^{j-2} \MM^{(k_{i})}_{h_{i-1},h_{i}}(x,n_{i},\de_{i}) 
\Biggr) 
 \MM^{(k_{j-1})}_{h_{j-2},\ol{h}}(x,n_{j-1},\calL) \times\\
&\qquad\qquad\qquad \times
\Biggl( \prod_{i=j}^{j+q} \MM^{(k_{i})}_{\ol{h},\ol{h}}(x,n_{i},\partial) \Biggr) 
 \MM^{(k_{j+q+1})}_{\ol{h},h_{j+q+1}}(x,n_{j+q+1},\calL)
\Biggl( \prod_{i=j+q+2}^{p}
\MM^{(k_{i})}_{h_{i-1},h_{i}}(x,n_{i},\de_{i}) \Biggr)
\end{aligned}
\end{equation*}
\begin{equation*}
\begin{aligned}
& \qquad = \sum_{k''\le k}
\sum_{\substack{k_{j}+\ldots+k_{j+q} =k''}}
\Biggl( \prod_{i=j}^{j+q} \MM^{(k_{i})}_{\be,\be}(x,n_{i},\partial) \Biggr) \times\\
&\qquad\qquad\qquad \times
\sum_{\substack{k_{1}+\ldots+k_{j-1}+k_{k+q+1}+\ldots+k_{p} = k - k'\\ h_{1},\ldots,h_{j-2},h_{j+q+1},\ldots,h_{p-1}=\be,B}}
\Biggl( \prod_{i=1}^{j-2} \MM^{(k_{i})}_{h_{i-1},h_{i}}(x,n_{i},\de_{i})  \Biggr) \times\\
&\qquad\qquad\qquad \times
\sum_{\ol{h}=\be,B}  \MM^{(k_{j-1})}_{h_{j-2},\ol{h}}(x,n_{j-1},\calL) \MM^{(k_{j+q+1})}_{\ol{h},h_{j+q+2}}
(x,n_{j+q+1},\calL)
 \Biggl( \prod_{i=j+q+2}^{p}
 \MM^{(k_{i})}_{h_{i-1},h_{i}}(x,n_{i},\de_{i}) \Biggr) ,
\end{aligned}
\end{equation*}
and the last sum is zero by the same argument used for (\ref{eq:c.5});
note that we used (\ref{eq:6.5b}) to extract a common factor
$\MM^{(k_{i})}_{\be,\be}(x,n_{i},\partial)$ in the third line.

We say that a cluster $T$ is a \emph{fake cluster} on scale $n$ if it is a connected subgraph
of a tree with only one entering line $\ell_{T}'$ and one exiting line $\ell_{T}$
such that (i) all lines in $T$ have scale $\le n$ and there is at least one line on $T$ 
with scale $n$ and (ii) the lines $\ell_{T}$ and $\ell_{T}'$ carry the same momentum;
note that a fake cluster can fail to be a self-energy cluster only because there the scales of the
external lines have no relation with $n$ (and hence it can even fail
to be a cluster). Denote by $\gotS^{*k}_{m,u,e}$ the set of fake clusters
with order $k$, scale $m$ and such that $h_{\ell_{T}'}=e$ and $h_{\ell_{T}}=u$.

\begin{figure}[ht] 
\centering 
\ins{032pt}{-13pt}{$\ol{n}_{0},h_{0}$} 
\ins{042pt}{-30pt}{$\ell_{0}$} 
\ins{116pt}{-46pt}{$T_{1}$} 
\ins{076pt}{-22pt}{$n_{1},k_{1},\de_{1}$} 
\ins{130pt}{-13pt}{$\ol{n}_{1},h_{1}$} 
\ins{140pt}{-30pt}{$\ell_{1}$} 
\ins{214pt}{-46pt}{$T_{2}$} 
\ins{171pt}{-22pt}{$n_{2},k_{2},\de_{2}$} 
\ins{224pt}{-13pt}{$\ol{n}_{2},h_{2}$} 
\ins{234pt}{-30pt}{$\ell_{2}$} 
\ins{266pt}{-26pt}{$\ldots\ldots\ldots$}
\ins{316pt}{-13pt}{$\ol{n}_{p\!-\!1},h_{p\!-\!1}$} 
\ins{332pt}{-30pt}{$\ell_{p\!-\!1}$} 
\ins{407pt}{-46pt}{$T_{p}$} 
\ins{365pt}{-22pt}{$n_{p},k_{p},\de_{p}$} 
\ins{420pt}{-13pt}{$\ol{n}_{p},h_{p}$} 
\ins{430pt}{-30pt}{$\ell_{p}$} 
\includegraphics[width=6.0in]{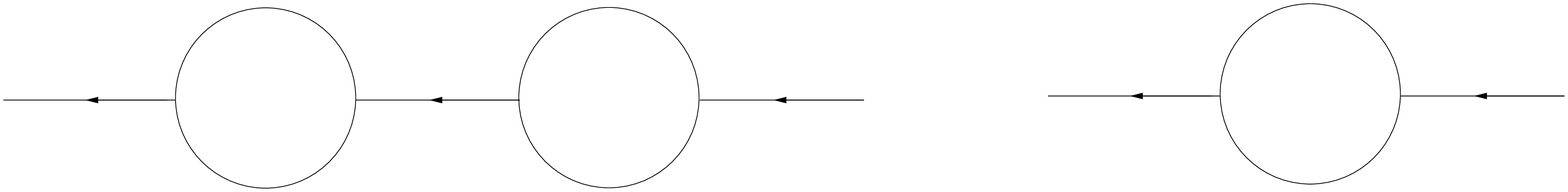} 
\vskip.3truecm 
\caption{A $*$-chain $C\in\gotC^{*}(k;h,h';\ol{n}_{0},\ldots,\ol{n}_{p})$:
the labels satisfy the constraints listed in the text.}
\label{fig:c2} 
\end{figure}

In (\ref{eq:c.4}) we can expand
\begin{equation} \nonumber
\MM_{h_{i-1},h_{i}}^{(k_{i})}(x,n_{i},\de_{i}) 
= \!\!\!\!\!\! \sum_{T_{i}\in \gotS^{*k_{i}}_{h_{i-1},h_{i}}(n_{i},\de_{i})} \!\!\!\!\!\!
\Val_{T_{i}}(x,\de_{i}) , \qquad \qquad i=1,\ldots,p ,
\end{equation}
where we have set
\begin{equation}\label{eq:c.6}
\gotS^{*k}_{u,e}(n,\de):=\begin{cases}
\displaystyle{\bigcup_{m\ge-1} \gotS^{*k}_{m,u,e}} &\qquad \de=\calL,\partial,\partial^{2},\\
\displaystyle{\bigcup_{m>n}\gotS^{k*}_{m,u,e}}, & \qquad \de=\RR,
\end{cases}
\end{equation}
for all $k\ge 0$, $n\ge 0$ and $u,e\in\{\be,B\}$, and defined
\begin{equation}\label{eq:c.7}
\Val_{T}(x,\de):=\begin{cases}
\Val_{T}(0), & \qquad \de =\calL , \\
x\,\partial\Val_{T}(0), & \qquad \de=\partial,\\
-\Val_{T}(x), & \qquad \de = \RR, \\
\displaystyle{ x^{2}\int_{0}^{1}{\rm d}\tau(1-\tau)\partial^{2}_{x}\Val_{T}(\tau x)},
& \qquad \de = \partial^{2} ,
\end{cases}
\end{equation}
with $\Val_{T}(x)$ defined as for self-energy clusters in Section \ref{sec:6}.
Denote by $\gotC^{*}(k;h,h';\ol{n}_{0},\ldots,\ol{n}_{p})$
the set of fake clusters $\{T_{1},\ldots,T_{p}\}$ with $T_{i}\in\gotS^{*k_{i}}_{h_{i-1},h_{i}}(n_{i},\de_{i})$
for any choice of the labels $\{k_{i},n_{i},\de_{i}\}_{i=1}^{p}$ and
$\{h_{i}\}_{i=0}^{p}$ with the following constraints (see Figure \ref{fig:c2}):\\
(i) $k_{1}+\ldots+k_{p}=k$,\\
(ii) $n_{i} < \min\{\ol{n}_{i-1},\ol{n}_{i}\}$ for $i=1,\ldots,p$,\\
(iii) $h_{0}=h$, $h_{p}=h'$, \\
(iv) if $\de_{i}=\calL$ for $i=2,\ldots,p-1$, then $\de_{i-1},\de_{i+1}\ne \calL$, \\
(v) if $\de_{i}=\partial$ for $i=1,\ldots,p$, then $h_{i-1}=h_{i}$,\\
(vi) if $\de_{j}=\de_{j+1}=\ldots=\de_{j+q}=\partial$ for some
$j\in\{2,\ldots,p-1\}$ and some $q\in\{0,\ldots,p-1-j\}$ and $\de_{j-1}=\calL$,
then $\de_{j+q+1}\ne\calL$.

We call \emph{$*$-chain} any set
$C\in\gotC^{*}(k;h,h';\ol{n}_{0},\ldots,\ol{n}_{p})$ and $*$-\emph{links} the
sets $T_{1},\ldots,T_{p}$. 
By the discussion between (\ref{eq:c.4}) and (\ref{eq:c.6}), we can write
(\ref{eq:c.4}) -- and hence the sum in (\ref{eq:6.7}) -- as
\begin{equation}\label{eq:c.8}
\sum_{C \in \gotC^{*}(k;h,h';\ol{n}_{0},\ldots,\ol{n}_{p})} \Val_{C}(x),
\end{equation}
where
\begin{equation}\label{eq:c.9}
\Val_{C}(x):=\prod_{i=1}^{p}\Val_{T_{i}}(x,\de_{i}).
\end{equation}
With each $*$-chain $C$ summed over in (\ref{eq:c.8}) we associate a
\emph{depth} label $D(C)=0$; if $C=\{T_{1},\ldots,T_{p}\}$ we associate
with each $T_{i}$ the same depth label as $C$, i.e.
$D(T_{i})=D(C)=0$ for $i=1,\ldots,p$; the introduction of such a label
is due to the fact that we are performing an iterative construction
and we want to keep track of the iteration step by means of the depth label.

Given a $*$-chain $C=\{T_{1},\ldots,T_{p}\}$, for all $i=1,\ldots,p$ and
all $\ell\in L(T_{i})$
there exist $q\ge1$ \resonances $T_{i}=T^{(0)}_{i}\supset T^{(1)}_{i}
\supset\ldots\supset T^{(q-1)}_{i}$, with $T^{(j)}_{i}$ a maximal
\resonance inside $T^{(j-1)}_{i}$ for all $j=0,\ldots,q-1$ and
$T^{(q-1)}_{i}$ is the minimal \resonance containing $\ell$.
Note that both $q$ and the \resonances $T^{(1)},\ldots,
T^{(q)}$ depend on $\ell$, even though we are not making explicit such a
dependence. We call $\{T^{(j)}_{i}\}_{j=0}^{q}$ the \emph{cloud} of $\ell$
and $\{T^{(j)}_{i}\}_{j=1}^{q}$ the \emph{internal cloud} of $\ell$.
Of course if $q=0$ the internal cloud of $\ell$ is the empty set.

In (\ref{eq:c.9}) consider first a factor  $\Val_{T_{i}}(x,\de_{i})$
with $\de_{i}=\calL,\RR$. Assign a label $\gotd_{T}=1$ with
each maximal \resonance $T$ contained inside $T_{i}$.
Denote by $\gotC_{0}(T_{i})$ the set of maximal subchains $C'$ contained inside $T_{i}$.
For each $C_{j}\in \gotC_{0}(T_{i})$
there are labels $k^{(i)}_{j},h^{(i)}_{j},{h^{(i)\prime}_{j}},
\ol{n}^{(i)}_{j,0},\ldots,\ol{n}^{(i)}_{j,p_{j}}$ such
that $C_{j} \in\gotC_{1}(k^{(i)}_{j};h^{(i)}_{j},{h^{(i)\prime}_{j}};
\ol{n}^{(i)}_{j,0},\ldots,\ol{n}^{(i)}_{j,p_{j}})$.
Call $\To_{i}$ the set of nodes and lines obtained from $T_{i}$ by
removing all nodes and lines of the subchains in $\gotC_{0}(T_{i})$
and $\gotF(T_{i})$ the family of all possible sets $T_{i}'\in
\gotS^{*k_{i}}_{h_{i-1},h_{i}}(n_{i},\de_{i})$ obtained from $T_{i}$ by
replacing, for all $j=1,\ldots,|\gotC_{0}(T_{i})|$, each subchain
$C_{j}$ with any subchain $C'_{j} \in
\gotC_{1}(k^{(i)}_{j};h^{(i)}_{j},{h^{(i)\prime}_{j}};\ol{n}^{(i)}_{j,0},\ldots,
\ol{n}^{(i)}_{j,p_{j}})$. 
Note that $\To'_{i}=\To_{i}$ for all $T'_{i}\in\gotF(T_{i})$. 
If we sum together all contributions in $\gotF(T_{i})$
we obtain
\begin{equation}\label{eq:c.10}
\sum_{T'_{i}\in\gotF(T_{i})}\Val_{T'}(x,\de_{i}) = a(\de_{i}) \!\!\!\!\!\!\!\!\!\!\!\!
\sum_{\substack{C'_{j}\in\gotC_{1}(k^{(i)}_{j};h^{(i)}_{j},{h^{(i)\prime}_{j}};
\ol{n}^{(i)}_{j,0},\ldots,\ol{n}^{(i)}_{j,p_{j}}) \\ 1\le j\le |\gotC_{0}(T_{i})|}}
\!\!\!\!\!\!\!\!\!\!\!\!\!\!\!
\Val_{\To_{i}}(x,\de_{i}) \prod_{j=1}^{|\gotC_{0}(T_{i})|} \Val_{C'_{j}}(x_{\ell_{C'},\de_{i}}),
\end{equation}
where
\begin{equation}\label{eq:c.11}
\Val_{\To_{i}}(x,\de_{i}):=
\Biggl(\prod_{v\in N(\To_{i})}\calF_{v}\Biggr)
\Biggl(\prod_{\ell\in L(\To_{i})}\calG_{n_{\ell}}(x_{\ell,\de_{i}})\Biggr)
\end{equation}
and
\begin{equation}\label{eq:c.12}
a(\de)=\begin{cases}
1 , & \quad \de=\calL, \\
- 1 , & \quad \de=\RR,
\end{cases}
\qquad\qquad
x_{\ell,\de}=\begin{cases}
\oo\cdot\nn_{\ell}^{0}, & \quad \de=\calL, \\
\oo\cdot\nn_{\ell}, & \quad \de=\RR,
\end{cases}
\end{equation}
for all $\ell\in L(T')$ with $T'\in\gotF(T_{i})$.

Now consider a set $T_{i}$ in (\ref{eq:c.9}) with $\de_{i}=\partial,\partial^{2}$. Write
\begin{equation}\label{eq:c.13}
\Val_{T_{i}}(x,\partial)=x\sum_{\ell\in L(T_{i})}\partial_{x}
\calG_{\ell}(x_{\ell}^{0})
\Biggl(\prod_{v\in N(T_{i})}\calF_{v}\Biggr)
\Biggl(\prod_{\substack{\ell'\in L(T_{i}) \\ \ell'\ne\ell}}\calG_{\ell'}(x_{\ell'}^{0}) \Biggr) ,
\end{equation}
with $x_{\ell}^{0}:=\oo\cdot\nn_{\ell}^{0}$, and
$\Val_{T_{i}}(x,\partial^{2})$ as in the last line of (\ref{eq:c.7}), with
\begin{equation} \label{eq:c.14}
\begin{aligned} 
\partial^{2}_{x}\Val_{T_{i}}(\tau x)& = \!\!\!\!\!\!\!
\sum_{\ell_{1} \neq \ell_{2}\in L(T)} \left( \partial_{x} \calG_{n_{\ell_{1}}}
(x_{\ell_{1}}(\tau)) \right)
\left( \partial_{x} \calG_{n_{\ell_{2}}}(x_{\ell_{2}}(\tau)) \right)
\Biggl( \prod_{\ell\in L(T) \setminus\{\ell_{1},\ell_{2}\}} \!\!\!\!\!\!\!\!\!\!\!
\calG_{n_{\ell}}(x_{\ell}(\tau)) \Biggr)
\Biggl( \prod_{v\in N(T)} \!\!\!\! \calF_{v} \Biggr) \\
& + \sum_{\ell_{1} \in L(T)} \left( \partial_{x}^{2} \calG_{n_{\ell_{1}}}
(x_{\ell_{1}}(\tau)) \right)
\Biggl( \prod_{\ell\in L(T) \setminus\{\ell_{1}\}} \!\!\!\!\!\!\!\!
\calG_{n_{\ell}}(x_{\ell}(\tau)) \Biggr)
\Biggl( \prod_{v\in N(T)} \!\!\!\! \calF_{v} \Biggr) ,
\end{aligned}
\end{equation}
where $x_{\ell}(\tau)=x_{\ell}^{0}+ \tau x$ if
$\ell\in\calP_{T}$ and $x_{\ell}=x_{\ell}^{0}$ otherwise.
To simplify the notations we associate with each line $\ell$ a label
$d_{\ell}=0,1,2$, which denotes the number of derivatives acting
on the corresponding propagator, and set
\begin{equation} \label{eq:c.15}
\ol{\calG}_{\ell}(x) = \partial_{x}^{d_{\ell}} \calG_{n_{\ell}}(x) ;
\end{equation}
then we rewrite
\begin{subequations}
\begin{align}
&\Val_{T_{i}}(x,\partial) =x \!\!\!
\sum_{\ell_{1}\in L(T_{i})} \Biggr( \prod_{\ell\in L(T_{i})} \!\!
\ol{\calG}_{\ell}(x_{\ell}^{0}) \Biggl)
\Biggr( \prod_{v\in N(T_{i})} \!\! \calF_{v} \Biggr),
\label{eq:c.16a} \\
&\partial_{x}^{2} \Val_{T_{i}}(\tau x) = \!\!\!\!\!\!\!
\sum_{\ell_{1} , \ell_{2}\in L(T_{i})} \Biggr( \prod_{\ell\in L(T_{i})} \!\!
\ol{\calG}_{\ell}(x_{\ell}(\tau)) \Biggl)
\Biggr( \prod_{v\in N(T_{i})} \!\! \calF_{v} \Biggr) ,
\label{eq:c.16b}
\end{align}
\label{eq:c.16}
\end{subequations}
\vskip-.3truecm
\noindent  
with the constraint $\sum_{\ell\in L(T)} d_{\ell}=1,2$ for $\de_{i}=\partial,
\partial^{2}$ respectively.

Each summand in (\ref{eq:c.16a}) can be regarded as the value of a fake
cluster $T_{i}$ in which the propagator of a line $\ell_{1}$ has been
differentiated. Given a \resonance $T$ contained inside $T_{i}$,
we set $\gotd_{T}=0$ if either (1) $T$ belongs to the cloud of $\ell_{1}$
or (2) $T$ is a \resonance with either
(2a) $\ell_{T}$ non-resonant and $\ell'_{T}=\ell_{T'}$ for some $T'$ belonging
to the cloud of $\ell_{1}$ or
(2b) $\ell'_{T}$ is non-resonant and $\ell_{T}=\ell'_{T''}$ for some
$T''$ belonging to the cloud of $\ell_{1}$.
Moreover we set $\gotd_{T}=1$ if $T$ is any other maximal \resonance in
$T_{i}$ or in $T'$ with $\gotd_{T'}=0$.
If $\gotC_{0}(T_{i})=\{C_{1},C_{2},\ldots\}$,  each $C_{j}$ belongs to
$\gotC_{1}(k^{(i)}_{j};h^{(i)}_{j},
h^{(i)\prime}_{j};\ol{n}^{(i)}_{j,0},\ldots,\ol{n}^{(i)}_{j,p_{j}})$ for suitable values of the labels.
Call $\To_{i}$ the set of nodes and lines obtained from $T_{i}$ by
removing all nodes and lines of the subchains in $\gotC_{0}(T_{i})$
and $\gotF(T_{i})$ the family of all possible sets $T_{i}' \in
\gotS^{*k_{i}}_{h_{i-1},h_{i}}(n_{i},\de_{i})$ obtained from $T_{i}$ by
replacing, for all $j=1,\ldots,|\gotC_{0}(T_{i})|$, each subchain
$C_{j}$ with any subchain $C'_{j} \in
\gotC_{1}(k^{(i)}_{j};h^{(i)}_{j},{h^{(i)\prime}_{j}};\ol{n}^{(i)}_{j,0},\ldots,
\ol{n}^{(i)}_{j,p_{j}})$.
Note that $\To'_{i}=\To_{i}$ for all $T'_{i}\in\gotF(T_{i})$. 
Note also that both $\gotF(T_{i})$ and
$\To_{i}$ depend on the choice of the line $\ell_{1}$, although we are
not writing explicitly such a dependence.
By summing together all contributions obtained by choosing first the line
$\ell_{1}$ and hence the fake clusters belonging to the corresponding family
$\gotF(T_{i})$, we obtain for $\de_{i}=\partial$
\begin{equation}\label{eq:c.17}
\begin{aligned}
x\sum_{\ell_{1}\in L(T_{i})}&\sum_{T'_{i}\in\gotF(T_{i})}
\Biggr( \prod_{\ell\in L(T'_{i})} \!\!
\ol{\calG}_{\ell}(x_{\ell}^{0}) \Biggl)
\Biggr( \prod_{v\in N(T'_{i})} \!\! \calF_{v} \Biggr) \\
&= \sum_{\ell_{1}\in L(T_{i})} x \!\!\!\!\!
\sum_{\substack{C'_{j}\in\gotC_{1}(k^{(i)}_{j};h^{(i)}_{j},{h^{(i)\prime}_{j}};
\ol{n}^{(i)}_{j,0},\ldots,\ol{n}^{(i)}_{j,p_{j}}) \\ 1\le j\le |\gotC_{0}(T_{i})|}}
\!\!\!\!\!\!\!\!\!\! \Val_{\To_{i}}(x,\partial) 
\prod_{j=1}^{|{\gotC}_{0}(T_{i})|} \Val_{{C'}_{j}} (x^{0}_{\ell_{C_{j}}}) ,
\end{aligned}
\end{equation}
with
\begin{equation}\label{eq:c.18}
\Val_{\To_{i}}(x,\partial)= \Biggr( \prod_{\ell\in L(\To_{i})} \!\!
\ol{\calG}_{\ell}(x_{\ell}^{0}) \Biggl)
\Biggr( \prod_{v\in N(\To_{i})} \!\! \calF_{v} \Biggr) .
\end{equation}

We deal in a similar way with (\ref{eq:c.16b}). Indeed, each summand can
be regarded as the value of a fake cluster $T_{i}$ in which the propagators of
two lines $\ell_{1},\ell_{2}$ (possibly coinciding) have been differentiated.
As in the previous case we associate a label $\gotd_{T}=0$ with (1) the
\resonances $T$ of the clouds of both $\ell_{1},\ell_{2}$ and (2) the
\resonances $T$ such that either
(2a) $\ell_{T}$ is non-resonant and $\ell'_{T}=\ell_{T'}$
for some $T'$ belonging to the cloud of $\ell_{1}$ or $\ell_{2}$ or
(2b) $\ell'_{T}$ is non-resonant and $\ell_{T}=\ell'_{T'}$ for some
$T'$ belonging to the cloud of $\ell_{1}$ or $\ell_{2}$ or
(2c) both $\ell_{T},\ell'_{T}$ are resonant and $\ell_{T}=\ell'_{T'}$,
$\ell'_{T}=\ell_{T''}$ for some $T',T''$ belonging to the clouds of $\ell_{1},\ell_{2}$.
Moreover we associate a label $\gotd_{T}=1$ with the maximal \resonances $T$ in $T_{i}$ or in any
$T'$ with $\gotd_{T'}=0$.
Then we reason as in the previous case, by defining 
$\To_{i}$ as the set of nodes and lines obtained from $T_{i}$ by
removing all nodes and lines of the subchains in $\gotC_{0}(T_{i})$
and $\gotF(T_{i})$ the family of all possible sets $T_{i}' \in
\gotS^{*k_{i}}_{h_{i-1},h_{i}}(n_{i},\de_{i})$ obtained from $T_{i}$ by
replacing, for all $j=1,\ldots,|\gotC_{0}(T_{i})|$, each subchain $C_{j}$ with any subchain
$C'_{j} \in \gotC_{1}(k^{(i)}_{j};h^{(i)}_{j},{h^{(i)\prime}_{j}};\ol{n}^{(i)}_{j,0},\ldots,\ol{n}^{(i)}_{j,p_{j}})$.
Then for $\de_{i}=\partial^{2}$ we obtain
\begin{equation}\label{eq:c.19}
\begin{aligned}
&x^{2}\int_{0}^{1}{\rm d}\tau(1-\tau)
\sum_{\ell_{1}\ell_{2}\in L(T_{i})}\sum_{T'_{i}\in\gotF(T_{i})}
\Biggr( \prod_{\ell\in L(T'_{i})} \!\!
\ol{\calG}_{\ell}(x_{\ell}(\tau)) \Biggl)
\Biggr( \prod_{v\in N(T'_{i})} \!\! \calF_{v} \Biggr)\\
&=\sum_{\ell_{1}\ell_{2}\in L(T_{i})} x^{2}
\int_{0}^{1}{\rm d}\tau(1-\tau) \!\!\!\!\!\!\!\!\!\! \!\!\!\!\!
\sum_{\substack{C'_{j}\in\gotC_{1}(k^{(i)}_{j};h^{(i)}_{j},{h^{(i)\prime}_{j}};
\ol{n}^{(i)}_{j,0},\ldots,\ol{n}^{(i)}_{j,p_{j}}) \\ 1\le j\le|\gotC_{0}(T_{i})|}}
\!\!\!\!\! \!\!\!\!\!
\Val_{\To_{i}}(x,\partial^{2})
\prod_{j=1}^{|{\gotC}_{0}(T_{i})|}\Val_{{C'}_{j}} (x_{\ell_{{C}_{j}}}(\tau)),
\end{aligned}
\end{equation}
with
\begin{equation}\label{eq:c.20}
\Val_{\To_{i}}(x,\partial^{2})=
\Biggr( \prod_{\ell\in L(\To_{i})} \!\!
\ol{\calG}_{\ell}(x_{\ell}(\tau)) \Biggl)
\Biggr( \prod_{v\in N(\To_{i})} \!\! \calF_{v} \Biggr),
\end{equation}
where $x_{\ell}(\tau)=x_{\ell}^{0}+\tau x$ if $\ell\in\calP_{T_{i}}$
and $x_{\ell}(\tau)=x_{\ell}^{0}$ otherwise.

By summarising we obtained
\begin{equation}\label{eq:c.21}
\begin{aligned}
& \sum_{C\in\gotC^{*}(k;h,h';\ol{n}_{0},\ldots,\ol{n}_{p})}\Val_{C}(x)=
\sum_{\{T_{1},\ldots,T_{p}\}\in\gotC^{*}(k;h,h';\ol{n}_{0},\ldots,\ol{n}_{p})}
\prod_{i=1}^{p}\Val_{T_{i}}(x,\de_{i}) \\
& \qquad = \!\!\!\!\!\!\!\!\!\sum_{\{T_{1},\ldots,T_{p}\}\in\gotC^{*}(k;h,h';\ol{n}_{0},\ldots,\ol{n}_{p})}
\Biggl( \prod_{\substack{i=1 \\ \de_{i}=\partial}}^{p} \sum_{\ell_{i,1}\in L(T_{i})} \Biggr)
\Biggl( \prod_{\substack{i=1 \\ \de_{i}=\partial^{2}}}^{p}
\sum_{\ell_{i,1},\ell_{i,2}\in L(T_{i})} \Biggr) \times \\
& \qquad \qquad \qquad \times\Biggl( \prod_{i=1}^{p} a(x_{\ell_{T_{i}}},\de_{i}) \Biggr)
\Biggl( \prod_{\substack{i=1 \\ \de_{i}=\partial^2}}^{p} 
\int_{0}^{1}{\rm d}\tau_{i} \left( 1-\tau_{i} \right) \Biggr) \times \\
& \qquad \qquad \qquad \times
\prod_{i=1}^{p} \frac{1}{|\gotF(T_{i})|}
\Val_{\To_{i}}(x_{\ell_{T_{i}}},\de_{i}) \!\!\!\!\!\!\!\!\!\!\!\!\!\!\!\!
\sum_{C_{j}\in\gotC_{1}(k^{(i)}_{j};h^{(i)}_{j},{h^{(i)\prime}_{j}};
\ol{n}^{(i)}_{j,0},\ldots,\ol{n}^{(i)}_{j,p_{j}})}\!\!\!\!\!
\!\!\!\!\!\!\!\!\!\!\!\!\!\!\!\!
\Val_{C_{j}}(x_{\ell_{C_{j}}}(\tau_{i}(\de_{i}))),
\end{aligned}
\end{equation}
with $x_{\ell_{T_{i}}}=x$ by construction and
\begin{equation}\label{eq:c.22}
\Val_{\To_{i}}(x,\de_{i}) =
\Biggr( \prod_{\ell\in L(\To_{i})} \!\!
\ol{\calG}_{\ell}(x_{\ell}(\tau(\de_{i}))) \Biggl)
\Biggr( \prod_{v\in N(\To_{i})} \!\! \calF_{v} \Biggr),
\end{equation}
where we have defined $x_{\ell}(\tau):=x_{\ell}^{0}+\tau x$ and
\begin{equation}\label{eq:c.23}
\tau_{i}(\de):=\begin{cases}
0, & \quad \de=\calL,\\
1, & \quad \de=\RR,\\
0, & \quad \de=\partial,\\
\tau_{i} , & \quad \de=\partial^{2},
\end{cases}
\qquad\qquad
a(x,\de):=\begin{cases}
1 , & \quad \de=\calL , \\
-1, & \quad \de=\RR , \\
x, & \quad \de=\partial , \\
x^{2}, & \quad \de=\partial^{2}.
\end{cases}
\end{equation}
The factors $1/|\gotF(T_{i})|$ have been introduced in (\ref{eq:c.21}) to avoid overcountings.
The last sums in (\ref{eq:c.21}) have the same form as the sum (\ref{eq:c.1}), 
so that we can iterate the procedure, by writing
\begin{equation}\label{eq:c.24}
\sum_{C_{j}\in\gotC_{1}(k^{(i)}_{j};h^{(i)}_{j},{h^{(i)\prime}_{j}};
\ol{n}^{(i)}_{j,0},\ldots,\ol{n}^{(i)}_{j,p_{j}})}\!\!\!\!\!
\!\!\!\!\!\!\!\!\!\!\!\!\!\!\!\!
\Val_{C_{j}}(x_{\ell_{C_{j}}}(\tau_{i}(\de_{i}))) =
\sum_{C_{j}\in\gotC^{*}(k^{(i)}_{j};h^{(i)}_{j},{h^{(i)\prime}_{j}};
\ol{n}^{(i)}_{j,0},\ldots,\ol{n}^{(i)}_{j,p_{j}})}\!\!\!\!\!
\!\!\!\!\!\!\!\!\!\!\!\!\!\!\!\!
\Val_{C_{j}}(x_{\ell_{C_{j}}}(\tau_{i}(\de_{i})))
\end{equation}
for $i=1,\ldots,p$ and $j=1,\ldots,|\gotC_{0}(T_{i})|$, with $\Val_{T}(x)$ defined as in (\ref{eq:c.9}).
Now we associate with each $*$-chain $C_{j}$ a depth label $D(C_{j})=1$,
and if $C_{j}=\{T^{(j)}_{1},\ldots,T^{(j)}_{p_{j}}\}$ we associate with each
$T^{(j)}_{i}$ the same depth label as $C_{j}$, i.e. $D(T^{(j)}_{i})=D(C_{j})$.
More generally, by pursuing the construction, in order to keep track of the
iteration step, we associate a depth label $D(C)=d$ with each $*$-chain $C$
which appears at the $d$-th step and the same depth label $D(T)=d$ with each
$*$-link $T$ of $C$.
Since at each step the order of the chains is decreased, sooner or later the
procedure stops.

To make the notation more uniform, for any $*$-link $T$ such that
$\gotC_{0}(T)=\emptyset$ we write $\To=T$.
Given any $*$-link $T$ with $\To\ne T$, if two lines $\ell,\ell'\in L(\To)$ 
are such that there exists a maximal link $T'$ in $T$ with $\ell_{T'}=\ell$ and $\ell'_{T'}=\ell'$,
we say that the two lines are \emph{consecutive} and we write $\ell'\prec\ell$.

At the end of the procedure described above we obtain a sum of terms
of the form
\begin{equation}\label{eq:c.25}
\Biggl(\prod_{\substack{i\in I\\ \de_{T_{i}}=\partial^{2}}}
\int_{0}^{1}{\rm d}\tau_{i}(1-\tau_{i})\Biggr)
\prod_{i\in I}a(x_{\ell_{T_{i}}}(\underline{\tau}),\de_{i})
\Val_{\To_{i}}(x_{\ell_{T_{i}}}(\underline{\tau}),\de_{i}),
\end{equation}
where the following notations have been used:\\
(i) $I=\{1,2,\ldots,N\}$ for some $N\in\NNN$;\\
(ii) $\{T_{i}\}_{i\in I}$ 
are $*$-links such that (1) for all $T_{i}$ with
$D(T_{i})=d$, $d\ge1$, there exists a $*$-link $T_{j}$, $j\in I$, with
$D(T_{j})=d-1$ and two consecutive lines $\ell'\prec\ell\in L(\To_{j})$ with
the same labels as $\ell'_{T_{i}},\ell_{T_{i}}$, respectively, and conversely
(2) for all $T_{j}$ with $D(T_{j})=d$, $d\ge0$,
and all pairs of consecutive lines $\ell'\prec\ell\in L(\To_{j})$,
there is a  $*$-link $T_{i}$, $i\in I$, with $D(T_{i})=d+1$, such that
$\ell'_{T_{i}},\ell_{T_{i}}$ have the same labels as $\ell',\ell$, respectively
(roughly we can imagine to `fill the holes'
of all $\To$ such that $D(T)=0$ with all $\To'$ with $D(T')=1$, then
`fill the remaining holes' with all $\To''$ with $D(T'')=2$ and so on
up to the $*$-links of maximal depth which have no `holes');\\
(iii) $\underline{\tau}=(\tau_{1}(\de_{1}),\ldots,\tau_{N}(\de_{N}))$, with
$\tau_{i}(\de_{i})$ defined as in (\ref{eq:c.23}),
and for all $\ell\in L(\To_{1})\cup\ldots\cup L(\To_{N})$ we have set
\begin{equation}\label{eq:c.26}
x_{\ell}(\underline{\tau}):=x_{\ell}^{0}+\tau_{i_{d}}(\de_{T_{i_{d}}})
\Biggl(x_{\ell_{d}}^{0}+ \tau_{i_{d-1}}(\de_{T_{i_{d-1}}})\Bigl(x_{\ell_{d-1}}^{0}+
\tau_{i_{d-2}}(\de_{T_{i_{d-2}}})(\ldots+ \tau_{i_{0}}(\de_{T_{i_{0}}})x)
\Bigr)\Biggr),
\end{equation}
where $T_{i_{d}}$ is the minimal $*$-link (with depth $d$) containing $\ell$
and $\ell_{j}$ is the line in the $*$-link $T_{i_{j-1}}$ with
depth $j-1$ corresponding to the entering line of $T_{i_{j}}$, for
$j=1,\ldots,d$.

Recall that each propagator is differentiated at most twice
and note that, for $T$ such that $\de_{T}=\partial$, there is a line
$\ell\in L(\To)$ with $d_{\ell}=1$. Then, when bounding the product of
propagators, instead of
\begin{equation}\label{eq:c.27}
\prod_{\ell\in L(T)}\frac{\g_{0}}{2}\al_{m_{n_{\ell}}}(\oo)^{-1}
\end{equation}
with $\g_{0}$ as in the proof of Lemma \ref{lem:4.3}, we obtain
the bound (\ref{eq:c.27}) times an extra factor
\begin{equation}\label{eq:c.28}
c_{1}\al_{m_{n_{\ell}}}(\oo)^{-1}|x_{\ell_{T}}|,
\end{equation}
for suitable constant $c_{1}$. Analogously, for $T$ such that
$\de_{T}=\partial^{2}$, either there are two lines $\ell_{1},\ell_{2}$ with
$d_{\ell_{1}}=d_{\ell_{2}}=1$ or one line $\ell_{1}$ with $d_{\ell_{1}}=2$;
in both cases we obtain (\ref{eq:c.27}) times an extra factor
\begin{equation}\label{eq:c.29}
c_{1}\al_{m_{n_{\ell_{1}}}}(\oo)^{-1}\al_{m_{n_{\ell_{2}}}}(\oo)^{-1}
|x_{\ell_{T}}|^{2}.
\end{equation}
On the other hand we have no gain factor coming from the $*$-links with label
$\de=\calL,\RR$, or from the \resonances $T'$ with $\gotd_{T'}=0$.
In order to deal with such lines we need some preliminary results.

Given a $*$-link $T$, define $\To$ as before and denote by $L_{R}(\To)$
the set of resonant lines in $\To$. Set\\
(i) $L_{NR}(\To):=L(\To)\setminus L_{R}(\To)$,\\
(ii) $L_{D}(\To):=\{\ell\in L_{R} (\To)\,:\,d_{\ell}>0\}$,\\
(iii) $L^{1}_{0}(\To):=\{\ell\in L_{R}(\To)\,:\, \ell=\ell_{T'} \mbox{ for some relevant self-energy cluster }
T'\subset T \mbox{ with }\gotd_{T'}=0\}$,\\
(iv) $L^{2}_{0}(\To):=\{\ell\in L_{R}(\To)\,:\,\ell=\ell'_{T'}
\mbox{ for some relevant self-energy cluster }T' \subset T\mbox{ with
}\de_{T'}=0 \}$;\\
(v) $L_{0}(\To):=L_{0}^{1}(\To) \cup L_{0}^{2}(\To)$;\\
(vi) $L^{*}_{R}(\To)=L_{D}(\To)\cup L_{0}(\To)$.\\
Of course $L_{D}(\To)=L_{0}(\To)=\emptyset$ if $\de_{T}=\calL,\RR$.
Given a $*$-link $T$ we say that $\ell\in L(\To)$ is \emph{maximal in} $T'$
if $T'$ is the minimal \resonance contained in $\To$
(with $\gotd_{T'}=0$) such that $\ell\in L(T')$.
If there is no such \resonance we say that $\ell$ is \emph{maximal in} $T$.
Given a $*$-link $T$ with $\de_{T}=\partial,\partial^{2}$, we
denote by $L_{M}(T)$ the set of lines which are maximal in $T$; for any
\resonance $T'$ with $\gotd_{T'}=0$ we denote by
$L_{M}(T')$ the set of lines which are maximal in $T'$.

\begin{lemma}\label{lem:c.1}
Let $T$ be any $*$-link with $\de_{T}=\partial,\partial^{2}$.
For all \resonances $T'$ contained in $\To$ one has
$q_{0}(T'):=|L_{M}(T')\cap L_{0}(\To)|\le4$. Moreover
$$
q_{1}(T'):=\sum_{\ell\in L_{M}(T')}d_{\ell}\le\min\{4-q_{0}(T'),2\}.
$$
The same hold if we replace $T'$ with
the $*$-link $T$.
\end{lemma}

\prova
If $\de_{T}=\partial$ there is at most one line $\ell\in L(\To)$ such that
$d_{\ell}=1$. Let $\{T_{j}\}_{j=0}^{m}$ be the cloud of $\ell$, where we denoted
$T_{0}=T$, so that $L_{0}(\To)\subseteq\{\ell_{T_{1}},\ell'_{T_{1}},
\ldots,\ell_{T_{m}},\ell'_{T_{m}}\}$. Then $T_{m}$ is the minimal \resonance
containing $\ell$ and $L_{M}(T_{j})\cap L_{0}(\To)\subseteq\{\ell_{T_{j+1}},
\ell'_{T_{j+1}}\}$, $j=0,\ldots,m-1$. Hence $q_{0}(T_{j})\le2$ and
$q_{1}(T_{j})=0$ for $j=0,\ldots,m-1$, while $q_{0}(T_{m})=0$
and $q_{1}(T_{m})= 1$.
If $\de_{T}=\partial^{2}$ and there is one line $\ell\in L(\To)$ with
$d_{\ell}=2$ one can reason as in the previous case, denoting by
$\{T_{j}\}_{j=0}^{m}$ the cloud of $\ell$ and hence obtaining
$q_{0}(T_{j})\le2$ and $q_{1}(T_{j})=0$ for $j=0,\ldots,m-1$, while
$q_{0}(T_{m})=0$ and $q_{1}(T_{m})= 2$.
If there are two lines $\ell_{1}',\ell'_{2}\in L(\To)$ with $d_{\ell_{1}'}
=d_{\ell'_{2}}=1$ one proceeds as follows. Call  $\{T_{j}^{(i)}\}_{j=0}^{m_{i}}$
the cloud of $\ell_{i}'$, $i=1,2$, with $T^{(1)}_{0}=T^{(2)}_{0}=T_{0}=T$,
and set
$$ r:=\max\{j\ge0\,:\, T^{(1)}_{j}=T^{(2)}_{j}=:T_{j}\}. $$
If $r=m_{1}=m_{2}$ then again $q_{0}(T_{j})\le2$ and $q_{1}(T_{j})=0$ for
$j=0,\ldots,r-1$, while $q_{0}(T_{r})=0$ and $q_{1}(T_{r})= 2$.
If $r=m_{1}<m_{2}$ then
$q_{0}(T_{j})\le2$ and $q_{1}(T_{j})=0$ for $j=0,\ldots,r-1$,
$q_{0}(T_{r})\le2$ and $q_{1}(T_{r})=1$, 
$q_{0}(T^{(2)}_{j})\le2$ and $q_{1}(T^{(2)}_{j})=0$ for
$j=r+1,\ldots,m_{2}$, and $q_{0}(T_{m_{2}})=0$ and $q_{1}(T_{m_{2}})= 1$.
Finally if $r<\min\{m_{1},m_{2}\}$ then
$L_{M}(T_{r})\cap L_{0}(\To)\subseteq\{\ell_{T^{(1)}_{r+1}},
\ell'_{T^{(1)}_{r+1}},\ell_{T^{(2)}_{r+1}},\ell'_{T^{(2)}_{r+1}}\}$, so that
$q_{0}(T_{j})\le2$ and $q_{1}(T_{j})=0$ for $j=0,\ldots,r-1$,
$q_{0}(T_{r})\le 4$ and $q_{1}(T_{r})=0$, while
$q_{0}(T^{(i)}_{j})\le2$ and $q_{1}(T^{(i)}_{j})=0$ for
$j=r+1,\ldots,m_{i}$, and $q_{0}(T_{m_{i}})=0$ and $q_{1}(T_{m_{i}})= 1$,
$i=1,2$.\EP

Define the multiplicity (function) of a non-injective map as the cardinality
of its pre-image sets \cite{Fr,R}.

\begin{lemma}\label{lem:c.2} 
Let $T$ be a $*$-link with $\de_{T}=\partial,\partial^{2}$. There exists an
application $\La:L^{*}_{R}(\To)\to
L_{NR}(\To)$ with multiplicity at most 2 such that $\ze_{\ell}=\ze_{\La(\ell)}$.
\end{lemma}

\prova
By Lemma \ref{lem:c.1} there are at most four lines
$\ell_{1},\ell_{2},\ell_{3},\ell_{4}\in L^{*}_{R}(\To)$ such that, if
$T'_{i}$ denote the minimal \resonance containing $\ell_{i}$, then
$T'_{1}=\ldots=T'_{4}$. Moreover by Remark \ref{rmk:6.8} if $\ell$
is a resonant line, then the minimal \resonance containing $\ell$
contains also two non-resonant lines $\ell'_{1},\ell'_{2}$ with
the same minimum scale as $\ell$.
Therefore the assertion follows.\EP

Now consider a $*$-link $T$ contributing to (\ref{eq:c.25}) with largest
depth, say $D$. Since $T$ does not contain any resonant line, by
Lemma \ref{lem:6.6b} and Remark \ref{rmk:6.6c} we have
\begin{equation}\label{eq:c.30}
|\Val_{T_{i}}(x_{\ell_{T_{i}}}(\underline{\tau}),\RR)|\le c_{2}^{k(T_{i})}e^{-\xi
K(T_{i})/2}\le c_{3} c_{2}^{k(T_{i})} |x_{\ell_{T_{i}}}(\underline{\tau})|^{2},
\end{equation}
for some positive constants $c_{2}$ and $c_{3}$; we have also used that
$K(T_{i})\ge 2^{m_{n_{\ell_{T_i}}}-1}$ for $\de_{T_{i}}=\RR$ and
$|x_{\ell_{T_{i}}}(\underline{\tau})|\ge \al_{m_{n_{\ell_{T_i}}}}(\oo)$
if $\Psi_{n_{\ell_{T_i}}}(x_{\ell_{T_{i}}}(\underline{\tau}))\neq0$.
Therefore we can bound
\begin{equation}\label{eq:c.31}
|a(x_{\ell_{T}}(\underline{\tau}),\de_{T})||\Val_{T}(x_{\ell_{T}}(\underline{\tau}),\de_{T})|\le
\begin{cases}
c_{4}^{k(T)} , & \qquad \de=\calL,\\
c_{4}^{k(T)}|x_{\ell_{T}}(\underline{\tau})| , & \qquad \de_{T}=\partial,\\
c_{4}^{k(T)}|x_{\ell_{T}}(\underline{\tau})|^{2} , & \qquad \de=\partial^{2},\RR ,
\end{cases}
\end{equation}
for some constant $c_{4}$.
Now consider a $*$-link $T$ contributing to (\ref{eq:c.25}) with depth $D-1$.
For each resonant line $\ell\in L_{R}(\To)$ denote by $T'$ the minimal
\resonance containing $\ell$ (set $T'=T$ if there is no minimal \resonance
containing $\ell$). For all resonant lines $\ell'\in L_{M}(T')\setminus
L_{0}(\To)$, there is a subchain $C\in\gotC_{0}(T)$ for which
$\ell'$ is an internal chain-line, such that $C$ is uniquely associated (see the
comments after (\ref{eq:c.25})) with a $*$-chain $C^{*}=\{T_{i_{1}},\ldots,
T_{i_{p(C)}}\}$ with depth $D$ contributing to (\ref{eq:c.25}),
whose value can be bounded by
\begin{equation}\label{eq:c.32}
\Biggl|\prod_{j=1}^{p(C)}a(x_{\ell_{T_{j_{i}}}}(\underline{\tau}),\de_{i})
\Val_{\To_{i}}(x_{\ell_{T_{j_{i}}}}(\underline{\tau}),\de_{j_{i}})\Biggr|
\le |x_{\ell_{T_{j_{i}}}}(\underline{\tau})|^{p(C)-1}c_{5}^{k(C)},
\end{equation}
for some $c_{5}\ge 0$, and this can be obtained as follows.
If either $\de_{T_{i_{j}}}\ne \calL$ for all $j=1,\ldots,p(C)$ or
there is only one $j=1,\ldots,p(C)$ such that $\de_{T_{i_{j}}}=\calL$,
then (\ref{eq:c.32}) trivially follows from (\ref{eq:c.31}).
Otherwise if $1\le j<j'\le p(C)$ are such that $\de_{T_{i_{j}}}=
\de_{T_{i_{j'}}}=\calL$ then there is at least one $j''=j+1,\ldots,j'-1$
such that $\de_{T_{i_{j''}}}=\partial^{2},\RR$; recall the constraints
(i)--(vi) after (\ref{eq:c.7}). Then (\ref{eq:c.32}) follows.

Moreover by Lemma \ref{lem:c.1},
if there are $R_{n}^{\bullet}(T')$ resonant lines in $L_{M}(T')$ with minimum
scale $n$, we have an overall gain
$\sim \al_{m_{n}}(\oo)^{R_{n}^{\bullet}(T')-q_{0}(T')}$
and on the other hand the product of the propagators of such
resonant lines can be bounded proportionally to
$\al_{m_{n}}(\oo)^{-(R_{n}^{\bullet}(T')+q_{1}(T'))}$, with
$q_{0}(T')+q_{1}(T')\le4$.
Therefore, by Lemma \ref{lem:c.2},
if we replace the bound for the propagators of each non-resonant line
$\ell\in L_{M}(T')$, $\ze_{\ell}=n$, with $c_{6}\al_{m_{n}}(\oo)^{-3}$
for some positive constant $c_{6}$, we have
exactly a gain factor which is enough to compensate each propagator of
the resonant lines with minimum scale $n$ in $L_{M}(T')$.
But since we can reason in the same way for all $n$ and all resonant
lines in $T$, if we replace the bound for the propagators of each
$\ell\in L_{NR}(\To)$ with $c_{6}\al_{m_{n_{\ell}}}(\oo)^{-3}$
we obtain a gain proportional to
$\al_{m_{n_{\ell'}}}(\oo)^{1+d_{\ell'}}$ for any $\ell'\in L_{R}(\To)$,
and hence we can use again Lemma \ref{lem:6.4} in order to obtain
the bound (\ref{eq:c.31}) also for the $*$-links with depth $D-1$.

Then we pass to the $*$-links with depth $D-2$ and reason in the same way
as above and so on. When we arrive to the $*$-links with depth $0$
we only have to recall that they are associated with the maximal \resonances
$T$ of $\theta$, which all have label $\gotd_{T}=1$. Hence
in (\ref{eq:c.25}) we can bound
$$ \Biggl| \prod_{i=1}^{N}a(x_{\ell_{T_{i}}}(\underline{\tau}),\de_{i})
\Val_{\To_{i}}(x_{\ell_{T_{i}}}(\underline{\tau}),\de_{i})\Biggr|
\le c_{7}^{k(T_{1})+\ldots+k(T_{N})}|x|^{p-1}=c_{7}^{k}|x|^{p-1}, $$
for some positive constant $c_{7}$.
We have still to sum over all the possible contributions
of the form (\ref{eq:c.25}).
To take into account the scale labels $n_{\ell}$, $\ell\in L(\To_{1})
\cup\ldots\cup L(\To_{N})$ simply recall that for each momentum $\nn_{\ell}$
only 2 scale labels are allowed; see Remark \ref{rmk:3.8}. To sum over the
mode labels $\nn_{v}$, $v \in N(\To_{1})\cup\ldots\cup N(\To_{N})$ we can
neglect the constraints and use  a factor $e^{-(\xi/4)|\nn_{v}|}$ for each
$v$. Moreover the component labels $h_{\ell}$ are 2. 
Finally the sum over all possible unlabelled chains with order $k$ is
bounded by a constant to the power $k$. This completes the proof of
the bound (\ref{eq:6.7}).


\end{document}